\documentclass[11pt,oneside,notitlepage]{book}
\usepackage[margin=1.1in]{geometry}
\usepackage{amsmath,amssymb,amsthm,mathtools}
\usepackage{enumitem}
\usepackage[dvipsnames]{xcolor}
\usepackage[colorlinks=true,linkcolor=blue!50!black,citecolor=green!40!black]{hyperref}
\usepackage[nameinlink]{cleveref}

\numberwithin{equation}{section}
\allowdisplaybreaks
\theoremstyle{plain}
\newtheorem{theorem}{Theorem}[section]
\newtheorem{proposition}[theorem]{Proposition}
\newtheorem{lemma}[theorem]{Lemma}
\newtheorem{corollary}[theorem]{Corollary}
\newtheorem{conjecture}[theorem]{Conjecture}
\theoremstyle{definition}
\newtheorem{definition}[theorem]{Definition}
\newtheorem{remark}[theorem]{Remark}
\crefname{theorem}{Theorem}{Theorems}
\crefname{proposition}{Proposition}{Propositions}
\crefname{lemma}{Lemma}{Lemmas}
\crefname{corollary}{Corollary}{Corollaries}
\crefname{definition}{Definition}{Definitions}
\crefname{remark}{Remark}{Remarks}
\crefname{conjecture}{Conjecture}{Conjectures}

\newcommand{\A}{\mathbb{A}}
\newcommand{\Q}{\mathbb{Q}}
\newcommand{\Z}{\mathbb{Z}}
\newcommand{\C}{\mathbb{C}}
\newcommand{\cO}{\mathcal{O}}
\newcommand{\GL}{\mathrm{GL}}
\newcommand{\SL}{\mathrm{SL}}
\newcommand{\SO}{\mathrm{SO}}
\newcommand{\Ort}{\mathrm{O}}
\newcommand{\Mp}{\mathrm{Mp}}
\newcommand{\Mat}{\mathrm{Mat}}
\newcommand{\so}{\mathfrak{so}}
\newcommand{\cS}{\mathcal{S}}
\newcommand{\cW}{\mathcal{W}}
\newcommand{\cB}{\mathcal{B}}
\newcommand{\tGL}{\widetilde{\GL}}
\newcommand{\tr}{\operatorname{tr}}
\newcommand{\Tr}{\operatorname{Tr}}
\newcommand{\N}{\operatorname{N}}
\newcommand{\val}{\operatorname{val}}
\newcommand{\diag}{\operatorname{diag}}
\newcommand{\vol}{\operatorname{vol}}

\newcommand{\one}{\mathbf{1}}
\newcommand{\hil}[2]{(#1,#2)}
\newcommand{\T}{\mathsf{T}}
\renewcommand{\doteqdot}{\mathrel{\dot=}}
\newenvironment{abstract}
{\begin{center}\small\bfseries Abstract\end{center}\begin{quote}\small}
	{\end{quote}}
\newcommand{\blfootnote}[1]{%
	\begingroup\renewcommand{\thefootnote}{}\footnote{#1}%
	\addtocounter{footnote}{-1}\endgroup}

\title{\textbf{A relative trace formula for the orthogonal Gan-Gross-Prasad model and the proof of its fundamental lemma}}
\author{Zhengyu Mao \and Chen Wan}

\date{}

\begin{document}
\frontmatter
\maketitle

\blfootnote{Zhengyu Mao: Department of Mathematics and Computer Science, Rutgers University, Newark, NJ 07102, USA. \textit{Email:} \texttt{zmao@rutgers.edu}}
\blfootnote{Chen Wan: Department of Mathematics and Computer Science, Rutgers University, Newark, NJ 07102, USA. \textit{Email:} \texttt{chen.wan@rutgers.edu}}

\begin{abstract}
In this paper, we propose a comparison of relative trace formulas to study the
Gan--Gross--Prasad conjecture. We unfold the relative trace formulas, match their regular orbits, and formulate the corresponding fundamental lemma and smooth transfer conjecture. Finally, we prove the fundamental lemma for the unit element of the Hecke algebra.
\end{abstract}

{\let\clearpage\relax\tableofcontents}

\mainmatter

\chapter{Introduction}

\section{The Gan--Gross--Prasad conjecture}

Let $k$ be a number field with ring of adeles $\A$. Let $W\subset V$ be nondegenerate quadratic spaces over $k$ with $\dim V=\dim W+1$, put $G=\SO(V)\times\SO(W)$, and let $H\simeq\SO(W)$ be embedded diagonally in $G$. For a cuspidal automorphic representation $\pi\boxtimes\sigma$ of $G(\A)$, the \emph{Gross--Prasad period} is the linear form
\[\mathcal P_H(\phi)=\int_{[H]}\phi(h)\,dh,\qquad\phi\in\pi\boxtimes\sigma .\]
Gross and Prasad \cite{GP92,GP94} conjectured that, for tempered $\pi\boxtimes\sigma$, $\mathcal P_H$ is not identically zero on some member of the Vogan packet of $\pi\boxtimes\sigma$ (that is, on $\pi'\boxtimes\sigma'$ for some pair of pure inner forms $W'\subset V'$) if and only if the central value $L(\frac12,\pi\times\sigma)$ of the Rankin--Selberg $L$-function does not vanish. They predicted moreover which member of the packet carries the period, in terms of local root numbers. Gan, Gross and Prasad \cite{GGP} later extended the conjectures to all classical groups: unitary groups (where one considers $\mathrm U(n+1)\times\mathrm U(n)$, or more generally Bessel and Fourier--Jacobi periods for $\mathrm U(n)\times\mathrm U(m)$) and symplectic-metaplectic groups. The conjecture has a refined form. Ichino and Ikeda \cite{II} proposed, for orthogonal groups in the codimension one case, an explicit formula of the shape
\[\frac{|\mathcal P_H(\phi)|^2}{\langle\phi,\phi\rangle}=\frac{C}{|S_\pi|}\,\frac{\Delta_G\,L(\frac12,\pi\times\sigma)}{L(1,\pi,\mathrm{Ad})L(1,\sigma,\mathrm{Ad})}\prod_v\alpha_v(\phi_v),\]
where $S_\pi$ is the component group of the global $L$-parameter, $\Delta_G$ is a product of special values of abelian $L$-functions, and the $\alpha_v$ are normalized integrals of matrix coefficients over $H(k_v)$. Later, N.~Harris \cite{Harris} formulated the analogous conjecture for unitary groups and Y. Liu \cite{Liu} formulated the analogous conjecture for general non-reductive Bessel case.

For even special orthogonal groups the refined formula requires a qualification on the automorphic realization and a multiplicity correction;  the correct formulation is due to H. Xue \cite[Conjecture~6.2.1]{XueFJ}. In particular, \cite[Remark~6.2.2]{XueFJ} records the corresponding correction to Liu's Bessel-period conjecture \cite[Conjecture~2.5]{Liu}: in the multiplicity-two case one must specify an outer-invariant realization and include the extra factor $2$.

\subsection*{The unitary case}

For unitary groups the conjecture and its refined form are now known for tempered representations, by a comparison of relative trace formulas proposed by Jacquet and Rallis \cite{JR}. Jacquet and Rallis compare a relative trace formula on $\mathrm U(n+1)\times\mathrm U(n)$, integrated twice over the diagonal $\mathrm U(n)$, with one on $\GL_{n+1}(E)\times\GL_n(E)$ integrated over $\GL_n(E)$ and over $\GL_{n+1}(k)\times\GL_n(k)$; the latter detects the base changes $\Pi$ of $\pi$ through the Rankin--Selberg period and the Flicker--Rallis period. They stated the corresponding fundamental lemma and smooth transfer conjecture, and proved the fundamental lemma for $n\le2$. The subsequent history is as follows.
\begin{itemize}
	\item Yun \cite{Yun} proved the Jacquet--Rallis fundamental lemma over local fields of positive characteristic with $p$ large, by a geometric argument in the spirit of Ng\^o's proof of the fundamental lemma; Gordon's appendix to \cite{Yun} transfers it to $p$-adic fields with large residual characteristic.
	\item W.~Zhang \cite{ZhangFT} proved the smooth transfer conjecture over $p$-adic fields, by a Fourier transform argument and an induction on $n$, and deduced the global conjecture for tempered $\pi$ satisfying certain local conditions (for instance, supercuspidality at a split place). In \cite{ZhangRefined} he proved the refined (Ichino--Ikeda type) formula under similar conditions.
	\item Xue \cite{Xue} proved an archimedean approximation of smooth transfer, which removes part of the local conditions at the archimedean places.
	\item Beuzart-Plessis gave a purely local proof of the fundamental lemma valid for all residual characteristics \cite{BP}, based on the Fourier transform and an induction on the dimension; he also proved the comparison of local relative characters and the Plancherel formula needed in the refined conjecture \cite{BPchar,BPplanch}.
	\item Beuzart-Plessis, Liu, Zhang and Zhu \cite{BPLZZ} removed the local conditions in the stable case (when the base change of $\pi$ is cuspidal), by isolating the cuspidal part of the spectrum with multipliers.
	\item Beuzart-Plessis, Chaudouard and Zydor \cite{BPCZ} treated the endoscopic case, using the full coarse relative trace formulas of Zydor \cite{Zydor} and the singular transfer of Chaudouard--Zydor \cite{CZ}. This completes the proof of the global conjecture and of its refined version for tempered representations of unitary groups.
\end{itemize}

\subsection*{The orthogonal case}

For special orthogonal groups much less is known. The refined conjecture holds in low rank: for $\SO(3)\times\SO(2)$ it is Waldspurger's formula for toric periods \cite{Wal85}, and for $\SO(4)\times\SO(3)$ it is Harris-Kudla's and Ichino's triple product formula (\cite{HK}, \cite{Ich}). Both are proved through the theta correspondence and seesaw identities, which are not available in general. Waldspurger's formula has also been approached by relative trace formulas \cite{Jac86}; Krishna \cite{Krishna} gave a proof of it by a relative trace formula comparison of the type considered here, built from theta functions. As far as we know, no relative trace formula for the orthogonal Gan--Gross--Prasad problem in higher rank has been carried out. 

Meanwhile, one direction (nonvanishing of the $L$-value implies nonvanishing of a period on some member of the packet) was obtained by Ginzburg--Jiang--Rallis \cite{GJR} using automorphic descent, and in much greater generality by Jiang--Zhang \cite{JZ}.

\section{The relative trace formulas}

In this paper we propose a  relative trace formulas comparison for the orthogonal Gan--Gross--Prasad conjecture in the case $\SO(2n+1)\times\SO(2n)$, carry out their unfolding, formulate the smooth transfer conjecture that the comparison requires, and prove the corresponding fundamental lemma for the unit element of the Hecke algebra, for the regular orbits.

\subsection*{The setting}
In the usual orthogonal relative trace formula, the test functions are on $(\SO(V)\times\SO(W))(\A)$, with $\dim V=2n+1$ and $\dim W=2n$; \(\SO(W)\times\SO(W)\) acts on \(\SO(V)\times\SO(W)\) by left and right multiplication through the diagonal embedding. We enlarge this two-sided acting group to
\[
\mathcal H_W=\{(h_1,h_2)\in O(W)^2:\det h_1=\det h_2\}
\]
(\Cref{def:RTF}(a)), following the rank-one construction of Krishna \cite[\S4B]{Krishna}.
 Its connected part gives the usual $\SO(W)$-periods, while the full distribution also averages the simultaneous outer action. Its geometric terms are therefore $\Ort(W)$-orbital integrals, as required by the linear matching.

The linear side reflects the structure of these $L$-functions under functorial transfer to general linear groups \cite{Arthur}. A tempered $\pi$ on $\SO(2n+1)$ transfers to an automorphic representation $\Pi$ of $\GL_{2n}$ of symplectic type, which is characterized by a nonvanishing Shalika period \cite{JS}. A tempered $\sigma$ on $\SO(2n)$ transfers to a representation $\Sigma$ of $\GL_{2n}$ of orthogonal type; these have a pole of the symmetric square $L$-function, which is detected by integrating against a product of two exceptional theta forms on the double cover of $\GL_{2n}$ \cite{BG,Takeda}. Finally $L(s,\Pi\times\Sigma)$ is represented by the Rankin--Selberg integral against a mirabolic Eisenstein series on $\GL_{2n}$ \cite{JPSS}. Accordingly, the linear relative trace formula (\Cref{def:RTF}(b)) integrates the product of the kernels of two test functions on $\GL_{2n}$ against the Shalika character in one variable, the product $\Theta^+\Theta^-$ of two exceptional theta forms in another, and the Rankin--Selberg Eisenstein series in the diagonal variable. This setting is well known to experts; the difficulty lies in making the geometric side explicit.

\subsection*{One linear side, a sum of orthogonal sides}
The two relative trace formulas are not compared one to one. The linear distribution is associated with a quadratic character \(\eta=\eta_\delta\), whereas the orthogonal distribution depends on a quadratic space \(W\). The intended comparison is therefore with
\[
\sum_{\substack{[W]\\ \dim W=2n\\ \operatorname{disc}(W)=\delta}} I_W(f_W),
\]
where \(I_W\) is the distribution defined using the extended relative integration group \(\mathcal H_W\). The sum is over isometry classes of quadratic spaces. On the regular comparison chart, a linear datum \(z\) determines a unique isometry class \(W_z\) and an \(O(W_z)\)-conjugacy class of \(g_z\). Different data can produce nonisometric spaces with the same discriminant, which explains the need for the sum. This is consistent with the inner-form aspect of the Gross–Prasad conjecture; its precise spectral interpretation must also account for the outer action and the automorphic-realization qualifications stated above.


\subsection*{Why $\SO(2n+1)\times\SO(2n)$ suffices}
As in the unitary case, it is enough to treat one parity. A cuspidal representation of a general pair $\SO(2a+1)\times\SO(2b)$ consisting of an odd and an even special orthogonal group can be parabolically induced to $\SO(2n+1)\times\SO(2n)$ for suitable $n$, and the Gross--Prasad period of the resulting Eisenstein series is expressed through the (Bessel) period of the original pair. Therefore, once the comparison of the two relative trace formulas for $\SO(2n+1)\times\SO(2n)$ is established for the full spectral expansion, including the contributions of Eisenstein series, the conjecture follows for all pairs $\SO(2a+1)\times\SO(2b)$. This is why we only formulate the relative trace formulas in this setting.

\subsection*{The unfolding}
The geometric expansion of the orthogonal side is elementary: the $\mathcal H_W$-orbits in $\SO(V)\times\SO(W)$ correspond to $\Ort(W)$-conjugacy classes in $\SO(V)$, and the regular part unfolds into orbital integrals (\Cref{prop:orth-unfold}). On the linear side the unfolding process is less straightforward. The key point of the unfolding is the uniqueness of the theta kernel with respect to the Shalika unipotent radical: the full-rank twisted Jacquet module of the exceptional representation along the unipotent radical $U$ is one dimensional \cite[Theorem 3.1]{Kaplan}. This has two consequences. First, after expanding the two theta forms along $U$ and applying partial Poisson summation to the Eisenstein series, the four characters of $U$ (two theta frequencies $A$ and $B$, the Poisson character and the Shalika character) must combine to the trivial character. This forces the rank one relation
\[A-B+vr=I_n,\]
so that the geometric side is indexed by enhanced matrix data $z=(A,v,r)$ modulo $\GL_n$ (\Cref{lem:poisson}). Second, the global Fourier coefficients of full rank factor into products of \emph{local} Shalika functionals of the exceptional representations (\Cref{lem:factor}). Combined with the observation that the Mellin variable of the Eisenstein series supplies exactly the central direction missing from the Levi quotient of the Shalika subgroup, this turns the regular part of the linear side into a sum over rational orbits of products of local orbital integrals (\Cref{thm:linear-unfold}). The local orbital integrals thus defined is dependent on orbital representative up to a sign. To facilitate the proof of fundamental lemma, we further normalize it to produce orbital integrals that are valued independent of representatives (\Cref{lem:equivariance}). The convergence of all regular contributions, uniformly in the Mellin variable so that the central value $s=\frac12$ may be substituted, is proved in \Cref{app:convergence}.

\subsection*{Main results}
The two families of regular orbits match: a regular datum $z=(A,v,r)$ determines a quadratic space $W_z$ of dimension $2n$ and an element $g_z\in\SO(W_z\oplus Fe)$, and this gives a bijection between the regular orbits on the linear side and the regular orbits on all quadratic spaces of dimension $2n$ taken together, under which $\det A\det B$ corresponds to the discriminant of $W_z$ (\Cref{prop:matching}). With the normalized orbital integrals $J_z$ and $I_z$ of \S\ref{sec:FL} we formulate the smooth transfer conjecture (\Cref{conj:transfer}). Our main theorem is the fundamental lemma for the unit elements (\Cref{thm:FL}).

\medskip\noindent\textbf{Theorem} (\Cref{thm:FL})\textbf{.} \emph{Let $F$ be a nonarchimedean local field of characteristic zero and odd residual characteristic with residue field of cardinality $q>3$, and let $\psi$ have conductor $\cO$. Then for every regular $z=(A,v,r)$,
	\[J^\circ(z)=I^\circ(z),\]
	where $J^\circ$ is the normalized linear orbital integral of the basic data $\one_{\GL_{2n}(\cO)}$, $\xi^{+,\circ}\otimes\xi^{-,\circ}$ and $\one_{\cO^{2n}}$, and $I^\circ(z)$ is the number of selfdual lattices $L\subset W_z$ with $g_z(L\oplus\cO e)=L\oplus\cO e$ (zero if $W_z$ has no selfdual lattice).}\medskip

The linear orbital integral is a signed weighted sum over lattices, the weights being values of the spherical exceptional Shalika function; the orthogonal orbital integral is an unweighted lattice count. 

\section{Outline of the proof}

The general strategy follows Beuzart-Plessis' proof of the Jacquet--Rallis fundamental lemma \cite{BP} (which is inspired by Jacquet's early works on fundamental lemmas): one realizes the difference of the two basic orbital integrals in a representation of $\SL_2(F)$ (or of its double cover), in which the upper unipotent elements $\begin{pmatrix} 1 & t\\ 0& 1\end{pmatrix}$ act by multiplication by $\psi(t\nu)$ for an invariant $\nu$ and the Weyl element acts by a Fourier transform; one shows by induction that the difference is supported on $\{\nu\in\varpi\cO\}$ and is fixed by the Weyl element, and concludes that it vanishes. The orthogonal case is more complicated in two respects.

First, the induction cannot be run within a single family. Splitting off a vector of unit norm changes the dimension of the quadratic space by one, so the induction must pass alternately through quadratic spaces of even and odd dimension. We therefore introduce two Lie algebra models: an \emph{even} model (\Cref{def:even}), with a framing vector and covector, and an \emph{odd} model (\Cref{def:odd}), without framing. In \Cref{ch:decomposition} we decompose the group orbital integrals into these two parts: according to the residual spectrum of $A$ we split the data into a residually nilpotent part and a residually invertible part. By applying different Cayley transforms to these two parts, we will prove (\Cref{thm:factorization})
\[I^\circ(z)=I^e\cdot I^o,\qquad J^\circ(z)=J^e\cdot J^o .\]
On the orthogonal side this is a statement about lattice counts. On the linear side we must check that the signs of the weights are compatible with the factorization.

Second, the Fourier transform has a different origin in the two parities. In the even model it is simply the Fourier transform in the framing variables. In the odd model there is no framing variable, and the Fourier transform has to be extracted from the exceptional representation itself: we use a Heisenberg (Jacobi) quotient of $\theta_{2n}$ along a maximal parabolic, identified by a theorem of Kaplan with a Weil representation tensored with the rank two exceptional representation, and we show that the Shalika coefficients become Schwartz functions of the last row and column on which $\Mp_2(F)$ acts by a Weil representation (\Cref{prop:jacobi}).

In \Cref{ch:induction} we state the compatibility of the Fourier transform with transfer (\Cref{thm:fourier}) and, assuming it, prove the fundamental lemma. The two unit descents (\Cref{lem:even-odd,lem:odd-even}) reduce the even identity of size $n$ to the odd identity of size $n$, and the odd identity of size $n$ to the even identity of size $n-1$, so that the two Lie algebra identities are proved by an interlaced induction in the order $E_0,O_1,E_1,O_2,E_2,\dots$ (\Cref{thm:LieFL}). Combined with \Cref{thm:factorization} this proves \Cref{thm:FL}. In \Cref{ch:fourier} we prove \Cref{thm:fourier}: after integration over the centralizer, both linear orbital integrals become twisted hyperbolic integrals and the orthogonal ones become norm integrals over a quadratic extension, and an elementary binary kernel identity (\Cref{prop:binary}) shows that the Weyl element acts on both by the same kernel. 

Notation and the background on metaplectic covers and exceptional representations are collected in \Cref{ch:notation}; the analytic estimates and several technical lemmas on the metaplectic cover are in \Cref{app:convergence,app:technical}.

\section{What remains to be done}

To deduce the Gan--Gross--Prasad and Ichino--Ikeda conjectures for general $\SO(2a+1)\times\SO(2b)$, from the relative trace formulas of this paper, the remaining steps are the orthogonal analogues of the steps carried out in the unitary case:
\begin{enumerate}[label=(\arabic*)]
	\item the smooth transfer conjecture (\Cref{conj:transfer}) over $p$-adic fields, and its archimedean analogue (compare \cite{ZhangFT,Xue,BPchar});
	\item the fundamental lemma for the full spherical Hecke algebra, at least in the form needed to separate spectra, 
	\item the comparison of the full geometric sides, including the singular orbits, the two Fourier axes and the noncyclic terms of the linear side, through a suitable regularization of both relative trace formulas (compare \cite{Zydor,CZ});
	\item the spectral expansions of both sides, retaining the outer action and automorphic multiplicity spaces, the isolation of the cuspidal spectrum (compare \cite{BPLZZ}), and the treatment of the endoscopic contributions (compare \cite{BPCZ});
	\item the factorization of the global relative characters and the comparison of the local relative characters, which identifies the local factors in the Ichino--Ikeda formula (compare \cite{BPchar,BPplanch}).
\end{enumerate}
In each of these steps we expect the methods developed for unitary groups to apply to the orthogonal case with modifications similar to those in the present paper, notably the systematic use of the exceptional representations and of the interlaced even and odd models. We do not pursue these directions here.

In another direction, J. Xiao \cite{Xiao} showed that the fundamental lemma for endoscopic transfer (established by B. C. Ng\^o\cite{Ngo}) in the case of unitary group can be derived using the Jacquet-Rallis fundamental lemma.  The same holds for the orthogonal group case, and the derivation of the fundamental lemma for endoscopic transfer will be treated in a future work.

\section{Acknowledgement} The majority of the argument was produced by GPT Astra and is clearly inspired by ideas of Beuzart-Plessis and Jacquet. Astra’s approach was also likely influenced by the authors’ prompts and by their prior discussions with the AI on other problems involving relative trace formulas. The authors are responsible for verifying and streamlining the argument, and preparing the manuscript, with additional assistance from Claude.

The second author’s work is partially supported by the NSF grant DMS-2349836 and a Simons Travel Grant.

\chapter{Notations and Backgrounds}\label{ch:notation}

This chapter fixes the notation used throughout the paper and collects the facts about the metaplectic double covers of general linear groups and their exceptional representations that we need. Local objects are defined over a local field, global objects over a number field $k$ and its ring of adeles $\A$. Everything here is standard; the only results with proofs are the elementary properties of the conjugation sign $\epsilon$ (\Cref{lem:eps}).

\section{Fields, characters, measures and matrices}\label{sec:notation}

\paragraph{Local fields.}
Throughout, $F$ denotes a nonarchimedean local field of characteristic zero and odd residual characteristic, $\cO$ its ring of integers, $\varpi$ a uniformizer and $q$ the cardinality of the residue field; $|\cdot|$ is the normalized absolute value, $|\varpi|=q^{-1}$. We write $\hil ab$ for the quadratic Hilbert symbol of $F$, and $\hil ab_K$ for that of a finite extension $K$ of $F$ (componentwise for a finite \'etale algebra $K=\prod_iK_i$). For $\delta\in F^\times$ let $\eta_\delta=\hil{\cdot}{\delta}$ be the quadratic character of $F^\times$ attached to $F(\sqrt\delta)$. When an unramified situation is considered, the additive character $\psi$ of $F$ has conductor $\cO$. Additive Haar measures are selfdual with respect to the pairing indicated in each case; the Haar measure of $\GL_j(F)$ gives $\GL_j(\cO)$ volume one.

For $a\in F^\times$ let $\psi_a(x)=\psi(ax)$, and for $a_1,\dots,a_m\in F^\times$ write $\langle a_1,\dots,a_m\rangle$ for the diagonal quadratic form $a_1x_1^2+\dots+a_mx_m^2$ on $F^m$.

\paragraph{Global fields.}
Let $k$ be a number field, $\A$ its ring of adeles, $v$ a place of $k$ and $k_v$ the completion; for finite $v$ we write $\cO_v$ and $q_v$ for the ring of integers and the cardinality of the residue field. The absolute values $|\cdot|_v$ are normalized, and $|x|_\A=\prod_v|x_v|_v$, so that $|a|_\A=1$ for $a\in k^\times$ (product formula). Let $\psi=\otimes_v\psi_v$ be a nontrivial character of $k\backslash\A$; then $\psi_v$ has conductor $\cO_v$ for almost all $v$. We write $\hil ab_v$ for the quadratic Hilbert symbol of $k_v$; Hilbert reciprocity says $\prod_v\hil ab_v=1$ for $a,b\in k^\times$. For $c\in k^\times$ put $\eta_c=\prod_v\hil{\cdot}{c}_v$; by reciprocity it is a quadratic character of $k^\times\backslash\A^\times$, every quadratic character of $k^\times\backslash\A^\times$ is of this form, and $c$ is determined modulo $k^{\times2}$. Its local components are the characters $\eta_{c}$ of $k_v^\times$ defined above.

For an algebraic group $G$ over $k$ we write $[G]=G(k)\backslash G(\A)$ and, if $Z$ is its center, $[G]_Z=Z(\A)G(k)\backslash G(\A)$. Haar measures on $G(\A)$ are products of local Haar measures; for $G=\GL_j$ these are normalized at the finite places as in the local case, so that $\GL_j(\cO_v)$ has volume one. The additive measure on $\A$ is selfdual for $\psi$, so that $k\backslash\A$ has volume one. Quotient measures are taken compatibly.

\paragraph{Test functions.}
Locally, test functions on groups are smooth and compactly supported, and $\cS(F^m)$ is the space of Schwartz functions. Globally, test functions on $G(\A)$ are finite sums of factorizable smooth compactly supported functions, and Schwartz functions on $\A^m$ are Schwartz--Bruhat functions.

\paragraph{Matrices.}
Let $R$ be one of the rings $F$, $k$, $k_v$ or $\A$. Column vectors are elements of $R^j$, row vectors are elements of $R^{1\times j}$, and $\Mat_j(R)$ is the space of $j\times j$ matrices; $N_j$ is the group of upper triangular unipotent matrices. For $C,h\in\GL_j(R)$ we put $C^h=h^{-1}Ch$. For $X\in\Mat_{j}(R)$ and $h,C\in\GL_j(R)$ we use the following elements of $\GL_{2j}(R)$:
\[
u(X)=\begin{pmatrix} I_j & X\\ 0& I_j\end{pmatrix},\qquad d_h=\begin{pmatrix} h&0\\0&h\end{pmatrix},\qquad m_C=\begin{pmatrix} C&0\\0&I_j\end{pmatrix}.
\]
Note that $m_C u(X) m_C^{-1}=u(CX)$ and $d_h^{-1}m_Cd_h=m_{C^h}$.

For $C\in\GL_j(F)$ there are $k_1,k_2\in\GL_j(\cO)$ and unique integers $a_1\ge\dots\ge a_j$ with $C=k_1\diag(\varpi^{a_1},\dots,\varpi^{a_j})k_2$ (elementary divisors); the $a_i$ are the \emph{Smith exponents} of $C$. Thus $C$ is integral if and only if $a_j\ge0$, and $\val\det C=\sum_ia_i$. A matrix $C=(c_{ik})$ is \emph{upper Hessenberg} if $c_{ik}=0$ for $i>k+1$; companion matrices are upper Hessenberg with subdiagonal entries $c_{k+1,k}=1$. For an integral upper Hessenberg matrix $C$ with subdiagonal entries $1$, the minor obtained by deleting the first row and the last column is the determinant of a unipotent upper triangular matrix, hence equals $1$; so the Smith exponents of $C$ are $(\val\det C,0,\dots,0)$.

\section{The metaplectic cover}\label{sec:cover}

\paragraph{The local cover.}
For $j\ge1$ let $\tGL_j(F)$ be the block compatible double cover of Kaplan--M\"ollers \cite{KM}. As a set it is $\GL_j(F)\times\{\pm1\}$ with multiplication given by the Banks--Levy--Sepanski (BLS) cocycle $\sigma$, and $s(g)=(g,1)$ denotes the BLS section. The cocycle is normalized, $\sigma(1,g)=\sigma(g,1)=1$. We use the following standard facts \cite[\S1, \S2.1, Claims 3.7, 3.18, 3.27]{KM}, \cite{Kable}. They are stated for $F$ as in \S\ref{sec:notation}, i.e. nonarchimedean of characteristic zero and odd residual characteristic, and the items (C3)–(C5), (C7), (C8) are used only in this case. Items \ref{C:block}, \ref{C:unip}, and \ref{C:center} remain valid over every local field of characteristic zero, including archimedean and dyadic fields; the precise justification is recorded in the paragraph `Other completions' below.
\begin{enumerate}[label=(C\arabic*)]
	\item\label{C:block} (Block compatibility.) For $g_i,h_i\in\GL_{a_i}(F)$,
	$\sigma(\diag(g_1,g_2),\diag(h_1,h_2))=\sigma(g_1,h_1)\sigma(g_2,h_2)\hil{\det g_1}{\det h_2}$ \cite[\S1.2, (1.1)]{KM}.
	\item\label{C:unip} (Unipotent elements.) If $U$ is a unipotent algebraic $F$-subgroup of $\GL_j$, then $U(F)$ has a unique splitting into $\tGL_j(F)$, whose value at $u$ is called the canonical lift and written $\tilde u$. Conjugation by any element of $\tGL_j(F)$ carries canonical lifts to canonical lifts. On the upper triangular unipotent group $N_j(F)$ the canonical lift is $s$, and $s(ugu')=s(u)s(g)s(u')$ for $u,u'\in N_j(F)$, $g\in\GL_j(F)$ \cite[\S1.2]{KM}.
	\item\label{C:torus} (Torus.) $\sigma(\diag(a_1,\dots,a_j),\diag(b_1,\dots,b_j))=\prod_{i<k}\hil{a_i}{b_k}$.
	\item\label{C:weyl} (Weyl elements.) Let $w$ be the standard signed representative of a Weyl element, obtained from a reduced expression using the simple representatives whose $2\times2$ block is $\left(\begin{smallmatrix}0&-1\\1&0\end{smallmatrix}\right)$, and let $t=\diag(t_1,\dots,t_j)$. Then $\sigma(t,w)=1$ and $\sigma(w,t)=\prod_{i<k,\,w(i)>w(k)}\hil{-t_k}{t_i}$. Consequently $s(w)s(t)s(w)^{-1}=\sigma(w,t)\, s(wtw^{-1})$.
	\item\label{C:compact} (Compact splitting.) There is a splitting $\kappa:\GL_j(\cO)\to\tGL_j(F)$. It agrees with $s$ on $N_j(\cO)$, on the diagonal matrices in $\GL_j(\cO)$ and on the Weyl representatives of \ref{C:weyl}; it is compatible with the block embeddings $g\mapsto\diag(g,I)$ and $g\mapsto\diag(I,g)$. We do not assume that $\kappa=s$ on all of $\GL_j(\cO)$.
	\item\label{C:center} (Scalars.) More generally, in $\tGL_N(F)$ one has
	\[
	s(g)s(aI_N)s(g)^{-1}s(aI_N)^{-1}=\hil{\det g}{a}^{N-1}.
	\]
	In particular, for $N=2j$ this is $\hil{\det g}{a}$. Compare \cite[(2.6)]{Kable}; the uniform argument valid at every completion is given below.
	\item\label{C:squares} (Square tori.) If $t$ is diagonal with all entries in $F^{\times2}$, then $\sigma(t,g)=\sigma(g,t)=1$ for all $g\in\GL_j(F)$ \cite[Claim 3.7, Remark 3.8]{KM}.
	\item\label{C:KP} (Kazhdan--Patterson.) Let $S\in\GL_j(F)$ be regular semisimple and $K=F[S]=\prod_iK_i$, and regard $S$ as an element of $K^\times$: writing $P_S=\prod_iP_i$ with distinct irreducible $P_i$, we have $K=F[S]\cong F[T]/(P_S)\cong\prod_iK_i$ with $K_i=F[T]/(P_i)$, and $S$ corresponds to $(S_i)_i$, where $S_i$ is the image of $T$ in $K_i$. Then for every $t\in K^\times$ (an element of the centralizer of $S$ in $\GL_j(F)$)
	\[
	\frac{\sigma(S,t)}{\sigma(t,S)}=\hil{\det S}{\det t}\prod_i\hil{S_i}{t_i}_{K_i}
	\]
	\cite[Proposition 0.1.5]{KP}. Since $S$ and $t$ commute, the left side is the commutator $s(S)s(t)s(S)^{-1}s(t)^{-1}$ of their lifts.
\end{enumerate}
By \ref{C:block}, the maps $\GL_a\to\GL_{a+b}$, $C\mapsto\diag(C,I_b)$ and $C\mapsto \diag(I_b,C)$, lift to homomorphisms of covers sending $s(C)$ to $s(\diag(C,I_b))$, resp. $s(\diag(I_b,C))$; we write them $\tilde C\mapsto \tilde C\oplus 1$ and $\tilde C\mapsto 1\oplus\tilde C$.

For $C,h\in\GL_j(F)$ define the sign
\begin{equation}\label{eq:eps}
	\epsilon(C,h)=\frac{\sigma(C,h)}{\sigma(h,C^h)}\,\hil{\det C}{\det h}\in\{\pm1\}.
\end{equation}
The determinant symbol is included so that $\epsilon$ is a classical character along centralizers: if $S$ is regular semisimple and $t\in K^\times=F[S]^\times$, then $S^t=S$, so $\epsilon(S,t)=\hil{\det S}{\det t}\,\sigma(S,t)/\sigma(t,S)$, and by \ref{C:KP} the two determinant symbols cancel:
\begin{equation}\label{eq:eps-centralizer}
	\epsilon(S,t)=\prod_i\hil{S_i}{t_i}_{K_i}=\eta_{K[\sqrt S]/K}(t)\qquad(t\in K^\times).
\end{equation}
Here we use the notation $K[\sqrt S]=K[X]/(X^2-S)$.

For $\tilde C\in\tGL_j(F)$ lying over $C$ and $h\in \GL_j(F)$ put
\begin{equation}\label{eq:Taction}
	\T_h(\tilde C)=\hil{\det C}{\det h}\,\tilde h^{-1}\tilde C\tilde h,
\end{equation}
where $\tilde h$ is either lift of $h$; the result does not depend on that choice.

\begin{lemma}\label{lem:eps}
	\begin{enumerate}[label=(\alph*)]
		\item $\T_h$ is a continuous automorphism of $\tGL_j(F)$, $\T_{hk}=\T_k\circ\T_h$, and $\T_h(s(C))=\epsilon(C,h)\,s(C^h)$. Consequently $\epsilon(C,hk)=\epsilon(C,h)\epsilon(C^h,k)$.
		\item In $\tGL_{2j}(F)$, $s(m_C)s(d_h)=\epsilon(C,h)\,s(d_h)s(m_{C^h})$.
		\item $\epsilon(C,aI_j)=\hil{\det C}{a}$ for $a\in F^\times$.
		\item $\epsilon(\diag(C_1,C_2),\diag(h_1,h_2))=\epsilon(C_1,h_1)\epsilon(C_2,h_2)$. More precisely, for $\tilde C_1,\tilde C_2$ over $C_1\in\GL_a(F)$, $C_2\in\GL_b(F)$, $\T_{\diag(h_1,h_2)}(\tilde C_1\oplus1)=\T_{h_1}(\tilde C_1)\oplus1$ and $\T_{\diag(h_1,h_2)}(1\oplus \tilde C_2)=1\oplus\T_{h_2}(\tilde C_2)$.
	\end{enumerate}
\end{lemma}
\begin{proof}
	(a) Replacing $\tilde h$ by $-\tilde h$ does not change $\tilde h^{-1}\tilde C\tilde h$. Conjugation is an automorphism, and $\tilde C\mapsto \hil{\det C}{\det h}$ is a character of $\tGL_j(F)$ with values in the center; hence $\T_h$ is an automorphism. $\T_k(\T_h\tilde C)=\hil{\det C}{\det k}\hil{\det C}{\det h}(\tilde h\tilde k)^{-1}\tilde C(\tilde h\tilde k)=\T_{hk}(\tilde C)$ because $\tilde h\tilde k$ lies over $hk$. Next, $s(C)s(h)=\sigma(C,h)s(Ch)$ and $s(h)s(C^h)=\sigma(h,C^h)s(hC^h)=\sigma(h,C^h)s(Ch)$; therefore $s(h)^{-1}s(C)s(h)=\frac{\sigma(C,h)}{\sigma(h,C^h)}s(C^h)$, which is the formula for $\T_h(s(C))$. Applying $\T_{hk}=\T_k\T_h$ to $s(C)$ gives the cocycle identity.
	
	(b) By \ref{C:block}, $s(m_C)s(d_h)=\sigma(C,h)\hil{\det C}{\det h}s(\diag(Ch,h))$ and $s(d_h)s(m_{C^h})=\sigma(h,C^h)s(\diag(hC^h,h))$. Since $Ch=hC^h$, the assertion follows.
	
	(c) Take $h=aI_j$ in (b); then $d_h=aI_{2j}$ and $m_{C^h}=m_C$, so (b) says $s(m_C)s(aI_{2j})=\epsilon(C,aI_j)s(aI_{2j})s(m_C)$. Compare with \ref{C:center} for $g=m_C$.
	
	(d) By \ref{C:block},
	\[\frac{\sigma(\diag(C_1,C_2),\diag(h_1,h_2))}{\sigma(\diag(h_1,h_2),\diag(C_1^{h_1},C_2^{h_2}))}=\frac{\sigma(C_1,h_1)\sigma(C_2,h_2)}{\sigma(h_1,C_1^{h_1})\sigma(h_2,C_2^{h_2})}\cdot\hil{\det C_1}{\det h_2}\hil{\det h_1}{\det C_2}.\]
	Expanding $\hil{\det C_1\det C_2}{\det h_1\det h_2}$ bimultiplicatively, the two cross symbols occur twice in the product and cancel. This gives the first formula; the second follows by applying it to $s(C_1)\oplus 1=s(\diag(C_1,I_b))$ and using (a), and similarly for $1\oplus\tilde C_2$.
\end{proof}

\paragraph{Other completions.}
The definitions \eqref{eq:eps} and \eqref{eq:Taction} make sense over every completion $k_v$ of $k$. The standing field hypothesis in \cite[\S1.1]{KM} is characteristic different from $2$, not odd residual characteristic. Thus the block-compatibility identity \ref{C:block} and the upper-unipotent section identities used in \ref{C:unip} apply over archimedean and dyadic fields as well; see \cite[\S1.2, (1.1) and the paragraph following (1.2)]{KM}.


The scalar commutator \ref{C:center} is also uniform. For $a\in F^\times$ put
\[
c_a(g)=[\tilde g,\widetilde{aI_N}]\in\{\pm1\}.
\]
Since $aI_N$ is central in $\GL_N(F)$, $c_a$ is a character of $\GL_N(F)$. Iterating block compatibility from rank one gives, for $t=\diag(t_1,\ldots,t_N)$,
\[
\sigma_N(t,aI_N)=\prod_{i<k}\hil{t_i}{a},
\qquad
\sigma_N(aI_N,t)=\prod_{i<k}\hil{a}{t_k},
\]
and therefore
\[
c_a(t)=\prod_{i=1}^N\hil{t_i}{a}^{N-1}
=\hil{\det t}{a}^{N-1}.
\]
For an elementary unipotent $u_{ik}(z)$ one has $u_{ik}(z)=u_{ik}(z/2)^2$, hence $c_a(u_{ik}(z))=1$. Since diagonal and elementary unipotent matrices generate $\GL_N(F)$, this proves
\[
[\tilde g,\widetilde{aI_N}]=\hil{\det g}{a}^{N-1}
\]
over every local field of characteristic zero. For $N=2j$ this is exactly \ref{C:center}. Consequently the proof of \Cref{lem:eps}, which uses only \ref{C:block} and \ref{C:center} together with cocycle algebra and Hilbert-symbol identities, applies over $k_v$ for every place $v$. We write $\sigma_v$, $s_v$, $\epsilon_v$ and $\T_{h,v}$ for these local objects. 


\paragraph{The global cover.}
Fix a finite set $\Sigma_{\rm bad}$ of places containing the archimedean and dyadic places and every finite place at which the chosen compact splitting is not available. For $v\notin\Sigma_{\rm bad}$ let $\kappa_v:\GL_j(\cO_v)\to\tGL_j(k_v)$ be the chosen compact splitting. The global double cover $\tGL_j(\A)$ is the quotient of the restricted product of the local covers $\tGL_j(k_v)$, taken with respect to the subgroups $\kappa_v(\GL_j(\cO_v))$ for $v\notin\Sigma_{\rm bad}$, by the subgroup of elements $(\zeta_v)_v$, $\zeta_v\in\{\pm1\}$, with $\prod_v\zeta_v=1$. The rational points $\GL_j(k)$ lift by the rational splitting, which we normalize, as in \cite[\S1.3]{KM} and \cite[Propositions 1.7, 1.8]{Takeda}, to be $\gamma\mapsto\prod_v s_v(\gamma)$. This product makes sense because $\kappa_v=s_v$ on $\gamma$ for almost all $v$, and it is a homomorphism because $\prod_v\sigma_v(\gamma,\gamma')=1$ for $\gamma,\gamma'\in\GL_j(k)$. Combined with Hilbert reciprocity this gives
\begin{equation}\label{eq:productformula}
	\prod_v\epsilon_v(C,\gamma)=1\qquad (C,\gamma\in\GL_j(k)).
\end{equation}
For $C\in\GL_j(k)$ and $a\in\A^\times$, \Cref{lem:eps}(c) gives $\epsilon_v(C,a_vI_j)=\hil{\det C}{a_v}_v$, which is $1$ for almost all $v$, and $\prod_v\epsilon_v(C,a_vI_j)=\eta_{\det C}(a)$. Every unipotent algebraic $k$-subgroup $U$ of $\GL_j$ splits canonically into $\tGL_j(\A)$, by the product of the local canonical lifts of \ref{C:unip} (at almost all places these agree with $\kappa_v$ by \Cref{lem:unip-kappa}); we again write $\tilde u$ for the canonical lift of $u\in U(\A)$. On $N_j(k)$ it agrees with the rational splitting, by \ref{C:unip}.

\section{Exceptional representations and the spherical coefficient}\label{sec:exceptional}

Let $\theta_{2j}$ be the exceptional representation of $\tGL_{2j}(F)$ with trivial ordinary central character \cite{Kable,KM}. It is genuine. By \cite[Theorem 3.1]{Kaplan} the space of functionals $\ell$ with $\ell(\theta_{2j}(\tilde u(X))\xi)=\psi(\tr X)\ell(\xi)$ is one dimensional; we fix such an $\ell$ and call it the Shalika functional. Let $\GL_j^{\rm diag}=\{d_h:\ h\in\GL_j(F)\}\subset\GL_{2j}(F)$ be the diagonally embedded copy of $\GL_j$, and $\tGL_j^{\rm diag}$ its preimage in $\tGL_{2j}(F)$. By \ref{C:block}, $\sigma(d_h,d_{h'})=\sigma(h,h')^2\hil{\det h}{\det h'}=\hil{\det h}{\det h'}$, so $\tGL_j^{\rm diag}$ is the double cover of $\GL_j(F)$ defined by the cocycle $(h,h')\mapsto\hil{\det h}{\det h'}$, i.e.\ the pullback via the determinant of the double cover of $F^\times$ defined by the Hilbert symbol. There is a genuine character $\beta$ of $\tGL_j^{\rm diag}$ (the \emph{diagonal Shalika character}) with $\ell(\theta_{2j}(x)\xi)=\beta(x)\ell(\xi)$ for $x\in\tGL_j^{\rm diag}$. The same applies to the character $\psi^{-1}$; the corresponding objects are written $\ell^-$, $\beta^-$, and by \cite[\S1.6, \S4.2]{KM} one has $\beta(x)\beta^-(x)=1$ for every $x\in\tGL_j^{\rm diag}$. When the data are unramified, $\xi^\circ$ is the vector fixed by $\kappa(\GL_{2j}(\cO))$ and $\ell(\xi^\circ)=\ell^-(\xi^\circ)=1$. For a vector $\xi$ and $\tilde C\in\tGL_j(F)$ over $C$ we put
\begin{equation}\label{eq:Wcoef}
\cB_\xi(\tilde C)=|\det C|^{-\frac{2j-1}4}\,\ell\bigl(\theta_{2j}(\tilde C\oplus1)\xi\bigr),\qquad \cW_\xi(C)=\cB_\xi(s(C)),
\end{equation}
and $\cB=\cB_{\xi^\circ}$, $\cW=\cW_{\xi^\circ}$. The genuine function $\cB$ is the convenient object because it does not depend on a section; $\cW$ is its value on the BLS section. The analogous objects for $\ell^-$ are written $\cB^-,\cW^-$. We use the following spherical formula \cite[Lemma 4.2]{KM}:
\begin{enumerate}[label=(S)]
\item\label{S} $\cW(C)=0$ unless $C$ is integral with all Smith exponents even. If $d_1\ge\dots\ge d_j\ge0$, then $\cW(\diag(\varpi^{2d_1},\dots,\varpi^{2d_j}))=q^{\sum_i(i-1)d_i}$. The same holds for $\cW^-$.
\end{enumerate}
In \Cref{app:technical} we deduce from \ref{S} that $\cB=\cB^-$, that $\cB$ is two sided invariant under $\kappa(\GL_j(\cO))$ on its support, and that $\cB(s(C))=1$ when $C$ is an integral upper Hessenberg matrix with subdiagonal entries $1$ and $\val\det C\in 2\Z_{\ge0}$ (\Cref{lem:compact,lem:hessenberg}).

\chapter{Unfolding and the local comparison}\label{ch:unfolding}

This chapter sets up the comparison and reduces it to a local statement. In \S\ref{sec:unfolding} we write down the two relative trace formulas (\Cref{def:RTF}), define their regular parts, and unfold these into sums of products of local orbital integrals. On the orthogonal side this is an elementary orbit computation for $\mathcal H_W$ (\Cref{prop:orth-unfold}). On the linear side it requires a Fourier expansion of two exceptional theta forms, a partial Poisson summation in the Eisenstein series and the factorization of global Shalika coefficients; the orbits are indexed by enhanced matrix data $z=(A,v,r)$, and each local orbital integral carries the metaplectic sign $\mu(z,h)=\epsilon(A,h)\epsilon(B,h)$  with $B=A+vr-I_n$ (\Cref{thm:linear-unfold}). We use the notation and the background on metaplectic covers and exceptional representations of \Cref{ch:notation}; the analytic justification of the unfolding is postponed to \Cref{app:convergence}.

In \S\ref{sec:FL} we match the regular orbits on the two sides by an explicit map $z\mapsto(W_z,g_z)$ (\Cref{prop:matching}), normalize the linear orbital integrals by the transfer factor $\Delta(z)=\mu(z,C_v)$ attached to the cyclic frame (\Cref{lem:equivariance}), and state the unit fundamental lemma (\Cref{thm:FL}) and the smooth transfer conjecture. The proof of \Cref{thm:FL} occupies Chapters \ref{ch:decomposition}--\ref{ch:fourier}.

\section{The unfolding}\label{sec:unfolding}

Throughout this section $k$ is a number field, with the global notation of \S\ref{sec:notation} and the global double covers of \S\ref{sec:cover}.

\subsection*{The two relative trace formulas}

Let $W$ be a nondegenerate quadratic space over $k$ of dimension $2n$ and $V=W\oplus ke$ with $(e,e)=1$. For $h\in\Ort(W)$ write $\widehat h=\diag(h,1)\in\Ort(V)$. Put
\[
 G_W=\SO(V)\times\SO(W),\qquad
 \mathcal H_W=\{(h_1,h_2)\in\Ort(W)^2:\det h_1=\det h_2\}.
\]
The right action
\[
 (g_1,g_2)\cdot(h_1,h_2)
 =\bigl(\widehat h_1^{-1}g_1\widehat h_2,\,h_1^{-1}g_2h_2\bigr)
\]
preserves $G_W$. This incorporates the full $\Ort(W)$-symmetry while keeping the group variable in $\SO(V)$; compare \cite[\S4B]{Krishna}. We suppress the embedding hat for elements of $\SO(W)$ when no confusion is possible.

Use the measure $dh\,du$ in the coordinates $h_1=h\in\Ort(W)$, $h_2=hu$ with $u\in\SO(W)$, with local $\Ort(W)$-measures as in \S\ref{sec:FL}. At each place $v$, the measure on $\SO(W)(k_v)$ is twice the restriction of the measure on $\Ort(W)(k_v)$; in particular, $\vol\SO(L_{W,v})=1$ whenever $\vol\Ort(L_{W,v})=1$. Thus each local component group $\mu_2(k_v)$ has mass one. Quotients by rational points use counting measure; thus $\vol(\mu_2(k)\backslash\mu_2(\A))=\tfrac12$.


On the linear side let $S=\{d_hu(X)\}\subset\GL_{2n}$ be the Shalika subgroup with character $\psi_S(d_hu(X))=\psi(\tr X)$. Since $\det(\mathrm{Ad}(h)|_{\Mat_n})=1$, the measure $dh\,dX$ is a Haar measure on $S(\A)$. Let $\Theta^+$ and $\Theta^-$ be automorphic forms in the global exceptional representations of $\tGL_{2n}(\A)$ whose Shalika periods are taken with respect to $\psi$ and $\psi^{-1}$ respectively \cite[\S4]{KM}. The product $\Theta^+\Theta^-$ is a function on $\GL_{2n}(k)\backslash\GL_{2n}(\A)$, since the genuine signs cancel. Let $\eta$ be a quadratic character of $k^\times\backslash\A^\times$ and let $\Xi$ be the $\eta$-isotypic part of $\Theta^+\Theta^-$ under the center, so that $\Xi(ag)=\eta(a)\Xi(g)$ for $a\in\A^\times$. For $\Phi\in\cS(\A^{1\times 2n})$ let
\begin{equation}\label{eq:eisenstein}
E(g,s)=|\det g|^s\int_{k^\times\backslash\A^\times}\Bigl(\sum_{w\in k^{1\times 2n}}\Phi(twg)-\Phi(0)\Bigr)\eta(t)|t|^{2ns}\,d^\times t,\qquad \Re s>1,
\end{equation}
be the Rankin--Selberg Eisenstein series \cite[Lecture 5, Proposition 5.4]{Cogdell}. Replacing $t$ by $a^{-1}t$ shows $E(ag,s)=\eta(a)^{-1}E(g,s)$.

\begin{definition}[The two relative trace formulas]\label{def:RTF}
\begin{enumerate}[label=(\alph*)]
\item (Orthogonal side.) For $f\in C_c^\infty(G_W(\A))$ put
\[
 I(f)=\int_{[\mathcal H_W]}\sum_{\gamma\in G_W(k)}
 f(\gamma\cdot\mathbf h)\,d\mathbf h.
\]
An element $g\in\SO(V)$ is called \emph{regular} if $e,ge,\dots,g^{2n}e$ is a basis of $V$. The \emph{regular part} $I^{\rm reg}(f)$ restricts the sum to $\gamma=(\gamma_1,\gamma_2)$ with $\gamma_1\widehat\gamma_2^{-1}$ regular.
\item (Linear side.) Let $Z$ be the center of $\GL_{2n}$. For $f_1,f_2\in C_c^\infty(\GL_{2n}(\A))$ put
\[F_1(g)=\int_{Z(\A)}f_1(zg)\,dz,\qquad F_2(g)=\int_{Z(\A)}f_2(zg)\eta(z)\,dz,\]\[K_i(x,y)=\sum_{\gamma\in Z(k)\backslash\GL_{2n}(k)}F_i(x^{-1}\gamma y),\]
where $\eta(aI_{2n})=\eta(a)$, and
\[J(f_1,f_2,\Phi;s)=\int_{[\GL_{2n}]_Z}\int_{[S]_Z}\int_{[\GL_{2n}]_Z}K_1(x,a)K_2(x,y)E(x,s)\Xi(y)\psi_S(a)^{-1}\,dy\,da\,dx.\]
Its \emph{regular part} $J^{\rm reg}(s)$ is defined in \eqref{eq:Jreg} below.
\end{enumerate}
\end{definition}

\begin{remark}
	(Cyclic vector). We say $e$ is cyclic for $g$ when $e,ge,\dots,g^{2n}e$ is a basis of $V$.
\end{remark}

In each variable of $J$ the central characters cancel: $K_2(x,y)$ transforms by $\eta$ in $x$ and $\eta^{-1}$ in $y$, $E$ by $\eta^{-1}$ and $\Xi$ by $\eta$. Neither $I(f)$ nor $J$ converges for arbitrary data (see \Cref{rem:divergence}). We therefore work with regular parts, which are defined without assuming that $I$ or $J$ converges. On the orthogonal side, $I^{\rm reg}(f)$ is obtained by restricting the rational kernel sum, as in \Cref{def:RTF}(a). On the linear side, regularity is a condition on Fourier frequencies, which appear only after the theta forms are expanded; $J^{\rm reg}(s)$ is defined in \eqref{eq:Jreg} by retaining the regular frequencies inside the Shalika period, \emph{before} the integration over the noncompact quotient $[\GL_n]_Z$ and the id\`ele classes. The relation of $J^{\rm reg}$ with $J$ is explained after \eqref{eq:Jreg}. In (a) the quadratic space $W$ is fixed, whereas (b) depends on the character $\eta$ but on no quadratic space: with $\eta=\eta_\delta$, the linear side is compared with the sum of the orthogonal sides over all quadratic spaces $W$ of dimension $2n$ and discriminant $\delta$ (\Cref{rem:sumW}). Classes of test functions for which the full expressions $I$ and $J$ converge are given in \Cref{prop:A-orth,prop:A-full}.

\subsection*{The orthogonal side}

\begin{proposition}\label{prop:orth-unfold}
$I^{\rm reg}(f)$ converges absolutely, only finitely many rational $\Ort(W)(k)$-conjugacy classes of regular $g\in\SO(V)(k)$ contribute, and
\[
 I^{\rm reg}(f)=\sum_{[g]}\int_{\Ort(W)(\A)}
 f_W(\widehat h^{-1}g\widehat h)\,dh,\qquad
 f_W(x)=\int_{\SO(W)(\A)}f(x\widehat u,u)\,du,
\]
the sum being over regular $g\in\SO(V)(k)$ modulo conjugation by $\Ort(W)(k)$. For factorizable $f$, each summand is the product of the corresponding local $\Ort(W)(k_v)$-orbital integrals.
\end{proposition}

\begin{proof}
The map $(\gamma_1,\gamma_2)\mapsto\gamma_1\widehat\gamma_2^{-1}$ identifies the $\mathcal H_W(k)$-orbits with $\Ort(W)(k)$-conjugacy classes in $\SO(V)(k)$: every orbit has a representative $(g,1)$, and its remaining equivalence is simultaneous conjugation by $h_1=h_2\in\Ort(W)(k)$. The stabilizer of a regular $(g,1)$ is trivial, since an isometry fixing $e$ and commuting with $g$ fixes the basis $e,ge,\dots,g^{2n}e$. Unfolding its rational orbit therefore gives
\[
 \int_{\mathcal H_W(\A)}
 f(\widehat h_1^{-1}g\widehat h_2,h_1^{-1}h_2)\,d(h_1,h_2)
 =\int_{\Ort(W)(\A)}\int_{\SO(W)(\A)}
 f(\widehat h^{-1}g\widehat h\widehat u,u)\,du\,dh,
\]
by $h=h_1$, $u=h_1^{-1}h_2$. This is the claimed orbital integral. The calculation is first made for $|f|$; convergence, finiteness of the contributing rational orbits, and factorization are proved in \Cref{prop:A-orth}.
\end{proof}

\subsection*{The linear side}

The unfolding follows from four lemmas: two global identities (\Cref{lem:kernels,lem:poisson} below) and two local ones (\Cref{lem:factor,lem:covariance}).

The following unfolding Lemma is standard and elementary.
\begin{lemma}[Unfolding the kernels]\label{lem:kernels}
	Let $\Re s>1$ and put $\Xi_{f_2}(x)=\int_{\GL_{2n}(\A)}f_2(b)\,\Xi(xb)\,db$.
	\begin{enumerate}[label=(\alph*)]
		\item For every $x\in\GL_{2n}(\A)$, $\int_{[\GL_{2n}]_Z}K_2(x,y)\Xi(y)\,dy=\Xi_{f_2}(x)$, the integral being absolutely convergent.
		\item For every $a\in S(\A)$ the double integral
		\[P(a)=\int_{[\GL_{2n}]_Z}\int_{[\GL_{2n}]_Z}K_1(x,a)K_2(x,y)E(x,s)\Xi(y)\,dy\,dx\]
		converges absolutely, and $P(a)=\int_{\GL_{2n}(\A)}f_1(g^{-1})\,\Xi_{f_2}(ag)E(ag,s)\,dg$. The function $P$ is left invariant under $Z(\A)S(k)$.
		\item If the triple integral defining $J(f_1,f_2,\Phi;s)$ converges absolutely, then $J(f_1,f_2,\Phi;s)=\int_{[S]_Z}P(a)\psi_S(a)^{-1}\,da$. If moreover
		\[\int_{\GL_{2n}(\A)}\int_{[S]_Z}\bigl|f_1(g^{-1})\,\Xi_{f_2}(ag)E(ag,s)\bigr|\,da\,dg<\infty,\]
		then
		\[J(f_1,f_2,\Phi;s)=\int_{\GL_{2n}(\A)}f_1(g^{-1})\int_{[S]_Z}\Xi_{f_2}(ag)E(ag,s)\psi_S(a)^{-1}\,da\,dg.\]
	\end{enumerate}
\end{lemma}

We remark that the convergence hypotheses in the Lemma hold for the test functions of \Cref{prop:A-full}.

For notational simplicity we assume in the rest of this section that $\Xi_{f_2}$ is the $\eta$-isotypic part of a product $\Theta^+\Theta^-$ of factorizable forms; the general case is a finite sum or a compact integral of such, and the estimates of \Cref{app:convergence} are uniform for such integrals. We write $[\Mat_n]=\Mat_n(k)\backslash\Mat_n(\A)$, with measure of volume one; the map $X\mapsto u(X)$ identifies it with $[U]=U(k)\backslash U(\A)$. We put
\[\Theta^{+}_A(g)=\int_{[\Mat_n]}\Theta^+(\tilde u(X)g)\psi(-\tr AX)\,dX,\qquad \Theta^{-}_B(g)=\int_{[\Mat_n]}\Theta^-(\tilde u(X)g)\psi(\tr BX)\,dX\]
for $A,B\in\Mat_n(k)$. For $g\in\GL_{2n}(\A)$, a row $r$ and a column $v$ put
\begin{equation}\label{eq:partialFT}
	\widehat\Phi_g(r,v)=\int_{\A^{1\times n}}\Phi((r,x)g)\psi(-xv)\,dx.
\end{equation}

\begin{lemma}[Fourier expansion and partial Poisson summation]\label{lem:poisson}
	Let $g\in\GL_{2n}(\A)$, $\Re s>1$, and $C_Y(g)=\int_{[\Mat_n]}\Theta^+\Theta^-(\tilde u(X)g)\psi(-\tr YX)\,dX$. Then $C_Y(g)=\sum_{A-B=Y}\Theta^+_A(g)\Theta^-_B(g)$ and
	\begin{multline*}
		\int_{[\Mat_n]}(\Theta^+\Theta^-)(\tilde u(X)g)E(u(X)g,s)\psi(-\tr X)\,dX\\
		=|\det g|^s\int_{k^\times\backslash\A^\times}\eta(t)\Bigl(|t|^{n(2s-1)}\sum_{r\in k^{1\times n},\,v\in k^n}C_{I_n-vr}(g)\widehat\Phi_g(tr,t^{-1}v)-|t|^{2ns}\Phi(0)C_{I_n}(g)\Bigr)d^\times t .
	\end{multline*}
	In the first term the frequencies satisfy $A-B+vr=I_n$.
\end{lemma}
\begin{proof}
	Since $[\Mat_n]$ is compact with dual group $\Mat_n(k)$ (via the trace pairing), $\Theta^+(\tilde u(X)g)=\sum_A\Theta^+_A(g)\psi(\tr AX)$ and $\Theta^-(\tilde u(X)g)=\sum_B\Theta^-_B(g)\psi(-\tr BX)$, with absolutely convergent series (\Cref{lem:A-fourier}). Multiplying and integrating against $\psi(-\tr YX)$ gives the formula for $C_Y$.
	
	Write $w=(r,x)$ with rows $r,x\in k^{1\times n}$. Fix $t\in\A^\times$ and $r$, and let $f(x)=\Phi((tr,t(x+rX))g)$. With $y=t(x+rX)$, $dx=|t|^{-n}dy$,
	\[\hat f(v)=\int f(x)\psi(-xv)\,dx=|t|^{-n}\int\Phi((tr,y)g)\psi(-(t^{-1}y-rX)v)\,dy=|t|^{-n}\psi(rXv)\widehat\Phi_g(tr,t^{-1}v).\]
	Since $(r,x)u(X)=(r,x+rX)$, Poisson summation $\sum_xf(x)=\sum_v\hat f(v)$ gives
	\[\sum_{w\in k^{1\times2n}}\Phi(twu(X)g)=|t|^{-n}\sum_{r,v}\widehat\Phi_g(tr,t^{-1}v)\psi(rXv).\]
	Insert this into \eqref{eq:eisenstein} at $u(X)g$; note $|\det u(X)g|=|\det g|$. As $\psi(rXv)=\psi(\tr(vrX))$, integration over $[\Mat_n]$ against $(\Theta^+\Theta^-)(\tilde u(X)g)\psi(-\tr X)$ produces $C_{I_n-vr}(g)$ from the first term and $C_{I_n}(g)$ from the term $\Phi(0)$. The interchange of the $X$-integral with the $t$-integral and the sums is justified in \Cref{lem:A-fourier}.
\end{proof}

The first term of \Cref{lem:poisson} is a sum over $(A,B,v,r)$ with $A-B+vr=I_n$; equivalently, over triples $z=(A,v,r)$, the matrix $B$ being determined by $z$. This leads to the following definitions.
	\begin{definition}\label{def:regular}
		Let $K$ be a field. For $A\in\Mat_n(K)$, $v\in K^n$ and $r\in K^{1\times n}$ put
		\[B=A+vr-I_n,\qquad C_v=[v,Av,\dots,A^{n-1}v],\qquad R_r=\begin{pmatrix} r\\ rA\\ \vdots\\ rA^{n-1}\end{pmatrix}.\]
		The triple $z=(A,v,r)$ is called \emph{regular} if $A,B\in\GL_n(K)$ and $\det(R_rC_v)\neq0$. The set of regular triples is denoted $X^{\rm reg}(K)$. The group $\GL_n(K)$ acts on the right by
		\[z\cdot h=(A^h,h^{-1}v,rh),\]
		and $B(z\cdot h)=B^h$. The matrix $R_rC_v=(rA^{i+j}v)_{0\le i,j<n}$ is the Hankel matrix of $z$.
	\end{definition}
	
	\begin{lemma}\label{lem:free}
		Let $z\in X^{\rm reg}(K)$. Then $C_v,R_r\in\GL_n(K)$, $C_{h^{-1}v}(z\cdot h)=h^{-1}C_v$ and $R_{rh}(z\cdot h)=R_rh$. The stabilizer of $z$ in $\GL_n(K)$ is trivial, i.e.\ $\GL_n(K)$ acts freely on $X^{\rm reg}(K)$ and $h\mapsto z\cdot h$ is a bijection from $\GL_n(K)$ onto the orbit of $z$; we express this by saying that the regular orbits are \emph{free}. Moreover $z_c=z\cdot C_v$ is the unique point of the orbit of $z$ of the form $(A_c,e_1,r_c)$ with $A_c$ a companion matrix. Explicitly $A_c=C_v^{-1}AC_v$ is the companion matrix of the characteristic polynomial of $A$ (with subdiagonal entries $1$) and $r_c=(rv,rAv,\dots,rA^{n-1}v)$.
	\end{lemma}
	\begin{proof}
		$\det(R_rC_v)\ne0$ forces both factors to be invertible. The transformation rules follow from $(A^h)^ih^{-1}v=h^{-1}A^iv$ and $rh(A^h)^i=rA^ih$. If $z\cdot h=z$ then $h$ commutes with $A$ and fixes $v$, hence fixes the basis $A^iv$, so $h=1$. In the basis $C_v$ the vector $v$ becomes $e_1$, $A$ maps $e_i$ to $e_{i+1}$ for $i<n$, and $r$ becomes $rC_v$. Conversely, if $z\cdot h$ has vector $e_1$ and companion matrix, then $h^{-1}C_v=C_{e_1}(A^h)=I_n$, i.e.\ $h=C_v$.
	\end{proof}

	For $\delta\in k^\times$ let $X^{\rm reg}_\delta(k)$ be the set of $z\in X^{\rm reg}(k)$ with $\det A\det B\in\delta k^{\times2}$; we fix $\delta$ with $\eta=\prod_v\hil{\cdot}{\delta}_v$.
	 By \ref{C:center} and the global factorization below, $\Theta^+_A\Theta^-_B(ag)=\eta_{\det A\det B}(a)\Theta^+_A\Theta^-_B(g)$ for full rank $A,B$, where $\eta_c=\prod_v\hil\cdot c_v$; hence the $\eta$-isotypic projection retains among the full rank pairs $(A,B)$ those with $\det A\det B\in\delta k^{\times2}$ in the sum in the first term of \Cref{lem:poisson} . We define the regular part by retaining only $z\in X^{\rm reg}_\delta(k)$:
	 \begin{definition}
	\begin{multline}\label{eq:Jreg}
		J^{\rm reg}(s)=\int_{\GL_{2n}(\A)}f_1(g^{-1})\int_{[\GL_n]_Z}\int_{k^\times\backslash\A^\times}\eta(t)|t|^{n(2s-1)}\\
		\times\sum_{z\in X^{\rm reg}_\delta(k)}|\det d_hg|^s\,\Theta^+_A\Theta^-_B(d_hg)\,\widehat\Phi_{d_hg}(tr,t^{-1}v)\,d^\times t\,dh\,dg .
	\end{multline}
		 \end{definition}
	Here we wrote $a=u(X)d_h$ in $[S]_Z$, so that $\psi_S(a)=\psi(\tr X)$, $[S]_Z$ fibers over $[\GL_n]_Z$ with fiber $[U]$, and $da=dX\,dh$. The omitted terms are the two axes $r=0$ and $v=0$ together with the subtracted term $\Phi(0)C_{I_n}$, the singular frequencies $\det A\det B=0$, and the noncyclic data off the axes; their treatment is a separate problem.
	

\begin{remark}
	 The global statements of this paper concern $J^{\rm reg}$ and $I^{\rm reg}$. A comparison of the full distributions $J$ and $I$ would require a regularization of the boundary terms, which is not attempted here. Suppose that $J$ satisfies the two hypotheses of \Cref{lem:kernels}(c), and that the omitted terms, integrated over $[\GL_n]_Z$, $k^\times\backslash\A^\times$ and $\GL_{2n}(\A)$ as in \eqref{eq:Jreg}, converge absolutely. Then $J=J^{\rm reg}+J^{\rm bd}$, where $J^{\rm bd}$ is the integral of the omitted terms. For the test functions of \Cref{prop:A-full}, these hypotheses hold, the omitted terms vanish identically, and $J=J^{\rm reg}$. .
	 \end{remark}
	 
We will now unfold the distribution $J^{\rm reg}$.	 \Cref{thm:A-regular} proves that \eqref{eq:Jreg} converges absolutely for every $s$, so the order of its integrations is immaterial.

We pass from global Fourier coefficients to local functionals. For a local field and $z=(A,v,r)$ with $A,B$ invertible, vectors $\xi^\pm$ of the local exceptional representations $\theta^\pm$ and $\Phi\in\cS(F^{1\times2n})$, put (products of genuine expressions use one and the same lift $\tilde g$ of $g$)
\begin{equation}\label{eq:Lz}
L_z(g,s)=|\det g|^s\,\ell^+\bigl(\theta^+(s(m_A))\theta^+(\tilde g)\xi^+\bigr)\,\ell^-\bigl(\theta^-(s(m_B))\theta^-(\tilde g)\xi^-\bigr)\,\widehat\Phi_g(r,v),
\end{equation}
and $L_z(g)=L_z(g,\tfrac12)$.

\begin{lemma}[Factorization]\label{lem:factor}
Normalize the local Shalika functionals so that the global Shalika coefficient $\Theta^+_{I_n}(g)$ equals $\prod_v\ell_v(\theta^+_v(g_v)\xi^+_v)$ for $\Theta^+=\otimes\xi_v^+$, and similarly for $\Theta^-$. Then for $A,B\in\GL_n(k)$ and adelic $v,r$,
\[|\det g|^s\Theta^+_A(g)\Theta^-_B(g)\widehat\Phi_g(r,v)=\prod_vL_{z,v}(g_v,s).\]
\end{lemma}
\begin{proof}
By \cite[Theorem 4.22]{KM}, the global Shalika coefficient is nonzero and factorizable. Its nonarchimedean local factors are unique up to scalars by \cite[Theorem 3.1 and Remark 3.2]{Kaplan}; for the archimedean factors, see \cite[§4.2]{KM}. Hence $\Theta^+_{I_n}(g)=c\prod_v\ell_v(\theta_v^+(g_v)\xi^+_v)$ with a constant $c$ independent of the vector; at one place we rescale $\ell_v$ so that $c=1$, keeping $\ell_v(\xi_v^\circ)=1$ at the unramified places. Let $j(m_A)=\prod_vs_v(m_A)$ be the rational lift. By automorphy, $\Theta^+(\tilde u(X)g)=\Theta^+(j(m_A)\tilde u(X)g)=\Theta^+(\tilde u(AX)j(m_A)g)$, using \ref{C:unip}. Substituting $Y=AX$, which preserves $[U]$ and its measure since $|\det A|_\A=1$,
\[\Theta^+_A(g)=\int_{[U]}\Theta^+(\tilde u(Y)j(m_A)g)\psi(-\tr Y)\,dY=\Theta^+_{I_n}(j(m_A)g)=\prod_v\ell_v(\theta_v^+(s_v(m_A))\theta^+_v(g_v)\xi^+_v).\]
The same argument applies to $\Theta^-_B$, and $\widehat\Phi$ and $|\det g|^s$ factor for factorizable $\Phi$.
\end{proof}

	For $z=(A,v,r)$ with $A,B$ invertible and a place (or a local field) put
\begin{equation}\label{eq:mu}
	\mu(z,h)=\epsilon(A,h)\epsilon(B,h).
\end{equation}
By \Cref{lem:eps}(a), $\mu(z,hk)=\mu(z,h)\mu(z\cdot h,k)$, and by \Cref{lem:eps}(c), $\mu(z,aI_n)=\hil{a}{\det A\det B}$.

\begin{lemma}[Local covariance]\label{lem:covariance}
Over a local field, for $z=(A,v,r)$ with $A,B$ invertible, $X\in\Mat_n(F)$ and $h\in\GL_n(F)$,
\[L_z(u(X)g,s)=\psi(\tr X)L_z(g,s),\qquad L_z(d_hg,s)=\mu(z,h)|\det h|^{2s-1}L_{z\cdot h}(g,s).\]
\end{lemma}
\begin{proof}
By \ref{C:unip}, $s(m_A)\tilde u(X)s(m_A)^{-1}=\tilde u(AX)$, so the first factor of \eqref{eq:Lz} at $u(X)g$ is multiplied by $\psi(\tr AX)$; similarly the second by $\psi(-\tr BX)$. Since $(r,x)u(X)=(r,x+rX)$, the substitution $x\mapsto x-rX$ in \eqref{eq:partialFT} gives $\widehat\Phi_{u(X)g}(r,v)=\psi(rXv)\widehat\Phi_g(r,v)$. The product of the three characters is $\psi(\tr((A-B+vr)X))=\psi(\tr X)$.

For $d_h$ use \Cref{lem:eps}(b): $\theta^+(s(m_A))\theta^+(s(d_h))=\epsilon(A,h)\theta^+(s(d_h))\theta^+(s(m_{A^h}))$. Hence, with the lift $s(d_h)\tilde g$ of $d_hg$,
\[\ell(\theta^+(s(m_A))\theta^+(s(d_h)\tilde g)\xi^+)=\epsilon(A,h)\beta(s(d_h))\ell(\theta^+(s(m_{A^h}))\theta^+(\tilde g)\xi^+),\]
and similarly for $B$ with $\beta^-(s(d_h))$. The product $\beta(s(d_h))\beta^-(s(d_h))$ is $1$. Next, substituting $x=yh^{-1}$ in \eqref{eq:partialFT},
\[\widehat\Phi_{d_hg}(r,v)=\int\Phi((rh,xh)g)\psi(-xv)\,dx=|\det h|^{-1}\widehat\Phi_g(rh,h^{-1}v),\]
and $|\det d_hg|^s=|\det h|^{2s}|\det g|^s$. Since $z\cdot h=(A^h,h^{-1}v,rh)$ and $B(z\cdot h)=B^h$, the lemma follows.
\end{proof}

For a local field, $f\in C_c^\infty(\GL_{2n}(F))$ and the data of \eqref{eq:Lz} define the local kernel and the raw local orbital integral
\begin{equation}\label{eq:rawOI}
K(z,s)=\int_{\GL_{2n}(F)}f(g^{-1})L_z(g,s)\,dg,\qquad O_z(s)=\int_{\GL_n(F)}\mu(z,h)|\det h|^{2s-1}K(z\cdot h,s)\,dh .
\end{equation}
They converge absolutely for $z\in X^{\rm reg}(F)$ (\Cref{prop:A-local}). We write $K(z)=K(z,\frac12)$ and $O_z=O_z(\frac12)$.

\begin{theorem}[Linear unfolding]\label{thm:linear-unfold}
The integral \eqref{eq:Jreg} converges absolutely for every $s\in\C$, locally uniformly in $s$, and hence defines an entire function of $s$; for factorizable data
\[J^{\rm reg}(s)=\sum_{z}\prod_vO_{z,v}(s),\]
the sum being over $z\in X_\delta^{\rm reg}(k)$ modulo $\GL_n(k)$, with $O_{z,v}$ formed with $f_{1,v}$, $\xi^\pm_v$ and $\Phi_v$. In particular $J^{\rm reg}(\frac12)=\sum_z\prod_vO_{z,v}$.
\end{theorem}

\begin{proof}
We first carry out the computation formally and then indicate where it is justified.

\emph{Step 1: the local form of the integrand.} By \Cref{lem:factor} the integrand of \eqref{eq:Jreg} (without $f_1$) at $(t,h)$ equals
\[\eta(t)|t|^{n(2s-1)}\sum_{z}\prod_vL_{z\cdot t,v}(d_hg,s),\qquad z\cdot t=(A,t^{-1}v,tr).\]
By \Cref{lem:eps}(c), $\mu_v(z,t)=\hil{t_v}{\det A\det B}_v$, so $\prod_v\mu_v(z,t)=\eta(t)$ on $X_\delta^{\rm reg}(k)$. By \Cref{lem:covariance} applied at each place with the scalar $t$, $\prod_vL_{z,v}(d_td_hg,s)=\eta(t)|t|^{n(2s-1)}\prod_vL_{z\cdot t,v}(d_hg,s)$. Thus the integrand is $\Lambda(th)$, where
\[\Lambda(h)=\sum_{z\in X^{\rm reg}_\delta(k)}\prod_vL_{z,v}(d_hg,s).\]

\emph{Step 2: combining the centre with $[\GL_n]_Z$.} For $\gamma\in\GL_n(k)$, \Cref{lem:covariance} and \eqref{eq:productformula} give $\prod_vL_{z,v}(d_\gamma d_hg,s)=\prod_vL_{z\cdot\gamma,v}(d_hg,s)$, because $|\det\gamma|_\A=1$. As $z\mapsto z\cdot\gamma$ permutes $X_\delta^{\rm reg}(k)$, $\Lambda$ is left $\GL_n(k)$-invariant. With compatible quotient measures, $\int_{[\GL_n]_Z}\int_{k^\times\backslash\A^\times}\Lambda(th)\,d^\times t\,dh=\int_{[\GL_n]}\Lambda(h)\,dh$.

\emph{Step 3: unfolding the rational orbits.} By \Cref{lem:free} the action of $\GL_n(k)$ on $X^{\rm reg}_\delta(k)$ is free, so
\[\int_{[\GL_n]}\Lambda(h)\,dh=\sum_{[z]}\int_{\GL_n(\A)}\prod_vL_{z,v}(d_hg,s)\,dh.\]

\emph{Step 4: factorization.} By \Cref{lem:covariance}, $L_{z,v}(d_hg,s)=\mu_v(z,h)|\det h|_v^{2s-1}L_{z\cdot h,v}(g,s)$. Integrating against $f_1(g^{-1})=\prod_vf_{1,v}(g_v^{-1})$ and exchanging the $g$- and $h$-integrals gives $\prod_vO_{z,v}(s)$.

\emph{Justification.} All integrands above are dominated by the nonnegative majorant constructed in \Cref{thm:A-regular}, which is integrable for every $s$, locally uniformly; therefore Tonelli and Fubini justify Steps 1--4, the sum over orbits converges normally, and the result is entire. The local integrals converge by \Cref{prop:A-local}.
\end{proof}



\section{Matching, the fundamental lemma and smooth transfer}\label{sec:FL}

In this section all objects are local, over $F$.

\subsection*{The quadratic form attached to a framed matrix}

\begin{lemma}\label{lem:Q0}
Let $S\in\Mat_j(F)$, $x\in F^j$, $y\in F^{1\times j}$ and assume that the Hankel matrix $(yS^{a+b}x)_{0\le a,b<j}$ is invertible. Then there is a unique symmetric matrix $H$ with
\[S^{\mathsf t}H=HS,\qquad Hx=y^{\mathsf t},\]
namely $H(S^ax,S^bx)=yS^{a+b}x$; it is nondegenerate. If $g\in\GL_j(F)$, the matrix attached to $(S^g,g^{-1}x,yg)$ is $g^{\mathsf t}Hg$.
\end{lemma}
\begin{proof}
Invertibility of the Hankel matrix $R_yC_x$ (notation as in \Cref{def:regular}) forces $C_x=[x,\dots,S^{j-1}x]$ to be invertible, so the $S^ax$ form a basis and $H$ may be defined by the displayed formula; its Gram matrix in this basis is the Hankel matrix, so $H$ is symmetric and nondegenerate. For $a,b<j$, $H(S^ax,S\cdot S^bx)=yS^{a+b+1}x=H(S\cdot S^ax,S^bx)$ when $a+1,b+1<j$; if $b=j-1$ write $S^jx=\sum_c\lambda_cS^cx$ by Cayley--Hamilton, and both sides equal $\sum_c\lambda_cyS^{a+c}x$, using the same expansion for $yS^{a+j}x$. Hence $S$ is $H$-selfadjoint, i.e.\ $S^{\mathsf t}H=HS$. Also $H(x,S^bx)=yS^bx$ for all $b$, i.e.\ $x^{\mathsf t}H=y$. Uniqueness: these two identities determine $H(S^ax,S^bx)=x^{\mathsf t}HS^{a+b}x=yS^{a+b}x$. The last statement follows from uniqueness, since $g^{\mathsf t}Hg$ satisfies the two identities for the transformed data.
\end{proof}
Note that $C_x^{\mathsf t}H C_x$ is just the Hankel matrix $(yS^{a+b}x)_{0\le a,b<j}$.

For $z=(A,v,r)\in X^{\rm reg}(F)$ let $Q_0=Q_0(z)$ be the matrix of \Cref{lem:Q0} for $(A,v,r)$:
\begin{equation}\label{eq:Q0}
A^{\mathsf t}Q_0=Q_0A,\qquad Q_0v=r^{\mathsf t}.
\end{equation}
Then also $B^{\mathsf t}Q_0=Q_0B$, because $B^{\mathsf t}Q_0=(A^{\mathsf t}+r^{\mathsf t}v^{\mathsf t}-I)Q_0=Q_0A+r^{\mathsf t}r-Q_0$ and $Q_0vr=r^{\mathsf t}r$.

\subsection*{Matching of regular orbits}

For $z\in X^{\rm reg}(F)$ let $E=F^n$ and
\begin{equation}\label{eq:Wz}
W_z=E\oplus E,\qquad Q_{W_z}=\diag(Q_0A,\,-Q_0B),\qquad V_z=W_z\oplus Fe,\ (e,e)=1,
\end{equation}
and
\begin{equation}\label{eq:gz}
g_z=\begin{pmatrix}2A-I&-2B&2v\\-2A&I+2B&-2v\\-2rA&2rB&1-2rv\end{pmatrix}.
\end{equation}
The forms are symmetric by \eqref{eq:Q0} and nondegenerate because $A,B$ are invertible. The signed discriminant of $W_z$ is
\begin{equation}\label{eq:disc}
(-1)^n\det Q_{W_z}=\det(Q_0)^2\det A\det B\equiv\det A\det B\pmod{F^{\times2}}.
\end{equation}

\begin{proposition}[Matching]\label{prop:matching}
	\begin{enumerate}[label=(\alph*)]
		\item $g_z\in\SO(V_z)$, $\det(g_z+I)=2^{2n+1}\det A$, and $e$ is cyclic for $g_z$.
		\item The Cayley transform $X_z=(g_z-I)(g_z+I)^{-1}$ is
		\[X_z=\begin{pmatrix}0&-A^{-1}B&A^{-1}v\\-I&0&0\\-r&0&0\end{pmatrix}=\begin{pmatrix}C&w\\-w^{\mathsf t}Q_{W_z}&0\end{pmatrix},\]\[C=\begin{pmatrix}0&-A^{-1}B\\-I&0\end{pmatrix},\ w=\begin{pmatrix}A^{-1}v\\0\end{pmatrix},\]
		with $C\in\so(W_z)$ invertible and $w$ cyclic for $C$.
\item The map $z\mapsto(W_z,g_z)$ induces a bijection from $X^{\rm reg}(F)/\GL_n(F)$ onto the set of isometry classes of pairs $(W,g)$, where $W$ is a nondegenerate quadratic space of dimension $2n$ and $g\in\SO(W\oplus Fe)$ is such that $g+I$ is invertible and the Cayley pair $(C,w)$ of $g$ is cyclic; isometries are taken in $\Ort(W)$, extended by $e\mapsto e$.
\item If $e$ is cyclic for $g\in\SO(W\oplus Fe)$, then $g+I$ is invertible and the Cayley pair of $g$ is cyclic. Consequently, by (a)--(c), $z\mapsto(W_z,g_z)$ induces a bijection from $X^{\rm reg}(F)/\GL_n(F)$ onto the set of isometry classes of pairs $(W,g)$ with $\dim W=2n$ and $e$ cyclic for $g$, i.e.\ onto the regular orthogonal orbits of \S\ref{sec:unfolding}, taken over all quadratic spaces $W$ of dimension $2n$.
\end{enumerate}
\end{proposition}

\begin{proof}
(a) Define $j:E\to V_z$ by $j(x)=(x,x,rx)$. For $(y_1,y_2,t)\in V_z$, using $r^{\mathsf t}=Q_0v$,
\[\bigl(j(x),(y_1,y_2,t)\bigr)=x^{\mathsf t}Q_0Ay_1-x^{\mathsf t}Q_0By_2+(rx)t=x^{\mathsf t}Q_0(Ay_1-By_2+vt).\]
So the adjoint of $j$, with respect to $Q_0$ on $E$, is $j^*(y_1,y_2,t)=Ay_1-By_2+vt$, and $j^*j=A-B+vr=I_n$. Therefore $P=jj^*$ is the orthogonal projection onto $j(E)$ and $s_1=I-2P$ is an isometry of determinant $(-1)^n$. Also $s_0=\diag(-I,I,1)$ is an isometry of determinant $(-1)^n$. Writing $u=Ay_1-By_2+vt$ one computes $s_0s_1(y_1,y_2,t)=(-y_1+2u,\,y_2-2u,\,t-2ru)$, which is \eqref{eq:gz}. Hence $g_z=s_0s_1\in\SO(V_z)$. In $g_z+I$ add the first block row to the second; the second becomes $(0,2I,0)$, so $\det(g_z+I)=2^n\det\left(\begin{smallmatrix}2A&2v\\-2rA&2-2rv\end{smallmatrix}\right)=2^{2n+1}\det A$ after adding $r$ times the first row to the last. Cyclicity of $e$ is proved after (b).

(b) We check $X_z(g_z+I)=g_z-I$. The second and third block rows of $X_z(g_z+I)$ are $-(2A,-2B,2v)$ and $-r(2A,-2B,2v)$, which are the corresponding rows of $g_z-I$. The first block row is $-A^{-1}B(-2A,2I+2B,-2v)+A^{-1}v(-2rA,2rB,2-2rv)$, whose entries are
\[2A^{-1}(B-vr)A=2(A-I),\quad -2A^{-1}(I+B-vr)B=-2B,\quad 2A^{-1}(B+I-vr)v=2v,\]
using $B-vr=A-I$. Next, $-w^{\mathsf t}Q_{W_z}=(-(A^{-1}v)^{\mathsf t}Q_0A,0)=(-v^{\mathsf t}Q_0,0)=(-r,0)$, since $A^{-\mathsf t}Q_0=Q_0A^{-1}$. $C$ is skewadjoint because the Cayley transform of an isometry is skewadjoint, and it is invertible because $A^{-1}B$ is. Finally $C^{2i}w=((A^{-1}B)^iA^{-1}v,0)$ and $C^{2i+1}w=(0,-(A^{-1}B)^iA^{-1}v)$; since $A^{-1}B=I-A^{-1}+(A^{-1}v)r$ differs from $I-A^{-1}$ by an operator with image in $F\,A^{-1}v$, induction shows that $A^{-1}v,\dots,(A^{-1}B)^iA^{-1}v$ and $A^{-1}v,\dots,A^{-i}A^{-1}v$ span the same space. The latter vectors span $E$ for $i=n-1$ because $v$ is cyclic for $A$, hence for $A^{-1}\in F[A]$. So $w$ is cyclic for $C$.

Cyclicity of $e$ for $g_z$: $F[g_z]=F[X_z]$, and the span $M_i$ of $e,X_ze,\dots,X_z^ie$ equals $Fe+\mathrm{span}(w,\dots,C^{i-1}w)$, by induction using $X_z(C^aw)=C^{a+1}w-(w^{\mathsf t}Q_{W_z}C^aw)e$. Hence $M_{2n}=V_z$.

(c) If $z'=z\cdot h$ then $Q_0(z')=h^{\mathsf t}Q_0h$ by \Cref{lem:Q0}, so $\diag(h,h,1)$ is an isometry from $V_{z'}$ to $V_z$ carrying $g_{z'}$ to $g_z$. Thus the map is defined on orbits. Conversely, let $(W,g)$ be as in (c), and let $(C,w)$ be its Cayley pair. Since $C$ is skewadjoint, $(C^aw,C^bw)=(-1)^a(w,C^{a+b}w)$, and $(w,C^{c}w)=0$ for odd $c$. Hence $E_{\rm ev}=\mathrm{span}(C^{2i}w)$ and $E_{\rm odd}=\mathrm{span}(C^{2i+1}w)$ are orthogonal; by cyclicity each has dimension $n$, and both are nondegenerate. $C$ is invertible: it is cyclic, so $\dim\ker C\le1$, while $\ker C$ is the radical of the alternating form $(x,Cy)$ and hence has even codimension in the even-dimensional space $W$, so $\dim\ker C$ is even. Using the basis $C^{2i}w$ of $E_{\rm ev}$ and the basis $-C^{2i+1}w$ of $E_{\rm odd}$, we get
\[C=\begin{pmatrix}0&T\\-I&0\end{pmatrix},\qquad w=\begin{pmatrix}w_0\\0\end{pmatrix},\qquad Q_W=\diag(Q_1,Q_1T)\]
for an invertible $T$ with $w_0$ cyclic, and a symmetric nondegenerate $Q_1$ (the form on $E_{\rm odd}$ is $(-Cx,-Cy)=-(x,C^2y)=(x,Q_1T\,\cdot)$). Since $g+I$ is invertible, $I-X$ is invertible, and $\det(I-X)=\det(I+T+w_0r)$ for $r=w_0^{\mathsf t}Q_1$ (compute the determinant by eliminating the last row and the second block column). Put
\[A=(I+T+w_0r)^{-1},\quad v=Aw_0,\quad B=-AT,\quad Q_0=Q_1A^{-1}.\]
Then $A+vr-I=A(I+w_0r)-I=A(A^{-1}-T)-I=-AT=B$; $Q_0=Q_1+Q_1T+r^{\mathsf t}r$ is symmetric because $Q_1T$ is; $Q_0A=Q_1$ is symmetric, i.e.\ $A^{\mathsf t}Q_0=Q_0A$; $Q_0v=Q_1w_0=r^{\mathsf t}$; $Q_0A=Q_1$ and $-Q_0B=Q_1T$; and $-A^{-1}B=T$, $A^{-1}v=w_0$. So the triple $z=(A,v,r)$ reproduces $(W,Q_W,C,w)$ through \eqref{eq:Wz} and (b), hence $(W,g)$. It is regular: $v$ is cyclic for $A$ because $w_0$ is cyclic for $T$, hence for $A^{-1}=I+T+w_0r$ by the rank one argument of (b), and $F[A]=F[A^{-1}]$; the Hankel matrix equals $C_v^{\mathsf t}Q_0C_v$, which is invertible. Two triples giving isometric pairs are in one orbit: an isometry carrying $(C,w)$ to $(C',w')$ carries the bases $C^aw$ to $C'^aw'$, so the matrices $T,w_0,Q_1$ coincide and so do the triples up to the change of basis. This proves (c).

(d) Let $e$ be cyclic for $g\in\SO(W\oplus Fe)$ and let $U$ be its generalized $-1$-eigenspace. This space is nondegenerate, and $\dim U$ is even: all eigenvalues other than $\pm1$ occur in reciprocal pairs and $\det g=1$. If $U\ne0$, cyclicity gives $\dim\ker((g+I)|_U)=1$. Since $(g-I)|_U$ is invertible, $Y=(g+I)(g-I)^{-1}|_U$ is skewadjoint of rank $\dim U-1$, which is odd. This is impossible because $(u,u')\mapsto(u,Yu')$ is alternating and hence has even rank. Thus $U=0$ and $g+I$ is invertible. The span calculation at the end of the proof of (b), read backwards, then shows that the Cayley pair $(C,w)$ is cyclic: $F[g]=F[X]$, and the span of $e,Xe,\dots,X^ie$ is $Fe+\mathrm{span}(w,\dots,C^{i-1}w)$, so cyclicity of $e$ for $g$ forces cyclicity of $w$ for $C$. Conversely, by (a) and (b), $e$ is cyclic for every $g_z$. Together with (c) this proves the last assertion.
\end{proof}

By \eqref{eq:disc}, the matching respects discriminants: $\det A\det B$ is the signed discriminant of $W_z$.

\subsection*{Orbital integrals}

\paragraph{Orthogonal side.} For $f\in C_c^\infty(\SO(W\oplus Fe))$ and $g$ as in \Cref{prop:matching}(c), put
\begin{equation}\label{eq:orthOI}
	O^{\rm orth}(g,f)=\int_{\Ort(W)}f(\widehat h^{-1}g\widehat h)\,dh,
\end{equation}
where $\widehat h=\diag(h,1)$ fixes $e$. If $W$ has a selfdual lattice $L_W$, the measure gives $\Ort(L_W)$ volume one; otherwise we fix a Haar measure. This is precisely the local orbital integral arising from \Cref{prop:orth-unfold}. It depends only on the $\Ort(W)$-conjugacy class of $g$, which is the class determined by the linear matching. For a family $f=(f_W)$, indexed by the isometry classes of $2n$-dimensional quadratic spaces, put $I_z(f)=O^{\rm orth}(g_z,f_{W_z})$, using any isometry of $W_z$ with the chosen representative of its isometry class.


\paragraph{Linear side.} Let $\theta^\pm$ be the local exceptional representations of \S\ref{sec:exceptional} (for $\psi$ and $\psi^{-1}$), $\Omega=\theta^+\otimes\theta^-$, which is an honest representation of $\GL_{2n}(F)$ by $\Omega(g)=\theta^+(\tilde g)\otimes\theta^-(\tilde g)$. For $f\in C_c^\infty(\GL_{2n}(F))$, $\omega\in\Omega$ and $\Phi\in\cS(F^{1\times 2n})$ let $K_{f,\omega,\Phi}(z)$ and $O_z(f,\omega,\Phi)$ be given by \eqref{eq:rawOI} at $s=\frac12$, with $\ell\otimes\ell^-$ applied to $\omega$ in place of the product in \eqref{eq:Lz}.

The orbital integral $O_z$ is not a function on orbits; its value depends on the representative $z$. The following Lemma allows us to further normalize the orbital integral.
\begin{lemma}\label{lem:equivariance}
For $z\in X^{\rm reg}(F)$ and $k\in\GL_n(F)$, $O_{z\cdot k}=\mu(z,k)^{-1}O_z$. Put
\[\Delta(z)=\mu(z,C_v)=\epsilon(A,C_v)\epsilon(B,C_v).\]
Then $\Delta(z\cdot h)=\mu(z,h)^{-1}\Delta(z)$, $\Delta(z_c)=1$, and
\[\Delta(z)^{-1}O_z=O_{z_c}=\int_{\GL_n(F)}\frac{K(z\cdot h)}{\Delta(z\cdot h)}\,dh .\]
\end{lemma}
\begin{proof}
By the cocycle identity $\mu(z\cdot k,h)=\mu(z,k)^{-1}\mu(z,kh)$, so $O_{z\cdot k}=\mu(z,k)^{-1}\int\mu(z,kh)K(z\cdot kh)\,dh=\mu(z,k)^{-1}O_z$ after $h\mapsto k^{-1}h$. By \Cref{lem:free}, $C_{h^{-1}v}(z\cdot h)=h^{-1}C_v$, so $\Delta(z)=\mu(z,h\cdot h^{-1}C_v)=\mu(z,h)\Delta(z\cdot h)$. For $z_c$ the frame is $I_n$. The first displayed equality is the first statement with $k=C_v$; the second follows from $\mu(z,h)=\Delta(z)/\Delta(z\cdot h)$.
\end{proof}

\begin{definition}
	We define the \emph{normalized linear orbital integral}
	\begin{equation}\label{eq:Jz}
		J_z(f,\omega,\Phi)=|\det A\det B|^{-\frac{2n-1}4}\,\Delta(z)^{-1}O_z(f,\omega,\Phi).
	\end{equation}
\end{definition} 
It is invariant under $z\mapsto z\cdot h$, (an essential fact used in the proof of fundamental lemma). Over a number field and for $z\in X^{\rm reg}(k)$, $\prod_v\Delta_v(z)=1$ by \eqref{eq:productformula} and $\prod_v|\det A\det B|_v=1$, so \Cref{thm:linear-unfold} reads
\[J^{\rm reg}(\tfrac12)=\sum_{[z]}\prod_vJ_{z,v}(f_{1,v},\omega_v,\Phi_v).\]
By \Cref{prop:orth-unfold,prop:matching}, the extended orthogonal expansion is
\[
 \sum_W I_W^{\rm reg}(f_W)
 =\sum_{[z]\in X^{\rm reg}_\delta(k)/\GL_n(k)}
   \prod_v I_{z,v}(f_{W_z,v}),
\]
where $f_{W_z,v}$ in the local orbital integral denotes the collapsed test function. The sum over $W$ runs over the rational quadratic spaces of dimension $2n$ and discriminant $\delta$. 

\begin{remark}[One linear side, a sum of orthogonal sides]\label{rem:sumW}
	\begin{enumerate}[label=(\alph*)]
		\item \emph{The global comparison.} Fix $\delta\in k^\times/k^{\times2}$ and take $\eta=\eta_\delta$ on the linear side. On the orthogonal side the test data are then a family $f=(f_W)_W$, consisting of one test function $f_W$ on $(\SO(V)\times\SO(W))(\A)$, $V=W\oplus ke$, for each quadratic space $W$ over $k$ of dimension $2n$ and discriminant $\delta$, and the regular part $J^{\rm reg}(\frac12)$ is compared with
		\[\sum_WI_W^{\rm reg}(f_W).\]
		Here $I_W$ is the orthogonal distribution defined using $\mathcal H_W$ in Definition~\ref{def:RTF}(a).
		The reason is the matching of orbits. For rational $z$ the symmetric matrix $Q_0$ of \Cref{lem:Q0} is rational, so $W_z$ is a quadratic space over $k$, of discriminant $\det A\det B$ modulo squares. By \Cref{prop:matching}, applied over $k$, the regular rational orbits on the linear side with $\det A\det B\in\delta k^{\times2}$ are in bijection with the disjoint union, over all such $W$, of the regular rational $\Ort(W)(k)$-orbits: each linear orbit corresponds to an orbit on exactly one $W$. A single space $W$ only sees the linear orbits with $W_z\simeq W$.
		\item \emph{Local matching.} Choose the factorizable family consistently from local test data indexed by local isometry classes, so that, under the chosen identifications, $f_{W,v}$ depends only on $W_v$. The local comparison at $z$ is between $J_{z,v}$ and the orthogonal orbital integral on $(W_z)_v$. This is why the smooth transfer conjecture (\Cref{conj:transfer}) is formulated for families indexed by all local quadratic spaces of discriminant $\delta$, and why the fundamental lemma (\Cref{thm:FL}) is a statement about all local quadratic spaces simultaneously: if $W_z$ has no selfdual lattice, then $I^\circ(z)=0$ by convention, and \Cref{thm:FL} asserts $J^\circ(z)=0$.
		\item \emph{Relation to the connected distribution.} Throughout the comparison, $I_W$ denotes the distribution for the extended relative integration group $\mathcal H_W$, not the usual double $[\SO(W)]$-period $I_W^{\SO}$ with product measures. With the component measures fixed in \S\ref{sec:unfolding}, $I_W$ is obtained from the connected distribution by adelic component averaging. In particular, if each local test function is invariant under simultaneous conjugation by an isometry of $W_v$ of determinant $-1$, then
		\[I_W^{\SO,\rm reg}(f_W)=2I_W^{\rm reg}(f_W).\]
		\item \emph{Only finitely many $W$ contribute.} Let $\Sigma$ be a finite set of places containing the archimedean places and all places excluded by our unramified hypotheses (in particular, the dyadic places), such that outside $\Sigma$ the data are basic and $\delta$ is a unit. For $v\notin\Sigma$ there are at most two isometry classes of quadratic spaces over $k_v$ of dimension $2n$ and discriminant $\delta$, and exactly one of them admits a selfdual lattice (\Cref{lem:selfdual}); on the other the basic function is zero. Hence only spaces $W$ with $W_v$ of the first kind for all $v\notin\Sigma$ contribute. There are finitely many possibilities for $(W_v)_{v\in\Sigma}$, and by the Hasse--Minkowski theorem a global quadratic space is determined up to isometry by its localizations. So the sum over $W$ is finite.
	\end{enumerate}
\end{remark}

\begin{remark}[Dependence on the section]\label{rem:section}
The normalized orbital integral depends on the choice of the BLS section, through an orbit function. Suppose the elements $s(m_C)$ in \eqref{eq:Lz} are replaced by $s'(m_C)=\beta(C)s(m_C)$ for a function $\beta:\GL_n(F)\to\{\pm1\}$, and $\epsilon$ is replaced by the sign $\epsilon'$ defined by $s'(m_C)s(d_h)=\epsilon'(C,h)s(d_h)s'(m_{C^h})$ (so that \Cref{lem:covariance} continues to hold), namely $\epsilon'(C,h)=\beta(C)\beta(C^h)\epsilon(C,h)$. Put $b(z)=\beta(A)\beta(B)$. Then $K'(z)=b(z)K(z)$, $\mu'(z,h)=b(z)b(z\cdot h)\mu(z,h)$, hence $O'_z=b(z)O_z$, $\Delta'(z)=b(z)b(z_c)\Delta(z)$, and
\[J'_z=b(z_c)\,J_z .\]
The orthogonal orbital integrals involve no section. Consequently \Cref{thm:FL} is a statement relative to the BLS section fixed in \S\ref{sec:cover}; for another section it would hold only after multiplying $J^\circ$ by the orbit function $b(z_c)$.
\end{remark}

\subsection*{The basic functions}

Assume $\psi$ has conductor $\cO$, $q>3$ and all data are unramified. The linear basic data are
\[f^\circ=\one_{\GL_{2n}(\cO)},\qquad\omega^\circ=\xi^{+,\circ}\otimes\xi^{-,\circ},\qquad\Phi^\circ=\one_{\cO^{1\times2n}}.\]
For each quadratic space $W$ of dimension $2n$ we put $f^\circ_W=\one_{\SO(L_W\oplus\cO e)}$ if $W$ has a selfdual lattice $L_W$ and $f^\circ_W=0$ otherwise. This is the function $f_W$ of \Cref{prop:orth-unfold} attached to $\one_{\SO(L_W\oplus\cO e)}\otimes\one_{\SO(L_W)}$, since $\vol\SO(L_W)=1$.

\begin{lemma}\label{lem:basic}
\begin{enumerate}[label=(\alph*)]
\item For $z\in X^{\rm reg}(F)$ let $J^\circ(z)=J_z(f^\circ,\omega^\circ,\Phi^\circ)$. Then
\[J^\circ(z)=\sum_{\Lambda}w(\Lambda),\qquad w(h\cO^n)=\Delta(z)^{-1}\mu(z,h)\cW(A^h)\cW(B^h),\]
where $\Lambda$ runs over the lattices $\Lambda\subset E$ with $A\Lambda\subset\Lambda$, $v\in\Lambda$, $r\Lambda\subset\cO$. The summand depends only on $\Lambda=h\cO^n$, and the sum is finite.
\item Let $I^\circ(z)=I_z(f^\circ)$. Then
\[I^\circ(z)=\#\{L\subset W_z:\ L=L^\#,\ g_z(L\oplus\cO e)=L\oplus\cO e\},\]
where $L^\#=\{x\in W_z:(x,L)\subset\cO\}$. The set is finite.
\end{enumerate}
\end{lemma}
\begin{proof}
(a) For $g\in\GL_{2n}(\cO)$ take the lift $\kappa(g)$; then $\theta^\pm(\kappa(g))\xi^{\pm,\circ}=\xi^{\pm,\circ}$, $|\det g|=1$ and $\Phi^\circ((r,x)g)=\Phi^\circ(r,x)$, so $\widehat{\Phi^\circ}_g(r,v)=\one_{\cO^{1\times n}}(r)\one_{\cO^n}(v)$ (the selfdual Fourier transform of $\one_{\cO^{n}}$ is $\one_{\cO^n}$). Thus $L_z(g)$ is constant on $\GL_{2n}(\cO)$, and
\[K(z)=|\det A\det B|^{\frac{2n-1}4}\,\cW(A)\cW^-(B)\,\one_{\cO^n}(v)\one_{\cO^{1\times n}}(r).\]
By \Cref{lem:compact}, $\cW^-=\cW$. Since $|\det A^h\det B^h|=|\det A\det B|$, inserting $K(z\cdot h)$ into \eqref{eq:rawOI} and \eqref{eq:Jz} gives $J^\circ(z)=\Delta(z)^{-1}\int\mu(z,h)\cW(A^h)\cW(B^h)\one(h^{-1}v)\one(rh)\,dh$. By \Cref{lem:compact} the function $h\mapsto\epsilon(C,h)\cW(C^h)$ is right $\GL_n(\cO)$-invariant, and so are the indicators. Hence the integral is a sum over cosets $h\GL_n(\cO)$, i.e.\ over lattices $\Lambda=h\cO^n$, each of volume one. The integrand is zero unless $A^h$ is integral, $h^{-1}v\in\cO^n$ and $rh\in\cO^{1\times n}$, i.e.\ $A\Lambda\subset\Lambda$, $v\in\Lambda$, $r\Lambda\subset\cO$. These conditions imply $C_v\cO^n\subset\Lambda$ and $R_r\Lambda\subset\cO^n$, so $\Lambda$ lies in the finite interval $C_v\cO^n\subset\Lambda\subset R_r^{-1}\cO^n$.

(b) If $W_z$ has no selfdual lattice both sides are zero. Otherwise, $h^{-1}g_zh\in\SO(L_W\oplus\cO e)$ iff $g_z$ preserves $hL_W\oplus\cO e$. The cosets $h\Ort(L_W)$ correspond bijectively to the selfdual lattices $hL_W$, because $\Ort(W)$ acts transitively on selfdual lattices (\Cref{lem:selfdual}). Each coset has volume one. Finally, $e$ is cyclic for $g_z$ (\Cref{prop:matching}(a)), so $M_z=\sum_{i=0}^{2n}\cO g_z^ie$ is a lattice in $V_z$. Every counted $L$ gives a $g_z$-stable lattice $N=L\oplus\cO e$ containing $e$, hence $M_z\subset N$; and $N$ is selfdual, hence $N\subset M_z^\#$. There are only finitely many lattices between $M_z$ and $M_z^\#$.
\end{proof}

\begin{theorem}[Unit fundamental lemma]\label{thm:FL}
Assume $F$ has characteristic zero, odd residual characteristic and $q>3$, and that $\psi$ has conductor $\cO$. Then for every $z\in X^{\rm reg}(F)$,
\[J^\circ(z)=I^\circ(z).\]
\end{theorem}

The proof occupies Chapters \ref{ch:decomposition}--\ref{ch:fourier}. In \Cref{ch:decomposition} both sides are factored into an even and an odd Lie algebra orbital integral; in \Cref{ch:induction} the Lie algebra identities are proved by induction, assuming the Fourier compatibility \Cref{thm:fourier}, which is proved in \Cref{ch:fourier}.

\subsection*{Smooth transfer}

Fix $\delta\in F^\times/F^{\times2}$ and let $X^{\rm reg}_\delta(F)$ be the set of $z$ with $\det A\det B\in\delta F^{\times2}$. Let $\Omega_\delta$ be the $\eta_\delta$-isotypic part of $\Omega$ under the center, $\eta_\delta=\hil{\cdot}{\delta}$; the functional $q_z(\xi^+\otimes\xi^-)=\ell(\theta^+(s(m_A))\xi^+)\,\ell^-(\theta^-(s(m_B))\xi^-)$ attached to $z$ satisfies $q_z\circ\Omega(aI_{2n})=\hil{a}{\det A\det B}q_z$ by \ref{C:center} and $\beta\beta^-=1$, so $q_z$ vanishes on $\Omega_\delta$ unless $z\in X^{\rm reg}_\delta(F)$.

\begin{conjecture}[Smooth transfer]\label{conj:transfer}
For every finite sum $t$ of triples $(f,\omega,\Phi)$ with $f\in C_c^\infty(\GL_{2n}(F))$, $\omega\in\Omega_\delta$, $\Phi\in\cS(F^{1\times2n})$ there is a family $(f_W)$, $f_W\in C_c^\infty(\SO(W\oplus Fe))$, $W$ running over the quadratic spaces of dimension $2n$ and discriminant $\delta$, such that
\[J_z(t)=I_z(f)\qquad\text{for all }z\in X^{\rm reg}_\delta(F),\]
and conversely every such family is matched by some $t$.
\end{conjecture}

We do not address the conjecture in this paper. 

\chapter{The even and odd factors}\label{ch:decomposition}

In this chapter we reduce \Cref{thm:FL} to two Lie algebra statements. Fix $z=(A,v,r)\in X^{\rm reg}(F)$ and recall $B=A+vr-I_n$. If the invariants of $z$ are not integral, both $I^\circ(z)$ and $J^\circ(z)$ vanish. Otherwise we split $E=F^n$ according to the residual spectrum of $A$: the part where $A$ is residually nilpotent, of dimension $k$, and the part where $A$ is residually invertible, of dimension $n-k$. On the first part we use the Cayley transform of $-g_z$, on the second that of $g_z$ (we refer the reader to \Cref{sec:I} for the definition of these transforms); both have unit denominators. This produces an \emph{even} Lie datum of size $k$ (a quadratic space of dimension $2k$) and an \emph{odd} Lie datum of size $n-k+1$ (a quadratic space of dimension $2(n-k)+1$), and we prove
\[I^\circ(z)=I^e\cdot I^o,\qquad J^\circ(z)=J^e\cdot J^o .\]
The orthogonal factorization is a statement about unweighted lattice counts and is proved in \S\ref{sec:I}. The linear factorization (\S\ref{sec:J}) also has to keep track of the exceptional weights and of their signs.

Throughout, $F$ is as in \Cref{thm:FL}. A lattice in a finite dimensional $F$-vector space is a finitely generated $\cO$-submodule which spans it. We use repeatedly the following obvious fact: if $Y$ preserves a lattice $N$ and $\det Y\in\cO^\times$, then $YN=N$, so $Y^{-1}$ preserves $N$ as well; and if $YN\subset N$ then the characteristic polynomial of $Y$ lies in $\cO[T]$.

\section{The Lie models and the residual decomposition}\label{sec:models}

\subsection*{The two Lie models}

For a quadratic space $(W,Q)$, a skewadjoint operator $C$ and a vector $w$ we put
\[I_Q(C,w)=\#\{L\subset W:\ L=L^\#,\ CL\subset L,\ w\in L\}.\]
If $w$ is cyclic for $C$, every lattice counted contains $M=\sum_{i<\dim W}\cO C^iw$ and is contained in $M^\#$, so the count is finite.

\begin{definition}[Even model]\label{def:even}
An \emph{even datum of size $n$} is a triple $(S,x,r)$ with $S\in\GL_n(F)$, $x\in F^n$, $r\in F^{1\times n}$, whose Hankel matrix $(rS^{a+b}x)_{a,b<n}$ is invertible. Let $H$ be the symmetric matrix of \Cref{lem:Q0}, $S^{\mathsf t}H=HS$, $Hx=r^{\mathsf t}$, and on $F^n\oplus F^n$ put
\[Q^e=\diag(H,-HS),\qquad C^e=\begin{pmatrix}0&-S\\-I&0\end{pmatrix},\qquad w^e=\begin{pmatrix}x\\0\end{pmatrix}.\]
$Q^eC^e=\left(\begin{smallmatrix}0&-HS\\HS&0\end{smallmatrix}\right)$ is antisymmetric since $HS$ is symmetric, so $C^e$ is skewadjoint. The attached invariant $\nu$ is $\nu=rx=Q^e(w^e,w^e)$. Put $\Delta(S,x)=\epsilon(S,C_x)$ with $C_x=[x,Sx,\dots,S^{n-1}x]$ and
\[J^e(S,x,r)=\Delta(S,x)^{-1}\int_{\GL_n(F)}\epsilon(S,h)\cW(S^h)\one_{\cO^n}(h^{-1}x)\one_{\cO^{1\times n}}(rh)\,dh,\]\[I^e(S,x,r)=I_{Q^e}(C^e,w^e).\]
For $n=0$ both are $1$.
\end{definition}

\begin{definition}[Odd model]\label{def:odd}
An \emph{odd datum of size $n$} is a matrix
\[S=\begin{pmatrix}S_0&b\\c&d\end{pmatrix}\in\GL_n(F),\qquad S_0\in\Mat_{n-1}(F),\]
such that the Hankel matrix $(cS_0^{a+b}\,b)_{a,b<n-1}$ is invertible; $S_0$ may be singular. Let $H$ be the symmetric matrix of \Cref{lem:Q0} for $(S_0,b,c)$, and on $F^{n-1}\oplus F^{n-1}\oplus F$ put
\[Q^o=\begin{pmatrix}H&0&0\\0&-HS_0&-Hb\\0&-c&-d\end{pmatrix},\qquad C^o=\begin{pmatrix}0&-S_0&-b\\-I&0&0\\0&0&0\end{pmatrix},\qquad w^o=\begin{pmatrix}0\\0\\1\end{pmatrix}.\]
$Q^o$ is symmetric since $(Hb)^{\mathsf t}=c$. The invariant is $\nu=d=-Q^o(w^o,w^o)$. Let $f=e_n\in F^n$, $C_f=[f,Sf,\dots,S^{n-1}f]$, $\Delta(S)=\epsilon(S,C_f)$, and for $h\in\GL_{n-1}(F)$ put $\hat h=\diag(h,1)$. Define
\[J^o(S)=\Delta(S)^{-1}\int_{\GL_{n-1}(F)}\epsilon(S,\hat h)\cW(S^{\hat h})\,dh,\qquad I^o(S)=I_{Q^o}(C^o,w^o).\]
For $n=1$: $S=(d)$, $Q^o=(-d)$, $C^o=0$, $w^o=1$ and $J^o(S)=\cW(d)$.
\end{definition}

\begin{lemma}\label{lem:models-basic}
\begin{enumerate}[label=(\alph*)]
\item In the odd model, $\det Q^o=\pm\det(H)^2\det S$, so $Q^o$ is nondegenerate; $Q^oC^o$ is antisymmetric; $f$ is cyclic for $S$ and $w^o$ is cyclic for $C^o$. In the even model $w^e$ is cyclic for $C^e$.
\item (Invariance.) For $g\in\GL_n(F)$, $J^e(S^g,g^{-1}x,rg)=J^e(S,x,r)$ and $I^e(S^g,g^{-1}x,rg)=I^e(S,x,r)$. For $g\in\GL_{n-1}(F)$, $J^o(S^{\hat g})=J^o(S)$ and $I^o(S^{\hat g})=I^o(S)$.
\item (Lattice form.) $J^e(S,x,r)=\sum_\Lambda\Delta(S,x)^{-1}\epsilon(S,h)\cW(S^h)$ over lattices $\Lambda=h\cO^n$ with $S\Lambda\subset\Lambda$, $x\in\Lambda$, $r\Lambda\subset\cO$; and $J^o(S)=\sum_\Lambda\Delta(S)^{-1}\epsilon(S,\hat h)\cW(S^{\hat h})$ over lattices $\Lambda=h\cO^{n-1}$ with $S(\Lambda\oplus\cO f)\subset\Lambda\oplus\cO f$. Each summand depends only on $\Lambda$ and both sums are finite.
\end{enumerate}
\end{lemma}
\begin{proof}
(a) Eliminating the last column of $Q^o$ against the middle block one finds $\det Q^o=\det H\cdot\det(-HS_0)\cdot(-d+cS_0^{-1}b)$ when $S_0$ is invertible, and $\det S=\det S_0(d-cS_0^{-1}b)$; in general the identity $\det Q^o=(-1)^{n}\det(H)^2\det S$ follows by continuity in $S_0$ (it is polynomial). The product $Q^oC^o$ has rows $(0,-HS_0,-Hb)$, $(HS_0,0,0)$, $(c,0,0)$, which is antisymmetric because $HS_0$ is symmetric and $c^{\mathsf t}=Hb$. For cyclicity of $f$: the projection of $S^{a+1}f$ to $F^{n-1}$ is $S_0^ab$ plus an $F$-combination of $S_0^ib$, $i<a$; since $b$ is cyclic for $S_0$, the vectors $f,Sf,\dots,S^{n-1}f$ span $F^n$. For $w^o$: $C^ow^o=(-b,0,0)$, $(C^o)^2(X,Y,z)=(S_0X,S_0Y+bz,0)$, so the $C^o$-span of $w^o$ contains $(S_0^ab,0,0)$ and $(0,S_0^ab,0)$ for all $a$, hence everything. The even case is similar ($(C^e)^2=\diag(S,S)$).

(b) If $(S,x,r)$ is replaced by $(S^g,g^{-1}x,rg)$, then $C_{g^{-1}x}(S^g)=g^{-1}C_x$; by the cocycle identity $\epsilon(S^g,h)=\epsilon(S,g)^{-1}\epsilon(S,gh)$ and $\Delta(S^g,g^{-1}x)=\epsilon(S,g)^{-1}\Delta(S,x)$, and the substitution $h\mapsto g^{-1}h$ proves the invariance of $J^e$. The same computation with $\hat g$ and $C_f(S^{\hat g})=\hat g^{-1}C_f(S)$ proves it for $J^o$. On the orthogonal side, by \Cref{lem:Q0} the new form is $g^{\mathsf t}Hg$, and $\diag(g,g)$ (resp. $\diag(g,g,1)$) is an isometry carrying the new model to the old one.

(c) This is the argument of \Cref{lem:basic}(a), using the right $\GL(\cO)$-invariance of $h\mapsto\epsilon(S,h)\cW(S^h)$ (\Cref{lem:compact}). In the odd case, $\cW(S^{\hat h})\ne0$ forces $S$ to preserve $\Lambda\oplus\cO f$; then $S_0\Lambda\subset\Lambda$, $b\in\Lambda$, $c\Lambda\subset\cO$, so $C_b(S_0)\cO^{n-1}\subset\Lambda\subset R_c(S_0)^{-1}\cO^{n-1}$, a finite interval.
\end{proof}

\subsection*{The integral locus}

\begin{definition}
$z=(A,v,r)\in X^{\rm reg}(F)$ is \emph{integral} if the characteristic polynomial $P_A$ of $A$ lies in $\cO[T]$ and $rA^iv\in\cO$ for all $i\ge0$.
\end{definition}

If $z$ is integral then $M=\cO[A]v$ is an $A$-stable lattice with basis $v,Av,\dots,A^{n-1}v$ (because $P_A$ is monic integral), and its Gram matrix for $Q_0$ is $(rA^{a+b}v)$, which is integral; so $M\subset M^\#=\{x:Q_0(M,x)\subset\cO\}$. Since $r(x)=Q_0(v,x)$ we have $r(M)\subset\cO$.

\begin{lemma}\label{lem:integral}
If $J^\circ(z)\ne0$ or $I^\circ(z)\ne0$, then $z$ is integral. More precisely:
\begin{enumerate}[label=(\alph*)]
\item every lattice $\Lambda$ occurring in \Cref{lem:basic}(a) satisfies $B\Lambda\subset\Lambda$ and $M\subset\Lambda\subset M^\#$;
\item every lattice $L$ occurring in \Cref{lem:basic}(b) satisfies $M\oplus M\oplus\cO e\subset L\oplus\cO e$.
\end{enumerate}
\end{lemma}
\begin{proof}
(a) $Bx=Ax+v\,r(x)-x\in\Lambda$ for $x\in\Lambda$. From $v\in\Lambda$ and $A\Lambda\subset\Lambda$ we get $M\subset\Lambda$; and for $x\in\Lambda$, $Q_0(A^iv,x)=rA^ix\in\cO$, so $x\in M^\#$. The matrix of $A$ in a basis of $\Lambda$ is integral, so $P_A\in\cO[T]$; and $rA^iv\in r(\Lambda)\subset\cO$.

(b) Put $N=L\oplus\cO e$, so $g_zN=N$ and $g_z^{-1}N=N$. Let $s_0=\diag(-I,I,1)$; from the proof of \Cref{prop:matching}, $g_z=s_0s_1$ with $s_0,s_1$ involutions, so $g_z^{-1}=s_1s_0=s_0g_zs_0$. Hence, by \eqref{eq:gz},
\[\tfrac14(g_z-g_z^{-1})e=(v,0,0)\in N,\qquad \tau=\tfrac14(g_z+g_z^{-1}+2I)=\begin{pmatrix}A&0&0\\0&I+B&-v\\0&rB&1-rv\end{pmatrix}\]
preserves $N$ (as $2\in\cO^\times$). Since $N$ is selfdual and $e,g_ze\in N$, $(e,g_ze)=1-2rv\in\cO$, so $rv\in\cO$; then $g_ze-(1-2rv)e-2(v,0,0)=(0,-2v,0)\in N$, so $(0,v,0)\in N$. Applying $\tau$ repeatedly to $(v,0,0)$ gives $(A^iv,0,0)\in N$, so $M\oplus0\oplus0\subset N$, and pairing with $(v,0,0)$ gives $rA^{i+1}v=(A^iv,0,0)\cdot(v,0,0)\in\cO$. If $(0,y,0)\in N$ then $\tau(0,y,0)=(0,(I+B)y,rBy)$, whose last coordinate is integral (pair with $e$); subtracting $(rBy)e$ gives $(0,(I+B)y,0)\in N$. Hence $0\oplus\cO[B]v\oplus0\subset N$. Finally $\cO[B]v=M$: $B^jv=A^jv+\sum_{i<j}\beta_{ji}A^iv$ with $\beta_{ji}\in\cO$, because $B(A^iv)=A^{i+1}v-A^iv+v\cdot rA^iv$ and the moments are integral; the change of generators is unitriangular. So $M$ is contained in the lattice $N\cap E$, $M$ is finitely generated and $A$-stable, and $P_A\in\cO[T]$.
\end{proof}

\subsection*{The residual decomposition}

From now on $z$ is integral. Reduce $P_A$ modulo $\varpi$ and write $\bar P_A(T)=T^k\bar P_u(T)$ with $\bar P_u(0)\ne0$. By Hensel's lemma $P_A=P_0P_u$ with monic $P_0,P_u\in\cO[T]$ reducing to $T^k$ and $\bar P_u$. Their reductions are coprime, so the resultant of $P_0,P_u$ is a unit, and there are $\phi_0,\phi_u\in\cO[T]$ with $\phi_0P_0+\phi_uP_u=1$. Put
\[\pi_0=\phi_u(A)P_u(A),\qquad\pi_u=\phi_0(A)P_0(A).\]
Since $P_0(A)P_u(A)=0$, these are complementary idempotents in $\cO[A]$. They preserve every $A$-stable lattice. Being polynomials in the $Q_0$-selfadjoint operator $A$, they are selfadjoint, so $E=E_0\oplus E_u$ ($E_0=\pi_0E$, $E_u=\pi_uE$) is $Q_0$-orthogonal. We choose a basis adapted to this decomposition; by \Cref{lem:models-basic}(b) and its analogue for $J^\circ,I^\circ$ (the definitions are invariant under $z\mapsto z\cdot h$) this is harmless. Then
\[A=\diag(A_0,A_u),\quad Q_0=\diag(Q_{00},Q_{0u}),\quad v=\begin{pmatrix}v_0\\v_u\end{pmatrix},\quad r=(r_0,r_u),\quad M=M_0\oplus M_u,\]
where $r_0,r_u$ are the restrictions of $r$ (since $r(x)=Q_0(v,x)$, $Q_{00}v_0=r_0^{\mathsf t}$ and $Q_{0u}v_u=r_u^{\mathsf t}$), and $M_i=\cO[A_i]v_i$ has basis $v_i,A_iv_i,\dots$; in particular $v_0$ is cyclic for $A_0$ and $v_u$ for $A_u$. We have:
\begin{enumerate}[label=(\roman*)]
\item $A_0$ is topologically nilpotent: $P_0(A_0)=0$ with $P_0\equiv T^k$, so $A_0^k\Lambda_0\subset\varpi\Lambda_0$ for every $A_0$-stable lattice $\Lambda_0$.
\item $A_0-I$ has unit determinant ($\det(A_0-I)=(-1)^kP_0(1)\equiv(-1)^k$), so $(A_0-I)^{\pm1}$ preserve every $A_0$-stable lattice.
\item $A_u$ has unit determinant ($\det A_u=\pm P_u(0)$), so $A_u^{\pm1}$ preserve every $A_u$-stable lattice.
\end{enumerate}

\begin{definition}[The even and odd factors]\label{def:factors}
The \emph{even factor} is the triple
\[S_e=A_0(A_0-I)^{-1}=I+(A_0-I)^{-1},\qquad x_e=(I-A_0)^{-1}v_0,\qquad r_0,\]
of size $k$, and the \emph{odd factor} is
\[S_o=\begin{pmatrix}I-A_u^{-1}&-A_u^{-1}v_u\\-r_u&r_0x_e-1\end{pmatrix}\in\GL_{n-k+1}(F).\]
If $k=0$ we put $I^e=J^e=1$ (and $x_e=0$, so that the corner entry of $S_o$ is $-1$); if $k=n$, $S_o=(r_0x_e-1)$.
\end{definition}

\begin{theorem}[Even--odd factorization]\label{thm:factorization}
Let $z\in X^{\rm reg}(F)$. If $z$ is not integral, $I^\circ(z)=J^\circ(z)=0$. If $z$ is integral, then with the data of \Cref{def:factors},
\[I^\circ(z)=I^e(S_e,x_e,r_0)\,I^o(S_o),\qquad J^\circ(z)=J^e(S_e,x_e,r_0)\,J^o(S_o).\]
\end{theorem}

Here $(S_e,x_e,r_0)$ is an even datum and $S_o$ is an odd datum (\Cref{lem:factor-identities}(b)--(d) below), so the right-hand sides are defined.

The first assertion in the Theorem is \Cref{lem:integral}. The orthogonal identity is proved in \S\ref{sec:I} and the linear one in \S\ref{sec:J}.

The two proofs use the following identities.

\begin{lemma}\label{lem:factor-identities}
\begin{enumerate}[label=(\alph*)]
\item With $T_u=\left(\begin{smallmatrix}A_u-I&v_u\\-r_u&1-r_0x_e\end{smallmatrix}\right)$ and blocks of sizes $(k,n-k,1)$,
\[L_0\begin{pmatrix}I_n&v\\0&1\end{pmatrix}\begin{pmatrix}B&0\\0&1\end{pmatrix}\begin{pmatrix}I_n&0\\-r&1\end{pmatrix}R_0=\diag(A_0-I,T_u),\]
where, in blocks of sizes $(k,n-k,1)$,
\[L_0=\begin{pmatrix}I_k&0&0\\0&I_{n-k}&0\\r_0(A_0-I)^{-1}&0&1\end{pmatrix},\qquad R_0=\begin{pmatrix}I_k&0&-(A_0-I)^{-1}v_0\\0&I_{n-k}&0\\0&0&1\end{pmatrix}.\] Moreover $S_o=\diag(A_u^{-1},1)\,T_u\,R$ with $R=\diag(I_{n-k},-1)$.
\item $\det B=\det(A_0-I)\det T_u$ and $\det S_o=-\det B/(\det(A_0-I)\det A_u)$; in particular $S_o$ is invertible, $\val\det S_e=\val\det A$ and $\val\det S_o=\val\det B$.
\item $(S_e,x_e,r_0)$ is an even datum, its form is $H_e=Q_{00}(I-A_0)$, and $-H_eS_e=Q_{00}A_0$.
\item $S_o$ is an odd datum; with $S_o=\left(\begin{smallmatrix}S'&b\\c&d\end{smallmatrix}\right)$, i.e.\ $S'=I-A_u^{-1}$, $b=-A_u^{-1}v_u$, $c=-r_u$, $d=r_0x_e-1$, its form is $H_o=Q_{0u}A_u$, and $-H_oS'=Q_{0u}(I-A_u)$, $-H_ob=r_u^{\mathsf t}$.
\end{enumerate}
\end{lemma}
\begin{proof}
(a) First, $\left(\begin{smallmatrix}I&v\\0&1\end{smallmatrix}\right)\diag(B,1)\left(\begin{smallmatrix}I&0\\-r&1\end{smallmatrix}\right)=\left(\begin{smallmatrix}B-vr&v\\-r&1\end{smallmatrix}\right)=\left(\begin{smallmatrix}A-I&v\\-r&1\end{smallmatrix}\right)$. In blocks $(k,n-k,1)$ this is $N=\left(\begin{smallmatrix}A_0-I&0&v_0\\0&A_u-I&v_u\\-r_0&-r_u&1\end{smallmatrix}\right)$. Right multiplication by $R_0$ replaces the third block column by itself minus the first times $(A_0-I)^{-1}v_0$, giving third column $(0,v_u,1+r_0(A_0-I)^{-1}v_0)=(0,v_u,1-r_0x_e)$; left multiplication by $L_0$ adds $r_0(A_0-I)^{-1}$ times the first block row to the third, giving third row $(0,-r_u,1-r_0x_e)$. The formula for $S_o$ is a direct multiplication.

(b) Take determinants in (a): all factors other than $\diag(B,1)$ and $\diag(A_0-I,T_u)$ are unipotent, $\det R=-1$. Also $\det S_e=\det A_0/\det(A_0-I)$ and $\det A=\det A_0\det A_u$; $\det(A_0-I),\det A_u$ are units.

(c) $H_ex_e=Q_{00}(I-A_0)(I-A_0)^{-1}v_0=Q_{00}v_0=r_0^{\mathsf t}$, and $H_eS_e=-Q_{00}(I-A_0)A_0(I-A_0)^{-1}=-Q_{00}A_0$ is symmetric; by uniqueness in \Cref{lem:Q0} the form is $H_e$, provided the Hankel matrix is invertible. This holds because $F[S_e]=F[A_0]$, so $x_e=-(A_0-I)^{-1}v_0$ is cyclic for $S_e$, and $H_e$ is nondegenerate (the Hankel matrix is $C_{x_e}^{\mathsf t}H_eC_{x_e}$).

(d) $H_ob=-Q_{0u}v_u=-r_u^{\mathsf t}=c^{\mathsf t}$, and $H_oS'=Q_{0u}(A_u-I)$ is symmetric. As in (c), $b$ is cyclic for $S'$ since $F[S']=F[A_u]$, and $H_o$ is nondegenerate.
\end{proof}

\section{The factorization of $I$}\label{sec:I}

In this section $N$ denotes a selfdual lattice in $V_z$ with $g_zN=N$ and $e\in N$. By the proof of \Cref{lem:basic}(b) and the following remark, $I^\circ(z)$ is the number of such $N$: if $L$ is selfdual in $W_z$, so is $L\oplus\cO e$; conversely, if $N$ is selfdual and $e\in N$, then $x-(x,e)e\in N\cap W_z$ for $x\in N$, so $N=(N\cap W_z)\oplus\cO e$, and $N\cap W_z$ is selfdual.

The proof of the orthogonal half of \Cref{thm:factorization} has three steps.
\begin{enumerate}[label=(\arabic*)]
\item \emph{Splitting the quadratic space} (\Cref{lem:Vsplit}). An explicit change of coordinates, carried out by matrix multiplication, shows that $V_z$ is the orthogonal direct sum of two $g_z$-stable subspaces $V_-$ and $V_u$, and computes the restrictions $q_\pm$, $g_\pm$ and the components $e_\pm$ of $e$. This is a statement about $F$-vector spaces.
\item \emph{Splitting the lattices} (\Cref{lem:charpoly,lem:Nsplit}). The characteristic polynomials of $g_-$ and $g_u$ are coprime modulo $\varpi$. Hence the projections onto $V_-$ and $V_u$ are polynomials in $g_z$ with coefficients in $\cO$, and every lattice $N$ as above splits along $V_z=V_-\oplus V_u$. This factors $I^\circ(z)$ into two lattice counts for the orthogonal group elements $g_-$ and $g_u$ (\eqref{eq:I-group-split}).
\item \emph{Cayley transforms} (\Cref{lem:cayley,lem:Cminus,lem:Cu}). The Cayley transforms of $-g_-$ and of $g_u$ have unit denominators. They replace the two group-level counts by counts for skewadjoint operators, which are those of the even and odd Lie models of \Cref{def:even,def:odd}.
\end{enumerate}

\subsection*{Step 1: splitting the quadratic space}

We write vectors of $V_z=E\oplus E\oplus F$ as columns $(x,y,t)$. In these coordinates the Gram matrix of $q$ and the matrix of $g_z$ are
\[
Q=\begin{pmatrix}Q_0A&0&0\\0&-Q_0B&0\\0&0&1\end{pmatrix},\qquad
g_z=\begin{pmatrix}2A-I&-2B&2v\\-2A&I+2B&-2v\\-2rA&2rB&1-2rv\end{pmatrix}.
\]
For an invertible matrix $P$, the substitution $(\text{old coordinates})=P\cdot(\text{new coordinates})$ replaces them by $P^{\mathsf t}QP$ and $P^{-1}g_zP$. The matrices $P$ below have entries in $F$; no integrality is claimed or needed, since all statements about lattices will be made through projectors in $\cO[g_z]$ (Step 2).

\paragraph{First change of coordinates.} Put $p=x+y$ and $s=t-ry$, i.e.\ $x=p-y$ and $t=s+ry$:
\[
P_1=\begin{pmatrix}I&-I&0\\0&I&0\\0&r&1\end{pmatrix},\qquad P_1^{-1}=\begin{pmatrix}I&I&0\\0&I&0\\0&-r&1\end{pmatrix}.
\]
Using $A^{\mathsf t}Q_0=Q_0A$, $Q_0v=r^{\mathsf t}$ and $B=A+vr-I$ one obtains, in the coordinates $(p,y,s)$,
\begin{equation}\label{eq:q-pys}
P_1^{\mathsf t}QP_1=\begin{pmatrix}Q_0A&-Q_0A&0\\-Q_0A&Q_0&r^{\mathsf t}\\0&r&1\end{pmatrix},
\end{equation}
\begin{equation}\label{eq:g-pys}
P_1^{-1}g_zP_1=\begin{pmatrix}-I&2I&0\\-2A&4A-I&-2v\\0&0&1\end{pmatrix}.
\end{equation}
In particular $P_1^{-1}g_zP_1$ fixes the coordinate $s$.

\paragraph{Reordering along $E=E_0\oplus E_u$.} Write $A=\diag(A_0,A_u)$, $Q_0=\diag(Q_{00},Q_{0u})$, $v=\binom{v_0}{v_u}$, $r=(r_0,r_u)$, decompose $p=(p_0,p_u)$, $y=(y_0,y_u)$, and order the coordinates as $(p_0,y_0\mid p_u,y_u,s)$. Then \eqref{eq:q-pys} and \eqref{eq:g-pys} read, respectively,
\[
\left(\begin{array}{cc|ccc}
Q_{00}A_0&-Q_{00}A_0&0&0&0\\
-Q_{00}A_0&Q_{00}&0&0&r_0^{\mathsf t}\\\hline
0&0&Q_{0u}A_u&-Q_{0u}A_u&0\\
0&0&-Q_{0u}A_u&Q_{0u}&r_u^{\mathsf t}\\
0&r_0&0&r_u&1
\end{array}\right),\]
\[
\left(\begin{array}{cc|ccc}
-I&2I&0&0&0\\
-2A_0&4A_0-I&0&0&-2v_0\\\hline
0&0&-I&2I&0\\
0&0&-2A_u&4A_u-I&-2v_u\\
0&0&0&0&1
\end{array}\right).
\]
The $E_0$-block and the $E_u$-block are already decoupled. The only coupling left is that of the $E_0$-block with the coordinate $s$, through the entries $r_0^{\mathsf t}$ (in the form) and $-2v_0$ (in the action).

\paragraph{Second change of coordinates.} Put $a=(A_0-I)^{-1}v_0=-x_e$, and substitute $p_0=p_0'+as$, $y_0=y_0'+as$:
\[
P_2=\begin{pmatrix}I&0&0&0&a\\0&I&0&0&a\\0&0&I&0&0\\0&0&0&I&0\\0&0&0&0&1\end{pmatrix},
\]
and $P_2^{-1}$ is the same matrix with $a$ replaced by $-a$. Finally, on the $E_u$-part put $X=p_u-y_u$, i.e.\ $p_u=X+y_u$:
\[
P_3=\begin{pmatrix}I&I&0\\0&I&0\\0&0&1\end{pmatrix}\quad\text{acting on }(X,y_u,s).
\]

We get the following splitting of the space $V_z$. The proof is a matrix calculation which we do not record here.
\begin{lemma}\label{lem:Vsplit}
In the coordinates $(p_0',y_0')$ and $(X,y_u,s)$ obtained by $P_1$, the reordering, $P_2$ and $P_3$, the space $V_z$ is the orthogonal direct sum of $g_z$-stable subspaces $V_-$ of dimension $2k$ and $V_u$ of dimension $2(n-k)+1$, on which (writing $p,y$ for $p_0',y_0'$ on $V_-$ and $y$ for $y_u$ on $V_u$)
\begin{align*}
q_-&=\begin{pmatrix}Q_{00}A_0&-Q_{00}A_0\\-Q_{00}A_0&Q_{00}\end{pmatrix}, &
g_-&=\begin{pmatrix}-I&2I\\-2A_0&4A_0-I\end{pmatrix},\\
q_u&=\begin{pmatrix}Q_{0u}A_u&0&0\\0&Q_{0u}(I-A_u)&r_u^{\mathsf t}\\0&r_u&1-r_0x_e\end{pmatrix}, &
g_u&=\begin{pmatrix}2A_u-I&2(I-A_u)&2v_u\\-2A_u&2A_u-I&-2v_u\\0&0&1\end{pmatrix}.
\end{align*}

 The vector $e$ decomposes as $e=e_-+e_u$ with $e_-=(-a,-a)$ and $e_u=(0,0,1)$; $q_-(e_-)=-r_0a$, $q_u(e_u)=1-r_0x_e$. Here $a=(A_0-I)^{-1}v_0=-x_e$.
\end{lemma}

\subsection*{Step 2: splitting the lattices}

\Cref{lem:Vsplit} is a decomposition of $F$-vector spaces. What makes the lattices split according to this decomposition is that the projections onto $V_-$ and $V_u$ are polynomials in $g_z$ with coefficients in $\cO$; such polynomials preserve every lattice stable under $g_z$. By the B\'ezout identity this holds as soon as the characteristic polynomials of $g_-$ and $g_u$ are coprime modulo $\varpi$. The next lemma proves this coprimality; it reflects the separation of the residual spectrum of $g_z$ into the residual eigenvalue $-1$ (on $V_-$) and residual eigenvalues different from $-1$ (on $V_u$), which corresponds to the separation of the residual spectrum of $A$ into $0$ and nonzero values.

\begin{lemma}\label{lem:charpoly}
$\det(T-g_z)=(T-1)\det((T+1)^2I-4TA)$. The characteristic polynomials of $g_-$ and $g_u$ are
\[\chi_-(T)=\det((T+1)^2I-4TA_0),\qquad\chi_u(T)=(T-1)\det((T+1)^2I-4TA_u);\]
they lie in $\cO[T]$, $\chi_-\equiv(T+1)^{2k}\pmod\varpi$ and $\chi_u(-1)=-2\cdot4^{n-k}\det A_u\in\cO^\times$. In particular the reductions of $\chi_-$ and $\chi_u$ modulo $\varpi$ are coprime.
\end{lemma}
\begin{proof}
In the coordinates $(p,y,s)$, by \eqref{eq:g-pys}, $T-g_z$ is block upper triangular with diagonal blocks $\left(\begin{smallmatrix}(T+1)I&-2I\\2A&(T+1)I-4A\end{smallmatrix}\right)$ and $T-1$. As the upper left block is scalar, the determinant of the $2\times2$ block matrix is $\det((T+1)((T+1)I-4A)+4A)=\det((T+1)^2I-4TA)$. The same computation applies to $g_-$ and $g_u$, using \Cref{lem:Vsplit}. Recall the decomposition $P_A=P_0P_u$. Since $P_0$ is monic of degree $k$, $\chi_-(T)=(4T)^kP_0((T+1)^2/4T)$, which lies in $\cO[T]$ and reduces to $(T+1)^{2k}$ because $\bar P_0=T^k$; similarly for $\chi_u$ with $P_u$. Finally $\chi_u(-1)=-2\det(4A_u)$, a unit because $A_u$ is residually invertible; so $\bar\chi_u$ does not vanish at $-1$ and is coprime to $\bar\chi_-=(T+1)^{2k}$.
\end{proof}

By \Cref{lem:charpoly} there are $\psi_-,\psi_u\in\cO[T]$ with $\psi_-\chi_-+\psi_u\chi_u=1$. Put $\Pi_-=\psi_u(g_z)\chi_u(g_z)$ and $\Pi_u=\psi_-(g_z)\chi_-(g_z)$. On $V_-$, $\chi_-(g_-)=0$ by the Cayley--Hamilton theorem, so $\Pi_u=0$ and $\Pi_-=I$ there; on $V_u$ the roles are reversed. Thus $\Pi_-$ and $\Pi_u$ are exactly the projections onto $V_-$ and $V_u$ of \Cref{lem:Vsplit}, now written as elements of $\cO[g_z]$.

\begin{lemma}\label{lem:Nsplit}
$N\mapsto(N\cap V_-,N\cap V_u)$ is a bijection from the set of selfdual $g_z$-stable lattices $N\ni e$ onto the set of pairs $(N_-,N_u)$ of selfdual lattices $N_-\subset V_-$, $N_u\subset V_u$ with $g_-N_-=N_-$, $e_-\in N_-$, $g_uN_u=N_u$, $e_u\in N_u$. The inverse is $(N_-,N_u)\mapsto N_-\oplus N_u$.
\end{lemma}
\begin{proof}
$\Pi_\pm N\subset N$, so $N=\Pi_-N\oplus\Pi_uN$ and $\Pi_\pm N=N\cap V_\pm$. As $V=V_-\perp V_u$, $N$ is selfdual iff both summands are; $e_\pm=\Pi_\pm e\in N$; and $g_z$-stability of $N$ is equivalent to $g_\pm$-stability of $N\cap V_\pm$. The converse is clear.
\end{proof}

Consequently
\begin{equation}\label{eq:I-group-split}
\begin{aligned}
I^\circ(z)={}&\#\{N_-\subset V_-:\ N_-=N_-^\#,\ g_-N_-=N_-,\ e_-\in N_-\}\\
&\cdot\#\{N_u\subset V_u:\ N_u=N_u^\#,\ g_uN_u=N_u,\ e_u\in N_u\}.
\end{aligned}
\end{equation}
This factorization uses only Steps 1 and 2. The two factors are lattice counts for the orthogonal group elements $g_-$ and $g_u$.

\subsection*{Step 3: Cayley transforms with unit denominators}

The even and odd Lie models are defined through skewadjoint operators rather than orthogonal group elements. To identify the two factors of \eqref{eq:I-group-split} with them we use Cayley transforms. For a Cayley transform to preserve stability conditions on lattices its denominator must be a unit; this is the content of the next lemma.

\begin{lemma}\label{lem:cayley}
Let $h$ be an invertible operator on a finite dimensional space $U$ whose characteristic polynomial $\chi_h$ lies in $\cO[T]$ and satisfies $\chi_h(-1)\in\cO^\times$. Put $Y=(h-I)(h+I)^{-1}$. Then for every lattice $N\subset U$, $hN\subset N$ if and only if $YN\subset N$. If $h$ is an isometry, $Y$ is skewadjoint.
\end{lemma}
\begin{proof}
$\det(h+I)=\pm\chi_h(-1)$ is a unit and $I-Y=2(h+I)^{-1}$ has unit determinant. If $hN\subset N$ then $h+I$ preserves $N$ and, having unit determinant, so does $(h+I)^{-1}$; hence $Y=I-2(h+I)^{-1}$ preserves $N$. Conversely, if $YN\subset N$ then $I-Y$ and its inverse preserve $N$, so $h=-I+2(I-Y)^{-1}$ does. If $h$ is an isometry, write $x=(h+I)x'$, $y=(h+I)y'$. Expanding and using $(hx',hy')=(x',y')$,
\[(Yx,y)=((h-I)x',(h+I)y')=(hx',y')-(x',hy'),\]\[(x,Yy)=((h+I)x',(h-I)y')=(x',hy')-(hx',y'),\]
so $(Yx,y)=-(x,Yy)$.
\end{proof}

The unit denominators come from the same residual separation as in Step 2. By \Cref{lem:charpoly}, all residual eigenvalues of $g_-$ are $-1$, and $-1\ne1$ residually because the residual characteristic is odd; so $\chi_{-g_-}(-1)=\pm\chi_-(1)\equiv\pm2^{2k}$ is a unit, and \Cref{lem:cayley} applies to $-g_-$. Also $\chi_u(-1)$ is a unit, so \Cref{lem:cayley} applies to $g_u$. We therefore use
\[C_-=(g_-+I)(g_--I)^{-1}\ (\text{the Cayley transform of }-g_-),\qquad C_u=(g_u-I)(g_u+I)^{-1}.\]

\begin{lemma}\label{lem:Cminus}
In the coordinates $(X,Y)=(y,p-y)$ on $V_-$,
\[C_-=\begin{pmatrix}0&-S_e\\-I&0\end{pmatrix},\qquad q_-=\diag(H_e,-H_eS_e),\qquad e_-=\begin{pmatrix}x_e\\0\end{pmatrix}.\]
Thus $(V_-,q_-,C_-,e_-)$ is the even model of $(S_e,x_e,r_0)$.
\end{lemma}
\begin{proof}
From $p=X+Y$, $y=X$ and the formula for $g_-$: the image of $(X,Y)$ has $y'=-2A_0(X+Y)+(4A_0-I)X$ and $p'=-(X+Y)+2X$, so $X'=(2A_0-I)X-2A_0Y$ and $Y'=p'-y'=2(I-A_0)X+(2A_0-I)Y$. Therefore
\[g_--I=2\begin{pmatrix}A_0-I&-A_0\\-(A_0-I)&A_0-I\end{pmatrix},\qquad g_-+I=2\begin{pmatrix}A_0&-A_0\\-(A_0-I)&A_0\end{pmatrix},\]
and $\left(\begin{smallmatrix}0&-S_e\\-I&0\end{smallmatrix}\right)(g_--I)=2\left(\begin{smallmatrix}S_e(A_0-I)&-S_e(A_0-I)\\-(A_0-I)&A_0\end{smallmatrix}\right)=g_-+I$ because $S_e(A_0-I)=A_0$. Substituting $p=X+Y$, $y=X$ in $q_-$ the cross terms cancel and $q_-=X^{\mathsf t}Q_{00}(I-A_0)X+Y^{\mathsf t}Q_{00}A_0Y$, which is $\diag(H_e,-H_eS_e)$ by \Cref{lem:factor-identities}(c). Finally $e_-$ has $X=y_0'=-a=x_e$ and $Y=p_0'-y_0'=0$.
\end{proof}

\begin{lemma}\label{lem:Cu}
In the coordinates $(X,y,s)$ on $V_u$,
\[C_u=\begin{pmatrix}0&-(I-A_u^{-1})&A_u^{-1}v_u\\-I&0&0\\0&0&0\end{pmatrix},\]
and $(V_u,q_u,C_u,e_u)$ is the odd model of $S_o$.
\end{lemma}
\begin{proof}
We check $C_u(g_u+I)=g_u-I$. The second and third rows are immediate ($g_u$ fixes $s$, so the last row of $g_u-I$ vanishes). The first row of $C_u(g_u+I)$ is $-(I-A_u^{-1})(-2A_u,2A_u,-2v_u)+A_u^{-1}v_u(0,0,2)=(2(A_u-I),2(I-A_u),2v_u)$, the first row of $g_u-I$. With the blocks of \Cref{lem:factor-identities}(d), $-S'=-(I-A_u^{-1})$ and $-b=A_u^{-1}v_u$, so $C_u=C^o$. The Gram matrix of $q_u$ is $\left(\begin{smallmatrix}Q_{0u}A_u&0&0\\0&Q_{0u}(I-A_u)&r_u^{\mathsf t}\\0&r_u&1-r_0x_e\end{smallmatrix}\right)=\left(\begin{smallmatrix}H_o&0&0\\0&-H_oS'&-H_ob\\0&-c&-d\end{smallmatrix}\right)=Q^o$, and $e_u=w^o$.
\end{proof}

\begin{proof}[Proof of the orthogonal half of \Cref{thm:factorization}]
By \eqref{eq:I-group-split}, $I^\circ(z)$ is the product of the numbers of lattices $N_-$ and $N_u$ described there. By \Cref{lem:cayley}, applied to $-g_-$ and to $g_u$, $g_\pm$-stability may be replaced by $C_\pm$-stability. By \Cref{lem:Cminus,lem:Cu}, the first factor is $I^e(S_e,x_e,r_0)$ and the second is $I^o(S_o)$.
\end{proof}

\section{The factorization of $J$}\label{sec:J}

In this section $\Lambda$ denotes a lattice in $E$ with $A\Lambda\subset\Lambda$, $v\in\Lambda$, $r\Lambda\subset\cO$; we call it contributing. Recall from \Cref{lem:basic} and \Cref{lem:models-basic}(c) that
\[J^\circ(z)=\sum_\Lambda w(\Lambda),\qquad J^e(S_e,x_e,r_0)=\sum_{\Lambda_0}w_e(\Lambda_0),\qquad J^o(S_o)=\sum_{\Lambda_u}w_o(\Lambda_u),\]
where $w_e(h_0\cO^k)=\Delta(S_e,x_e)^{-1}\epsilon(S_e,h_0)\cW(S_e^{h_0})$ and $w_o(h_u\cO^{n-k})=\Delta(S_o)^{-1}\epsilon(S_o,\hat h_u)\cW(S_o^{\hat h_u})$.

\subsection*{The lattices}

\begin{lemma}\label{lem:Lsplit}
$\Lambda\mapsto(\Lambda_0,\Lambda_u)=(\pi_0\Lambda,\pi_u\Lambda)$ is a bijection from the contributing lattices onto the pairs of lattices $\Lambda_0\subset E_0$, $\Lambda_u\subset E_u$ such that
\[S_e\Lambda_0\subset\Lambda_0,\quad x_e\in\Lambda_0,\quad r_0\Lambda_0\subset\cO,\qquad S_o(\Lambda_u\oplus\cO f)\subset\Lambda_u\oplus\cO f .\]
The inverse is $(\Lambda_0,\Lambda_u)\mapsto\Lambda_0\oplus\Lambda_u$. Moreover $M\mapsto(M_0,M_u)$, and $\cO[S_e]x_e=M_0$, $\cO[S_o]f=M_u\oplus\cO f$.
\end{lemma}
\begin{proof}
Since $\pi_0,\pi_u\in\cO[A]$ preserve $\Lambda$, $\Lambda=\Lambda_0\oplus\Lambda_u$ with $\Lambda_i=\Lambda\cap E_i$. Suppose $\Lambda$ contributes. Then $\Lambda_0$ is $A_0$-stable, $v_0=\pi_0v\in\Lambda_0$ and $r_0\Lambda_0\subset r\Lambda\subset\cO$. By (ii), $(A_0-I)^{\pm1}$ preserve $\Lambda_0$, so $S_e=I+(A_0-I)^{-1}$ preserves it and $x_e=-(A_0-I)^{-1}v_0\in\Lambda_0$. By (iii), $A_u^{\pm1}$ preserve $\Lambda_u$; therefore $S'=I-A_u^{-1}$ preserves $\Lambda_u$, $b=-A_u^{-1}v_u\in\Lambda_u$, $c=-r_u$ is integral on $\Lambda_u$, and $d=r_0x_e-1$ is integral because $x_e\in\Lambda_0$. These are the four blocks of $S_o$ acting on $\Lambda_u\oplus\cO f$.

Conversely, suppose the conditions hold. $S_e-I=(A_0-I)^{-1}$ preserves $\Lambda_0$ and has unit determinant, so $A_0-I$ preserves $\Lambda_0$; hence so does $A_0$, and $v_0=-(A_0-I)x_e\in\Lambda_0$. The upper left block of $S_o$ gives $(I-A_u^{-1})\Lambda_u\subset\Lambda_u$, so $A_u^{-1}$ preserves $\Lambda_u$, and (unit determinant) so does $A_u$; the off-diagonal blocks give $v_u=-A_ub\in\Lambda_u$ and $r_u\Lambda_u\subset\cO$. Hence $\Lambda=\Lambda_0\oplus\Lambda_u$ contributes.

$M=\cO[A]v$ is contributing and $\pi_iM=\cO[A_i]v_i=M_i$. $\cO[S_e]x_e\subset M_0$ because $S_e$ preserves $M_0$ and $x_e\in M_0$; conversely $\cO[S_e]x_e$ is $S_e$-stable, hence $A_0$-stable by the argument just given, and contains $v_0=-(A_0-I)x_e$, so it contains $M_0$. In the same way $\cO[S']b=M_u$. Finally $S_of=(b,d)$, and since $c,d$ are integral on $M_u\oplus\cO f$, induction gives $S_o^{j+1}f=((S')^jb+\sum_{i<j}\gamma_{ji}(S')^ib,\ *)$ with $\gamma_{ji}\in\cO$ and integral last coordinate; so the change from $f,S_of,\dots,S_o^{n-k}f$ to $f,b,S'b,\dots,(S')^{n-k-1}b$ is unitriangular over $\cO$, and $\cO[S_o]f=M_u\oplus\cO f$.
\end{proof}

\subsection*{The weights}

We now show $w(\Lambda)=w_e(\Lambda_0)w_o(\Lambda_u)$. The magnitudes are controlled by Smith normal forms, but the signs are not: the signs record conjugation of genuine lifts, not left--right equivalence of matrices. We therefore work with the genuine function $\cB$ of \eqref{eq:Wcoef} and the action $\T$ of \eqref{eq:Taction}. By \Cref{lem:eps}(a), for every $C$ and $h$,
\begin{equation}\label{eq:summand-T}
\epsilon(C,h)\cW(C^h)=\cB(\T_h(s(C))).
\end{equation}
For the data $C=A,B,S_e,S_o$ define the \emph{normalized lifts}
\[\tilde A=\epsilon(A,C_v)^{-1}s(A),\quad\tilde B=\epsilon(B,C_v)^{-1}s(B),\quad\tilde S_e=\Delta(S_e,x_e)^{-1}s(S_e),\quad\tilde S_o=\Delta(S_o)^{-1}s(S_o).\]
Then, using $\cW^-=\cW$ (\Cref{lem:compact}) and \eqref{eq:summand-T},
\begin{equation}\label{eq:w-B}
w(h\cO^n)=\cB(\T_h\tilde A)\cB(\T_h\tilde B),\quad w_e(h_0\cO^k)=\cB(\T_{h_0}\tilde S_e),\quad w_o(h_u\cO^{n-k})=\cB(\T_{\hat h_u}\tilde S_o).
\end{equation}
We use the following rules; (R3) and (R4) are proved in \Cref{app:technical}.
\begin{enumerate}[label=(R\arabic*)]
\item $\T_h$ is multiplicative and fixes $\pm1$ (\Cref{lem:eps}(a)).
\item $\T_{\diag(h_1,h_2)}(\tilde X\oplus1)=\T_{h_1}(\tilde X)\oplus1$ and $\T_{\diag(h_1,h_2)}(1\oplus\tilde Y)=1\oplus\T_{h_2}(\tilde Y)$ (\Cref{lem:eps}(d)).
\item For a unipotent $u$, $\T_h(\tilde u)$ is the canonical lift of $u^h$; if $u^h$ is integral, this is $\kappa(u^h)$ (\ref{C:unip}, \Cref{lem:unip-kappa}).
\item If $\tilde X$ lies over a matrix whose determinant has even valuation, then $\cB(\kappa(k_1)\tilde X\kappa(k_2))=\cB(\tilde X)$ for $k_1,k_2\in\GL(\cO)$, and $\cB(\tilde X\oplus1)=\cB(1\oplus\tilde X)=\cB(\tilde X)$, where on the left $\cB$ is the function for the bigger group (\Cref{lem:compact,lem:unitblock}).
\end{enumerate}
We write $x\doteqdot y$ if $x=\pm y$ with a sign that depends only on $z$ (not on the lattice).

\begin{proposition}\label{prop:phases}
There is $c(z)\in\{\pm1\}$ such that $w(\Lambda)=c(z)\,w_e(\Lambda_0)\,w_o(\Lambda_u)$ for every contributing $\Lambda$.
\end{proposition}
\begin{proof}
If $\val\det A$ is odd, then $w(\Lambda)=0=w_e(\Lambda_0)$ for all $\Lambda$, since by \ref{S} $\cW(C)=0$ when $\val\det C$ is odd, and $\val\det A^h=\val\det A=\val\det S_e$ (\Cref{lem:factor-identities}(b)). Likewise, if $\val\det B$ is odd, $w(\Lambda)=0=w_o(\Lambda_u)$. Assume now that both are even.

Fix a contributing $\Lambda=\Lambda_0\oplus\Lambda_u$ and choose $h=\diag(h_0,h_u)$ with $\Lambda_0=h_0\cO^k$, $\Lambda_u=h_u\cO^{n-k}$; put $\hat h=\diag(h,1)$ and $\hat h_u=\diag(h_u,1)$. Since $A_0-I$ and $A_u$ preserve $\Lambda_0$, $\Lambda_u$ with unit determinant, $(A_0-I)^{h_0}\in\GL_k(\cO)$ and $A_u^{h_u}\in\GL_{n-k}(\cO)$. Define signs $\alpha(h),\beta(h)$ by
\[\T_{h_u}(s(A_u))=\alpha(h)\,\kappa(A_u^{h_u}),\qquad\T_{h_0}(s(A_0-I))=\beta(h)\,\kappa((A_0-I)^{h_0}).\]
Since $s(A_0-I)s((A_0-I)^{-1})=\sigma(A_0-I,(A_0-I)^{-1})$ is a sign independent of $h$, (R1) gives $\T_{h_0}(s((A_0-I)^{-1}))\doteqdot\beta(h)\kappa((A_0-I)^{h_0})^{-1}$, and similarly $\T_{h_u}(s(A_u^{-1}))\doteqdot\alpha(h)\kappa(A_u^{h_u})^{-1}$.

\emph{The factor $A$.} By \ref{C:block}, $s(A)\doteqdot(s(A_0)\oplus1)(1\oplus s(A_u))$. By (R1), (R2),
\[\T_h\tilde A\doteqdot(\T_{h_0}s(A_0)\oplus1)\,(1\oplus\alpha(h)\kappa(A_u^{h_u})),\]
and $1\oplus\kappa(A_u^{h_u})=\kappa(\diag(I,A_u^{h_u}))$ by \ref{C:compact}. By (R4), $\cB(\T_h\tilde A)\doteqdot\alpha(h)\cB(\T_{h_0}s(A_0))$.

\emph{The factor $B$.} By \Cref{lem:factor-identities}(a), $\diag(B,1)=u_1\diag(A_0-I,T_u)u_2$ with $u_1=\left(\begin{smallmatrix}I&v\\0&1\end{smallmatrix}\right)^{-1}L_0^{-1}$ and $u_2=R_0^{-1}\left(\begin{smallmatrix}I&0\\-r&1\end{smallmatrix}\right)^{-1}$, each a product of two unipotent matrices. Hence $s(B)\oplus1=s(\diag(B,1))\doteqdot\tilde u_1\,(s(A_0-I)\oplus1)(1\oplus s(T_u))\,\tilde u_2$, where $\tilde u_i$ is the product of the canonical lifts of the two factors. The conjugates $u_i^{\hat h}$ of the four unipotent factors are integral: in the basis $\hat h$ this uses $v\in\Lambda$, $r\Lambda\subset\cO$, $r_0(A_0-I)^{-1}\Lambda_0\subset r_0\Lambda_0\subset\cO$ and $(A_0-I)^{-1}v_0\in\Lambda_0$. By (R1)--(R3),
\[\T_{\hat h}(s(B)\oplus1)\doteqdot\kappa(u_1^{\hat h})\,(\beta(h)\kappa((A_0-I)^{h_0})\oplus1)(1\oplus\T_{\hat h_u}s(T_u))\,\kappa(u_2^{\hat h}).\]
Using $\T_{\hat h}(s(B)\oplus1)=\T_hs(B)\oplus1$ and (R4) twice, $\cB(\T_h\tilde B)\doteqdot\beta(h)\cB(\T_{\hat h_u}s(T_u))$.

\emph{The factor $S_e$.} $s(S_e)=s(A_0(A_0-I)^{-1})\doteqdot s(A_0)s((A_0-I)^{-1})$, so $\T_{h_0}\tilde S_e\doteqdot\T_{h_0}s(A_0)\cdot\beta(h)\kappa((A_0-I)^{h_0})^{-1}$ and, by (R4), $\cB(\T_{h_0}\tilde S_e)\doteqdot\beta(h)\cB(\T_{h_0}s(A_0))$.

\emph{The factor $S_o$.} $s(S_o)\doteqdot(s(A_u^{-1})\oplus1)\,s(T_u)\,s(R)$ by \Cref{lem:factor-identities}(a). We have $R^{\hat h_u}=R$ and $\epsilon(R,\hat h_u)=1$: by \ref{C:block}, $\sigma(R,\hat h_u)=1$ and $\sigma(\hat h_u,R)=\hil{\det h_u}{-1}$, while $\hil{\det R}{\det\hat h_u}=\hil{-1}{\det h_u}$. Also $s(R)=\kappa(R)$ by \ref{C:compact}. Hence
\[\T_{\hat h_u}\tilde S_o\doteqdot(\alpha(h)\kappa(A_u^{h_u})^{-1}\oplus1)\,\T_{\hat h_u}s(T_u)\,\kappa(R),\]
and $\cB(\T_{\hat h_u}\tilde S_o)\doteqdot\alpha(h)\cB(\T_{\hat h_u}s(T_u))$.

The four matrices $A^h$, $B^h$, $S_e^{h_0}$, $S_o^{\hat h_u}$ have determinants of even valuation, so (R4) was legitimately applied. Multiplying, $w(\Lambda)\doteqdot\alpha(h)\beta(h)\cB(\T_{h_0}s(A_0))\cB(\T_{\hat h_u}s(T_u))\doteqdot w_e(\Lambda_0)w_o(\Lambda_u)$. The implied sign is a product of signs of the identities in the cover that we used, none of which depend on $\Lambda$.
\end{proof}

\begin{lemma}\label{lem:c=1}
$c(z)=1$.
\end{lemma}
\begin{proof}
If $\val\det A$ or $\val\det B$ is odd, both sides of \Cref{prop:phases} vanish identically and we may take $c(z)=1$. Otherwise evaluate at $\Lambda=M=C_v\cO^n$, which by \Cref{lem:Lsplit} corresponds to $(M_0,M_u)$ with $M_0=C_{x_e}(S_e)\cO^k$ and $M_u\oplus\cO f=C_f(S_o)\cO^{n-k+1}$. Since summands depend only on lattices, \eqref{eq:w-B} and the definition of the normalized lifts give
\[w(M)=\cB(\T_{C_v}\tilde A)\cB(\T_{C_v}\tilde B)=\cB(s(A^{C_v}))\cB(s(B^{C_v})),\]
and similarly $w_e(M_0)=\cB(s(S_e^{C_{x_e}}))$, $w_o(M_u)=\cB(s(S_o^{C_f}))$; the last uses that $\epsilon(S_o,\hat h_u)\cW(S_o^{\hat h_u})$ and $\epsilon(S_o,C_f)\cW(S_o^{C_f})$ coincide when $\hat h_u\cO^{n-k+1}=C_f\cO^{n-k+1}$, by right $\GL(\cO)$-invariance of \eqref{eq:summand-T}. Now $A^{C_v}$ is the companion matrix of $P_A$, which is integral; $B^{C_v}=A^{C_v}+e_1(rC_v)-I$ differs from it by $-I$ and by the integral first row $rC_v=(rv,\dots,rA^{n-1}v)$; $S_e^{C_{x_e}}$ and $S_o^{C_f}$ are companion matrices of integral characteristic polynomials (the matrices preserve $M_0$, resp.\ $M_u\oplus\cO f$). All four are integral upper Hessenberg matrices with subdiagonal entries $1$ and determinants of even valuation. By \Cref{lem:hessenberg} each $\cB$-value is $1$. Hence $1=c(z)\cdot1\cdot1$.
\end{proof}

\begin{proof}[Proof of the linear half of \Cref{thm:factorization}]
By \Cref{lem:Lsplit} the contributing $\Lambda$ correspond bijectively to pairs $(\Lambda_0,\Lambda_u)$ indexing the sums for $J^e(S_e,x_e,r_0)$ and $J^o(S_o)$ (\Cref{lem:models-basic}(c)). By \Cref{prop:phases,lem:c=1}, $w(\Lambda)=w_e(\Lambda_0)w_o(\Lambda_u)$. Summing independently over $\Lambda_0$ and $\Lambda_u$ gives $J^\circ(z)=J^e(S_e,x_e,r_0)J^o(S_o)$.
\end{proof}


\chapter{The fundamental lemma from Fourier compatibility}\label{ch:induction}

By \Cref{thm:factorization}, \Cref{thm:FL} follows from the two Lie algebra identities $J^e=I^e$ and $J^o=I^o$ of \Cref{def:even,def:odd}. We prove them simultaneously, following the induction strategy of Beuzart-Plessis \cite{BP}. At each induction step we fix an \emph{arbitrary} regular semisimple spectral parameter: the matrix $S$ in the even case, or the upper left block $S_0$ in the odd case.  
The key observation is that on the fixed spectral fibre, the difference $J-I$ belongs to a representation of $\SL_2(F)$ (even case) or $\Mp_2(F)$ (odd case). Upper unipotents act by multiplication by $\psi(t\nu)$ and the specified Weyl element fixes the difference. Unit descent and integral support force this difference to be supported where $\nu\in\varpi\cO$, and the group argument in \S\ref{sec:induction} then makes it vanish. Since the spectral parameter was arbitrary, the identity holds on the whole semisimple chart; local constancy extends it to every enhanced regular datum before the next induction step.

Section~\ref{sec:fibres} defines the spectral fibres, auxiliary test spaces and orbital maps. Section~\ref{sec:actions} introduces the group actions and constructs a common representation of the orbital images from the \emph{Fourier compatibility} theorem (\Cref{thm:fourier}). That theorem, proved in \Cref{ch:fourier}, says that the two orbital maps transform by the same integral operator on each fixed spectral fibre.

Throughout this chapter $F$ is as in \Cref{thm:FL}: odd residual characteristic, $q>3$, $\psi$ of conductor $\cO$.

\section{Spectral fibres and auxiliary orbital maps}\label{sec:fibres}

\subsection*{The Lie matching}

When defining the even and odd models, we have constructed cyclic orthogonal quadruples from linear data. We now show every such quadruple is isometric to an even or odd model, respectively, and two models are isometric exactly when their linear data are related by the corresponding \(\GL\)-action.
\begin{lemma}\label{lem:Lie-matching}
	\begin{enumerate}[label=(\alph*)]
		\item The map $(S,x,r)\mapsto(F^{2n},Q^e,C^e,w^e)$ induces a bijection from $\GL_n(F)$-orbits of even data of size $n$ onto isometry classes of quadruples $(W,Q,C,w)$ with $\dim W=2n$, $C\in\so(W,Q)$ invertible and $w$ cyclic for $C$.
		\item The map $S\mapsto(F^{2n-1},Q^o,C^o,w^o)$ induces a bijection from $\GL_{n-1}(F)$-orbits (acting by $S\mapsto S^{\hat g}$) of odd data of size $n$ onto isometry classes of quadruples $(W,Q,C,w)$ with $\dim W=2n-1$, $C\in\so(W,Q)$ and $w$ cyclic for $C$.
		\item Under these bijections the matrix $S$, resp.\ $S_0$, is read off from $C^2$: in (a), $\det(T-C^2)=P_S(T)^2$, and $C^2$ acts on the model by $\diag(S,S)$; in (b), $\det(T-C^2)=T\,P_{S_0}(T)^2$.
	\end{enumerate}
	The stabilizer of a cyclic $(C,w)$ in $\Ort(W)$ is trivial.
\end{lemma}
\begin{proof}
	Let $(W,Q,C,w)$ be given, with $w$ cyclic. As in the proof of \Cref{prop:matching}(c), $E_{\rm ev}={\rm span}(C^{2i}w)$ and $E_{\rm odd}={\rm span}(C^{2i+1}w)$ are orthogonal, nondegenerate, and $W=E_{\rm ev}\oplus E_{\rm odd}$.
	
	(a) Here both have dimension $n$. In the basis $C^{2i}w$ of $E_{\rm ev}$ and $-C^{2i+1}w$ of $E_{\rm odd}$, $C$ maps $E_{\rm ev}$ to $E_{\rm odd}$ by $-I$ and $E_{\rm odd}$ to $E_{\rm ev}$ by a matrix we call $-S$; $S$ is the matrix of $C^2$ on $E_{\rm ev}$, invertible because $C$ is. Let $H$ be the Gram matrix of $E_{\rm ev}$; the Gram matrix of $E_{\rm odd}$ is $(Cx,Cy)=-(x,C^2y)=-x^{\mathsf t}HSy$. With $x=e_1$ and $r=e_1^{\mathsf t}H$ we get $(S,x,r)$ with $(Q,C,w)=(Q^e,C^e,w^e)$; its Hankel matrix $C_x^{\mathsf t}HC_x$ is invertible. Conversely the even model of $(S,x,r)$ is cyclic by \Cref{lem:models-basic}(a). Any isometry between two such quadruples carries the bases $C^aw$ to each other, so the data differ by $g\in\GL_n(F)$; conversely $\diag(g,g)$ is an isometry between the models of $(S,x,r)$ and of its transform (\Cref{lem:models-basic}(b)).
	
	(b) Here $\dim E_{\rm ev}=n$ and $\dim E_{\rm odd}=n-1$. $C$ maps $E_{\rm odd}$ injectively into $E_{\rm ev}$ (it maps $C^{2i+1}w$ to $C^{2i+2}w$), with image $Y={\rm span}(C^2w,\dots,C^{2n-2}w)$, and $E_{\rm ev}=Y\oplus Fw$. Use the basis $\varepsilon_i=-C^{2i+1}w$ $(0\le i<n-1)$ of $E_{\rm odd}$, the basis $-C\varepsilon_i$ of $Y$, and $w$. Then $C(X)=(-X$ in $Y)$ for $X\in E_{\rm odd}$; $C(-C\varepsilon)=-C^2\varepsilon$, so on $Y$ the operator $C$ is $-S_0$ where $S_0$ is the matrix of $C^2|_{E_{\rm odd}}$; and $Cw=-\varepsilon_0=:-b$. Let $H$ be the Gram matrix of $E_{\rm odd}$, $c=b^{\mathsf t}H$ and $d=-Q(w)$. Then $(-C\varepsilon,-C\varepsilon')=-(\varepsilon,C^2\varepsilon')$ gives $-HS_0$, $(-C\varepsilon,w)=(\varepsilon,Cw)$ gives $-Hb$, and so $(Q,C,w)=(Q^o,C^o,w^o)$ for $S=\left(\begin{smallmatrix}S_0&b\\c&d\end{smallmatrix}\right)$. Since $b=e_1$ is cyclic for $S_0$ and $H$ is nondegenerate, the Hankel matrix $C_b^{\mathsf t}HC_b$ is invertible; $S$ is invertible because $\det Q^o=\pm\det(H)^2\det S$. Uniqueness up to $\hat g$ is proved as in (a). Conversely, odd models are cyclic by \Cref{lem:models-basic}(a).
	
	(c) Since $\det(T-C^2)$ is an isometry invariant, it suffices to compute it on the models. In the even model $(C^e)^2=\left(\begin{smallmatrix}0&-S\\-I&0\end{smallmatrix}\right)^2=\diag(S,S)$. In the odd model $C^o(X,Y,z)=(-S_0Y-bz,-X,0)$, so $(C^o)^2(X,Y,z)=(S_0X,\,S_0Y+bz,\,0)$, which is block upper triangular with diagonal blocks $S_0$, $S_0$ and $0$.
	
	Finally an isometry fixing $w$ and commuting with $C$ fixes all $C^aw$.
\end{proof}

\subsection*{The charts and the spectral fibres}
Let $\cB^e$ be the set of characteristic polynomials of regular semisimple $S\in\GL_n(F)$ and $\cB^o$ the set of characteristic polynomials of regular semisimple $S_0\in\GL_{n-1}(F)$ (a point if $n=1$). The \emph{even chart} $\mathcal A^e$ consists of the orbits of regular even data with $P_S\in\cB^e$; the \emph{odd chart} $\mathcal A^o$ consists of the orbits of regular odd data with $P_{S_0}\in\cB^o$. We call $P_S$, resp.\ $P_{S_0}$, the \emph{spectral parameter}. By \Cref{lem:Lie-matching}(c), on the orthogonal side it is the polynomial $P$ determined by
\[
 \det(T-C^2)=P(T)^2\quad\text{(even)},\qquad
 \det(T-C^2)=T\,P(T)^2\quad\text{(odd)}.
\]
For $P\in\cB^e$ or $\cB^o$, write $\mathcal A_P^e$ or $\mathcal A_P^o$ (or simply $\mathcal A_P$, once a parity is fixed) for the corresponding \emph{spectral fibre}. We always use the companion representative $C(P)$ with subdiagonal entries $1$: $S=C(P)$ in even parity and $S_0=C(P)$ in odd parity. For $m=\deg P>0$, put
\[
 K=F[T]/(P),\qquad \vartheta=T\bmod P,
\]
and identify $K$ with $F^m$ through the power basis $1,\vartheta,\ldots,\vartheta^{m-1}$. Multiplication by $\vartheta$ has matrix $C(P)$, identifying $K$ with $F[S]$, resp.\ $F[S_0]$. The element $1\in K$ is the fixed reference vector $x_0=e_1$ in even parity and $b_0=e_1$ in odd parity. These coordinates depend only on $P$; the orbital integrals still integrate over all conjugates of the companion representative.

\begin{remark}\label{rem:chart-density}
The charts are open and dense in the full enhanced regular loci. By the companion normal form (\Cref{lem:free}), regular even orbits have coordinates given by the coefficients of $P_S$ and the moments $rS^ix$ ($0\le i<n$), and regular odd orbits by the coefficients of $P_{S_0}$, the moments $cS_0^ib$ ($0\le i<n-1$), and $d$. The nonvanishing of the Hankel determinant and of $\det S$ defines an open subset of this affine coordinate space. The charts further require $\operatorname{disc}(P_S)\ne0$, resp.\ $\operatorname{disc}(P_{S_0})\det S_0\ne0$. A nonzero polynomial over $F$ cannot vanish on an open box, by the one-variable root bound applied successively to the coordinates. Hence these additional conditions define dense subsets in the local-field topology, both on orbit coordinates and on the full datum spaces. In the odd case $n=1$ the chart is the entire regular locus. 
\end{remark}

\subsection*{A fixed spectral fibre}
\phantomsection\label{par:fibre-coordinates}
We show a fixed spectral fibre in the charts can be identified with $K^\times$ (for even model) or $K^\times\times F^\times$ (odd model). Therefore we can introduce Fourier transform on functions on the spectral fibres.

\paragraph{Even fibre.} Regard the varying column $x$ as an element of $K$, and identify a row $r$ with the unique $\beta\in K$ such that
\[
 r(k)=\Tr_{K/F}(\beta k)\qquad(k\in K).
\]
The trace pairing is nondegenerate because $K$ is \'etale. The regularity condition on even datum is equivalent to $x,\beta\in K^\times$. The centralizer $K^\times$ acts by
\[
 (x,\beta)\longmapsto(t^{-1}x,\beta t),
 \qquad a=\beta x\in K^\times.
\]
The product $a$ is the complete orbit parameter, so $\mathcal A_P^e\simeq K^\times$. We denote by $\nu$ its trace:
\[
 \nu=rx=\Tr_{K/F}(a).
\]
\paragraph{Odd fibre.}
Assume $n\ge2$. Keep $S_0=C(P)$ fixed and write
\[
S=\begin{pmatrix}S_0&b\\ c&d\end{pmatrix}.
\]
Using the fixed power-basis identification $F^{n-1}\simeq K$ as before, let $u\in K$ correspond to the column $b$. Identify the row $c$ with the unique $\gamma\in K$ such that
\[
cS_0^{-1}(k)=\Tr_{K/F}(\gamma k)
\qquad(k\in K),
\]
where $k$ is viewed as a column through $K\simeq F^{n-1}$. Equivalently, since multiplication by $\vartheta^{-1}$ has matrix $S_0^{-1}$,
\[
c(\vartheta^{-1}k)=\Tr_{K/F}(\gamma k).
\]

The regularity condition on the odd datum is equivalent to
\[
u,\gamma\in K^\times.
\]
The centralizer $K^\times=Z_{\GL_{n-1}}(S_0)$ acts by
\[
(u,\gamma,d)\longmapsto(t^{-1}u,\gamma t,d),
\qquad t\in K^\times.
\]
Hence
\[
a=\gamma u\in K^\times,
\qquad
\tau=d-cS_0^{-1}b\in F^\times
\]
are complete orbit parameters. Moreover,
\[
cS_0^{-1}b=\Tr_{K/F}(a),
\qquad
d=\Tr_{K/F}(a)+\tau,
\qquad
\det S=\det S_0\,\tau.
\]
Thus
\[
\mathcal A_P^o\simeq K^\times\times F^\times,
\qquad
(a,\tau)\longmapsto [S].
\]
We define the scalar $\nu$ associated to the invariant $(a,\tau)$ as
\[
\nu=d=\Tr_{K/F}(a)+\tau.
\]

For later use choose representatives $\mathcal D$ for
$F^\times/F^{\times2}$ and, on the component indexed by
$\delta\in\mathcal D$, write
\[
\tau=\delta y^2,\qquad y\in F^\times.
\]
Instead of consider Fourier transform of functions in $(a,\tau)$, we will consider the Fourier transforms of functions in $(a,y)$.

For $n=1$ the spectral base is a point and the fibre is simply
$d\in F^\times$; there is no $K$-variable, and one writes
$d=\delta y^2$. In all subsequent formulas this case is obtained by
suppressing $K$, $a$, and $S_0$.

\paragraph{Orthogonal interpretation of the fibre.}
\phantomsection\label{par:orthogonal-fibre}
By Lemma~\ref{lem:Lie-matching}, each point $p\in\mathcal A_P$ determines an isometry class of quadruples $(W_p,Q_p,C_p,w_p)$, and every compatible quadruple arises in this way. We use this bijection to label orthogonal orbital integrals by the same parameters $a$, or $(a,\tau)$, as on the linear side. The model definitions give
\[
 \begin{aligned}
  \nu&=Q_p(w_p)=\Tr_{K/F}(a)&&\text{in even parity},\\
  \nu&=-Q_p(w_p)=\Tr_{K/F}(a)+\tau&&\text{in odd parity}.
 \end{aligned}
\]
The space $W_p$ may vary with $p$, so orthogonal test data form a family indexed by quadratic spaces. The explicit norm coordinates and centralizer calculations are deferred to \S\ref{sec:even} and \S\ref{sec:odd}.

\subsection*{Auxiliary test spaces and orbital maps}
Fix a parity and $P$. We first define the auxiliary spaces of test vectors and their orbital maps
\[
 R_P:V\longrightarrow\{\text{functions on }\mathcal A_P\},
\]
without specifying the rank-one group actions, which will be introduced in \S\ref{sec:actions}. On the linear side the normalization by $\Delta$ makes the orbital values independent of representatives (\Cref{lem:models-basic}(b)). Absolute convergence of the defining integrals is asserted in \Cref{thm:fourier}(1) and proved in \Cref{ch:fourier}.

\paragraph{Even case, linear side.} Put $V^e_{\rm lin}=\theta_{2n}\otimes\cS(F^n\oplus F^{1\times n})$, the algebraic tensor product. Its elements are finite sums of $\xi\otimes\phi$. With the fixed $S$ and $a=\beta x$ above, put
\[
 R_P(\xi\otimes\phi)(a)
 =\Delta(S,x)^{-1}\int_{\GL_n(F)}
 \epsilon(S,h)\cW_\xi(S^h)\phi(h^{-1}x,rh)\,dh.
\]
Any framing $(x,r)$ with invariant $a$ may be used in this definition.

\paragraph{Odd case, linear side.} Put $V^o_{\rm lin}=\theta_{2n}$. For a matrix $S$ with fixed block $S_0$ and invariants $(a,\tau)$, put
\[
 R_P(\xi)(a,\tau)
 =\Delta(S)^{-1}\int_{\GL_{n-1}(F)}
   \epsilon(S,\hat h)\cW_\xi(S^{\hat h})\,dh.
\]

\paragraph{Orthogonal side.} For each relevant $(W,Q)$ take $V_W=\cS(\so(W)\oplus W)$. If the quadruple attached to $p\in\mathcal A_P$ is realized on $W$, choose a representative $(W,Q,C,w)$ and put
\[
 R_{W,P}(\varphi_W)(p)
 =\int_{\Ort(W)}\varphi_W(h^{-1}Ch,h^{-1}w)\,dh,
\]
Otherwise set the value to zero. Changing the chosen isometry with $W$ conjugates $(C,w)$ by $\Ort(W)$ and leaves the integral unchanged. We also write this value as $R_{W,P}(\varphi_W)(C,w)$ when a representative is displayed. An orthogonal test datum is a family $\boldsymbol\varphi=(\varphi_W)_W$; denote its space by $V_{\rm orth}=\bigoplus_W V_W$. If $W_p$ is the space determined by $p$, then
\[
 R_{{\rm orth},P}(\boldsymbol\varphi)(p)
 =R_{W_p,P}(\varphi_{W_p})(p).
\]
Thus $R_{{\rm orth},P}=\sum_W R_{W,P}$, a sum over the finitely many isometry classes of the relevant dimension, with only one term contributing at each point. Measures are as in \S\ref{sec:FL}. We use $R_P$ notation for $R_{{\rm orth},P}$ and $R_{W_p,P}$ when the same formula applies to both.

\paragraph{Square-root coordinates in odd parity.}
For each fixed $\delta\in\mathcal D$, distinguish the orbital maps in $(a,\tau)$ from their pullbacks under $(a,y)\mapsto(a,\delta y^2)$:
\[
 \begin{aligned}
 R_P^\delta(\xi)(a,y)
   &:=R_P(\xi)(a,\delta y^2),\\
 R_{W,P}^\delta(\varphi_W)(a,y)
   &:=R_{W,P}(\varphi_W)(a,\delta y^2),\\
 R_{{\rm orth},P}^\delta(\boldsymbol\varphi)(a,y)
   &:=R_{{\rm orth},P}(\boldsymbol\varphi)(a,\delta y^2).
 \end{aligned}
\]
Here $a\in K^\times$ and $y\in F^\times$; for $n=1$ omit $a$. The superscript $\delta$ always indicates square-root coordinates, while the maps without it retain the orbit coordinate $\tau$. Since $y$ and $-y$ give the same $\tau$, these pullbacks are even in $y$. We use $R_P^\delta$ for both $R_{{\rm orth},P}^\delta$ and $ R_{W,P}^\delta(\varphi_W)$ when they equal. The pullbacks are initially defined only for $y\ne0$; the extensions to $y=0$ used in the Fourier formulas are supplied in \S\ref{sec:odd}.

\paragraph{Basic vectors.} The linear basic vectors are $v^\circ=\xi^\circ\otimes\one_{\cO^n\oplus\cO^{1\times n}}$ (even) and $v^\circ=\xi^\circ$ (odd). Their orbital images are $J^e|_{\mathcal A_P}$ and $J^o|_{\mathcal A_P}$. The orthogonal basic family $\boldsymbol\varphi^\circ$ has component
\[
 \varphi_W^\circ(C,w)=\one_{\{C:\,CL\subset L\}}(C)\one_L(w)
\]
when $W$ has a selfdual lattice $L$, and zero otherwise. Its orbital image is $I^e|_{\mathcal A_P}$ or $I^o|_{\mathcal A_P}$; the choice of $L$ is irrelevant by \Cref{lem:selfdual}.

\paragraph{Matching on the fixed fibre.}
A linear test vector $v$ and an orthogonal test family $\boldsymbol\varphi$ \emph{match on $\mathcal A_P$} if
\[
 R_P(v)=R_{{\rm orth},P}(\boldsymbol\varphi)
 \quad\text{on }\mathcal A_P.
\]
This is a relation between the orbital maps just defined. The next section equips the test spaces with group actions and states that this relation is preserved. In particular, our goal is to prove the matching of the two basic vectors, first on every spectral fibre and then on the full regular locus.

\section{Group actions and fibrewise Fourier compatibility}\label{sec:actions}

We now equip the test spaces of \S\ref{sec:fibres} with rank-one group actions. These actions are defined on the test vectors; their descent to a common space of orbital functions will follow from Fourier compatibility. Throughout, the spectral parameter $P$ remains fixed.

\subsection*{Weil operators and normalization}
Let $\Mp_2(F)$ be the inverse image of $\SL_2(F)$ in $\tGL_2(F)$. We write
\[n(t)=\begin{pmatrix}1&t\\0&1\end{pmatrix},\qquad w_0=\begin{pmatrix}0&1\\-1&0\end{pmatrix},\qquad h(a)=\begin{pmatrix}a&0\\0&a^{-1}\end{pmatrix},\]
and use the canonical lift $\tilde n(t)$ of $n(t)$ and the lift $\tilde w_0=s(w_0)$. For a nondegenerate symmetric bilinear form $Q$ on a finite dimensional $F$-space $U$ write $Q(u)=Q(u,u)$, let $du$ be selfdual for the pairing $\psi(2Q(u,u'))$, and put
\[\mathcal F_Qf(u')=\int_Uf(u)\psi(2Q(u,u'))\,du .\]
We use the Weil index $\gamma(Q)$ with the Gaussian normalization recorded in \S\ref{sec:oscillator}.
The Weil representation $\omega_Q$ of $\Mp_2(F)$ on $\cS(U)$ is determined by
\begin{equation}\label{eq:weil}
\omega_Q(\tilde n(t))f(u)=\psi(tQ(u))f(u),\qquad\omega_Q(\tilde w_0)f=\gamma(Q)\mathcal F_Qf .
\end{equation}
The compatibility of these formulas with the chosen lifts is justified in \S\ref{sec:oscillator}. If $\dim U$ is even, $\omega_Q$ factors through $\SL_2(F)$. For the hyperbolic form $q(x,r)=rx$ on $F^n\oplus F^{1\times n}$ the Weil index is $1$, the Fourier transform is
\[\mathcal F_{\rm hyp}\phi(x',r')=\int\phi(x,r)\psi(r'x+rx')\,dx\,dr,\]
and we write $\omega_{\rm hyp}=\omega_q$ for the corresponding Weil representation on $\cS(F^n\oplus F^{1\times n})$. It factors through $\SL_2(F)$, and by \eqref{eq:weil}
\[\omega_{\rm hyp}(n(t))\phi(x,r)=\psi(t\,rx)\,\phi(x,r),\qquad \omega_{\rm hyp}(w_0)\phi=\mathcal F_{\rm hyp}\phi .\]

\subsection*{Actions on the auxiliary test spaces}

\paragraph{Even linear action.} On $V^e_{\rm lin}$ let $\SL_2(F)$ act by
\[
 \rho(g)=1\otimes\omega_{\rm hyp}(g).
\]
Only the framing Schwartz variables $(x,r)$ are transformed; the exceptional-representation factor is unchanged under this action.

\paragraph{Odd linear action.} Let $\iota:\tGL_2(F)\to\tGL_{2n}(F)$ lift the embedding
\[
 \begin{pmatrix}\alpha&\beta\\\gamma&\delta\end{pmatrix}
 \longmapsto I_{2n}+(\alpha-1)E_{n,n}+\beta E_{n,2n}
       +\gamma E_{2n,n}+(\delta-1)E_{2n,2n}.
\]
On $V^o_{\rm lin}=\theta_{2n}$ let $\Mp_2(F)$ act by
\[
 \rho(g)=\theta_{2n}(\iota(g)).
\]
This action is defined on the exceptional representation itself.

\paragraph{Orthogonal action.} For a fixed quadratic space \((W,Q)\), the Weil representation attached to the dual pair \(\SL_2\times O(W)\) acts on \(\mathcal S(W)\) and commutes with the natural \(O(W)\)-action. We therefore let it act in the vector variable of
 $V_W=\cS(\so(W)\oplus W)$, and define
\[
 (\rho_W(g)\varphi_W)(C,w)=
 \begin{cases}
  [\omega_Q(g)\varphi_W(C,\cdot)](w),&\text{even parity},\\
  [\omega_{-Q}(g)\varphi_W(C,\cdot)](w),&\text{odd parity}.
 \end{cases}
\]
The action is in the full vector variable $w$, with $W,Q,C$ fixed. It commutes with $\Ort(W)$ and leaves the spectral parameter unchanged. On $V_{\rm orth}$ use the componentwise action. In even parity these actions factor through $\SL_2(F)$; in odd parity we use $\Mp_2(F)$. Write $\rho$ for the action on either source.

\subsection*{Fourier compatibility on a fixed fibre}
\paragraph{The fixed quadratic algebra.}
Put $E=K[U]/(U^2-\vartheta)$, also denoted $K[\sqrt S]$ in even parity and $K[\sqrt{S_0}]$ in odd parity. As in \ref{C:KP}, the square-root symbol denotes the class of $U$ in this quadratic \'etale algebra, not a chosen element of $K$. Factor \(P=\prod_iP_i\) into distinct irreducible factors over \(F\). Since \(P\) is separable,
\[
K=F[T]/(P)\simeq\prod_iK_i,\qquad K_i=F[T]/(P_i).
\]
Accordingly
\[
E=K[U]/(U^2-\vartheta)\simeq\prod_iE_i,
\qquad E_i=K_i[U]/(U^2-\vartheta_i).
\]
Write $\eta_i=\eta_{E_i/K_i}$.

 For the fixed $K,E,\eta,\psi_K$ define
\begin{equation}\label{eq:kernel}
 \mathcal K_\eta(s,a)
 =\int_{K^\times}\eta(t)\psi_K(ta+s/t)\,\frac{dt}{|t|}
 \qquad(s,a\in K^\times).
\end{equation}
Here $dt$ is $\psi_K$-selfdual and the integral is the limit over the annuli $|\val_{K_i}(t_i)|\le N$, which stabilizes (\Cref{lem:shell}). This is a Kloosterman-type kernel when the character is trivial and  Sali\'e-type kernel otherwise.

\begin{theorem}[Fibrewise Fourier compatibility]\label{thm:fourier}
Fix an arbitrary spectral parameter $P$ in either parity, with the identifications above.
\begin{enumerate}[label=(\arabic*)]
\item The orbital integrals of all test vectors on either side converge absolutely at every point of $\mathcal A_P$.
\item For every test vector $v$ on either side and $t\in F$,
\[
 R_P(\rho(\tilde n(t))v)=\psi(t\nu)R_P(v).
\]
Here $\nu=\Tr_{K/F}(a)$ in the even case and $\nu=\Tr_{K/F}(a)+\tau$ in the odd case. On the odd component $\delta$, the pulled-back identity is
\[
 R_P^\delta(\rho(\tilde n(t))v)(a,y)
 =\psi\!\left(t\bigl(\Tr_{K/F}(a)+\delta y^2\bigr)\right)
 R_P^\delta(v)(a,y).
\]
\item The Weyl element acts on the remaining invariants of the same fibre. In the even case,
\[
 R_P(\rho(\tilde w_0)v)(s)
 =\int_K R_P(v)(a)\,\mathcal K_\eta(s,a)\,da
 \qquad(s\in K^\times).
\]
In the odd case, on the fixed component $\delta$,
\[
 \begin{split}
 R_P^\delta(\rho(\tilde w_0)v)(s,y')
 ={}&\gamma(\psi_\delta)\int_K\int_F
 R_P^\delta(v)(a,y)\,\mathcal K_\eta(s,a)\\
 &\hspace{18mm}\cdot\psi(2\delta yy')\,|\delta|^{1/2}\,dy\,da,
 \end{split}
\]
for $s\in K^\times$ and $y'\in F^\times$. The integrals converge absolutely. In each factor $K_i$, the kernel vanishes unless $\eta_i(s_i)=\eta_i(a_i)$, equivalently $s_ia_i^{-1}\in\N(E_i^\times)$ (\Cref{lem:kernel-props}). 
\item If $v$ on the linear side and a family $\boldsymbol\varphi$ on the orthogonal side satisfy
\[
 R_P(v)=R_{{\rm orth},P}(\boldsymbol\varphi)
 \quad\text{on }\mathcal A_P,
\]
then their transforms under every $g\in\SL_2(F)$ (even case), resp.\ $g\in\Mp_2(F)$ (odd case), still satisfy this equality on $\mathcal A_P$. 
\item The basic test vectors $v^\circ$ and $\boldsymbol\varphi^\circ$ are fixed by $\tilde w_0$.
\end{enumerate}
\end{theorem}

\subsection*{The common representation}

Fix a parity and an arbitrary spectral parameter $P$. Write $\mathcal Y_P$ for the space of functions on $\mathcal A_P$, and $R_{{\rm lin},P}$, $R_{{\rm orth},P}$ for the two orbital maps to this space. Their source spaces are those of \S\ref{sec:fibres}, with the actions just defined. Let $G=\SL_2(F)$ in even parity and $G=\Mp_2(F)$ in odd parity. We need a common representation containing the difference of the two basic orbital functions \emph{on this fibre}. The following elementary lemma applies because \Cref{thm:fourier}(4) preserves every matching relation on $\mathcal A_P$.

\begin{lemma}\label{lem:pushout}
Let $\alpha:A\to Y$, $\beta:B\to Y$ be linear maps and let a group $G$ act linearly on $A$ and $B$. Let $M=\{(a,b):\alpha(a)=\beta(b)\}$. If $M$ is $G$-stable, then there is a unique action of $G$ on $\alpha(A)+\beta(B)\subset Y$ such that $\alpha$ and $\beta$ are $G$-equivariant.
\end{lemma}
\begin{proof}
Define $g\cdot(\alpha(a)+\beta(b))=\alpha(ga)+\beta(gb)$. This is well defined: if $\alpha(a)+\beta(b)=\alpha(a')+\beta(b')$ then $(a-a',b'-b)\in M$, so $(g(a-a'),g(b'-b))\in M$, i.e.\ $\alpha(ga)+\beta(gb)=\alpha(ga')+\beta(gb')$. It is clearly an action, and it is the only one making $\alpha,\beta$ equivariant.
\end{proof}

\begin{corollary}\label{cor:common}
There is a representation $\rho_P$ of $G$ on
\[
 \mathcal O_P=R_{{\rm lin},P}(V_{\rm lin})
             +R_{{\rm orth},P}(V_{\rm orth})\subset\mathcal Y_P
\]
making both orbital maps equivariant. On $\mathcal O_P$,
$\rho_P(\tilde n(t))F=\psi(t\nu)F$. The function
\[
 \Phi_P=R_{{\rm lin},P}(v^\circ)
           -R_{{\rm orth},P}(\boldsymbol\varphi^\circ)
       =(J-I)|_{\mathcal A_P}
\]
lies in $\mathcal O_P$ and is fixed by $\rho_P(\tilde w_0)$.
\end{corollary}
\begin{proof}
Apply \Cref{lem:pushout} with $Y=\mathcal Y_P$; its matching space is $G$-stable by \Cref{thm:fourier}(4). Part (2) of that theorem gives the upper-unipotent formula. Equivariance and part (5) give
\[
 \rho_P(\tilde w_0)\Phi_P
 =R_{{\rm lin},P}(\rho(\tilde w_0)v^\circ)
  -R_{{\rm orth},P}(\rho(\tilde w_0)\boldsymbol\varphi^\circ)
 =\Phi_P.
\]
\end{proof}

This construction is made separately for each $P$; it requires no common action linking different spectral fibres. Nor does it assert equality of the two orbital-image spaces, which would be a fibrewise smooth-transfer statement.

\section{Unit descent}\label{sec:descent}

The following identities are on the full enhanced regular loci. They need not stay in a spectral fibre or in the smaller semisimple chart; the induction will therefore use the preceding stage on its full regular locus.

\begin{lemma}[Even to odd]\label{lem:even-odd}
Let $(S,x,r)$ be an even datum of size $n\ge1$ with $\nu=rx\in\cO^\times$. Choose $g\in\GL_n(F)$ whose first $n-1$ columns span $\ker r$ and whose last column is $x$, and write $S^g=\left(\begin{smallmatrix}S_0&b\\c&d\end{smallmatrix}\right)$. Then $S^g$ is an odd datum of size $n$ and
\[J^e(S,x,r)=J^o(S^g),\qquad I^e(S,x,r)=I^o(S^g).\]
\end{lemma}
\begin{proof}
\emph{Normalization.} Replace $r$ by $\nu^{-1}r$. The linear integral does not change: $\one_{\cO^{1\times n}}(rh)$ is invariant under multiplication by a unit, and $\Delta(S,x)$ does not involve $r$. The form $H$ becomes $\nu^{-1}H$; the operator and vector do not change, and a unit multiple of a quadratic form has the same selfdual lattices. So we may assume $rx=1$. Then $g^{-1}x=e_n$ and $rg=e_n^{\mathsf t}$, and by \Cref{lem:models-basic}(b) we may assume $x=f=e_n$, $r=f^{\mathsf t}$.

\emph{Regularity.} Write $S=\left(\begin{smallmatrix}S_0&b\\c&d\end{smallmatrix}\right)$. To see that $(S_0,b,c)$ is regular: the projection of $S^{a+1}f$ to $F^{n-1}$ is $S_0^ab$ plus a combination of $S_0^ib$, $i<a$; since $f$ is cyclic for $S$ (the Hankel matrix of $(S,f,f^{\mathsf t})$ is invertible), $b$ is cyclic for $S_0$. The form $H$ of $(S,f,f^{\mathsf t})$ satisfies $H(f,\cdot)=f^{\mathsf t}$, so $H=\diag(H_0,1)$ with $H_0$ nondegenerate, and $H_0$ is the form of $(S_0,b,c)$ by \Cref{lem:Q0}. Its Gram matrix in the basis $S_0^ib$ is the Hankel matrix of $(S_0,b,c)$, which is therefore invertible.

\emph{Linear side.} A lattice $\Lambda$ with $f\in\Lambda$ and $f^{\mathsf t}\Lambda\subset\cO$ is $\Lambda'\oplus\cO f$ with $\Lambda'=\Lambda\cap F^{n-1}$: for $u\in\Lambda$, $u-(f^{\mathsf t}u)f\in\Lambda\cap F^{n-1}$. So $\Lambda=\hat h\cO^n$ with $\Lambda'=h\cO^{n-1}$, and the summand of $J^e(S,f,f^{\mathsf t})$ at $\Lambda$ (\Cref{lem:models-basic}(c)) is $\Delta(S,f)^{-1}\epsilon(S,\hat h)\cW(S^{\hat h})$, the condition $S\Lambda\subset\Lambda$ being the support condition of $\cW$. Since $\Delta(S,f)=\epsilon(S,C_f)=\Delta(S)$, this is the summand of $J^o(S)$ at $\Lambda'$.

\emph{Orthogonal side.} In coordinates $(X',X_n,Y)\in F^{n-1}\oplus F\oplus F^n$ of $F^n\oplus F^n$, $Q^e=\diag(H_0,1,-HS)$ and $w^e=\varepsilon:=(0,1,0)$, a vector of norm $1$ orthogonal to everything except its own coordinate. A selfdual $L\ni\varepsilon$ splits as $L=L'\perp\cO\varepsilon$ with $L'=L\cap\varepsilon^\perp$ selfdual in $\varepsilon^\perp=\{X_n=0\}$, and every selfdual $L'$ arises. On $\varepsilon^\perp$, with coordinates $(X',Y',Y_n)$, $Q^e$ restricts to $\diag(H_0,-HS)=Q^o$, because $HS=\left(\begin{smallmatrix}H_0S_0&H_0b\\c&d\end{smallmatrix}\right)$. Let $P$ be the orthogonal projection onto $\varepsilon^\perp$ (forgetting $X_n$). Then $C^e(X',0,Y)=(-S_0Y'-bY_n,\,-cY'-dY_n,\,-X',0)$, so $PC^e|_{\varepsilon^\perp}=C^o$, and $-C^e\varepsilon=(0,0,0,1)=w^o$. For $u\in L'$, $C^eu=PC^eu+Q^e(C^eu,\varepsilon)\varepsilon=C^ou+Q^e(u,w^o)\varepsilon$, because $Q^e(C^eu,\varepsilon)=-Q^e(u,C^e\varepsilon)=Q^e(u,w^o)$; the coefficient is integral when $w^o\in L'$. Also $C^e\varepsilon=-w^o\in\varepsilon^\perp$. Hence $C^eL\subset L$ iff $C^oL'\subset L'$ and $w^o\in L'$. This is a bijection between the lattices counted by $I^e(S,f,f^{\mathsf t})$ and by $I^o(S)$.
\end{proof}

\begin{lemma}[Odd to even]\label{lem:odd-even}
Let $S=\left(\begin{smallmatrix}S_0&b\\c&d\end{smallmatrix}\right)$ be an odd datum of size $n\ge1$ with $\nu=d\in\cO^\times$, and put $T=S_0-bd^{-1}c$. If $n=1$ then $J^o(S)=I^o(S)=1$. If $n\ge2$ then $(T,-b,-c)$ is an even datum of size $n-1$ and
\[J^o(S)=J^e(T,-b,-c),\qquad I^o(S)=I^e(T,-b,-c).\]
\end{lemma}
\begin{proof}
If $n=1$, $J^o(S)=\cW(d)=1$ by \Cref{lem:hessenberg} and $I^o(S)=1$ because $\cO$ is the unique selfdual lattice for the form $-dx^2$ and it contains $1$. Let $n\ge2$.

\emph{Data.} $\det S=d\det T$, so $T$ is invertible. The successive spans of $b,Tb,T^2b,\dots$ and of $b,S_0b,S_0^2b,\dots$ coincide, because $T^ab=S_0^ab+$(combination of $S_0^ib$, $i<a$) by induction; similarly for the row spans of $c$. So $b$ is cyclic for $T$, and the form $H$ of $(S_0,b,c)$ satisfies $HT=HS_0-c^{\mathsf t}d^{-1}c$, which is symmetric, and $H(-b)=(-c)^{\mathsf t}$. By \Cref{lem:Q0}, $H$ is also the form of $(T,-b,-c)$, so $(T,-b,-c)$ is regular.

\emph{Linear side.} For a lattice $\Lambda\subset F^{n-1}$,
\[S(\Lambda\oplus\cO f)\subset\Lambda\oplus\cO f\iff T\Lambda\subset\Lambda,\ b\in\Lambda,\ c\Lambda\subset\cO.\]
Indeed, integrality of $S$ gives $b\in\Lambda$, $c\Lambda\subset\cO$, $S_0\Lambda\subset\Lambda$, and then $bd^{-1}c$ is integral on $\Lambda$, so $T\Lambda\subset\Lambda$; the converse is the same computation. These are the contributing lattices of $J^o(S)$ and $J^e(T,-b,-c)$. Use the factorization
\[S=U\,\diag(T,d)\,L,\qquad U=\begin{pmatrix}I&bd^{-1}\\0&1\end{pmatrix},\quad L=\begin{pmatrix}I&0\\d^{-1}c&1\end{pmatrix}.\]
By \ref{C:block} and \ref{C:unip}, $s(S)\doteqdot\tilde U\,(s(T)\oplus1)\,(1\oplus s(d))\,\tilde L$, where $\doteqdot$ means equality up to a sign independent of the lattice. Fix a contributing $\Lambda=h\cO^{n-1}$ and $\hat h=\diag(h,1)$. The matrices $U^{\hat h},L^{\hat h}$ are integral unipotent ($h^{-1}b\in\cO^{n-1}$, $ch\in\cO^{1\times n-1}$), so by (R3) $\T_{\hat h}\tilde U=\kappa(U^{\hat h})$ and $\T_{\hat h}\tilde L=\kappa(L^{\hat h})$. By \Cref{lem:eps}(d), $\T_{\hat h}(1\oplus s(d))=1\oplus s(d)=\kappa(\diag(I,d))$ (\ref{C:compact}), and $\T_{\hat h}(s(T)\oplus1)=\T_hs(T)\oplus1$. If $\val\det S=\val\det T$ is odd, all summands on both sides vanish. Otherwise (R4) gives
\[\cB(\T_{\hat h}s(S))\doteqdot\cB(\T_hs(T)),\]
i.e., by \eqref{eq:summand-T}, $\Delta(S)^{-1}\epsilon(S,\hat h)\cW(S^{\hat h})=c\cdot\Delta(T,-b)^{-1}\epsilon(T,h)\cW(T^h)$ with a sign $c$ independent of $\Lambda$. To see $c=1$ we evaluate at one lattice where the summands do not vanish, if there is one (otherwise both sides are $0$). Then $T$, $b$, $c$ are integral on some $\Lambda$, so $M=\cO[T]b$ is a lattice contained in $\Lambda$, and $M=\cO[S_0]b$ by the span argument above (with $c(M)\subset\cO$). In the basis $C_{-b}(T)$ of $M$, $T$ is an integral companion matrix with determinant of even valuation, so the even summand at $M$ is $\cW(T^{C_{-b}(T)})=1$ by \Cref{lem:hessenberg}. As in the proof of \Cref{lem:Lsplit}, $\cO[S]f=M\oplus\cO f$ (here $c,d$ are integral on $M\oplus\cO f$), so $M\oplus\cO f=C_f(S)\cO^n$ and the odd summand at $M$ equals $\cB(\T_{C_f}(\Delta(S)^{-1}s(S)))=\cW(S^{C_f})=1$, again by \Cref{lem:hessenberg}, $S^{C_f}$ being an integral companion matrix. Hence $c=1$, and summing over lattices gives $J^o(S)=J^e(T,-b,-c)$.

\emph{Orthogonal side.} $w^o$ has norm $-d\in\cO^\times$, so a selfdual $L\ni w^o$ splits as $L=\cO w^o\perp L'$ with $L'=L\cap(w^o)^\perp$ selfdual, and conversely. In the coordinates $(X,Y,z)$, $(w^o)^\perp=\{z=-d^{-1}cY\}$; we use $(X,Y)$ as coordinates on it. Substituting $z=-d^{-1}cY$,
\[Q^o=X^{\mathsf t}HX-Y^{\mathsf t}HS_0Y-2Y^{\mathsf t}Hb\,z-dz^2=X^{\mathsf t}HX-Y^{\mathsf t}(HS_0-c^{\mathsf t}d^{-1}c)Y=X^{\mathsf t}HX-Y^{\mathsf t}HTY.\]
For $u=(X,Y,-d^{-1}cY)$, $C^ou=(-S_0Y+bd^{-1}cY,-X,0)=(-TY,-X,0)$, and its projection to $(w^o)^\perp$ has coordinates $(-TY,-X)$. So the compression of $C^o$ is $\left(\begin{smallmatrix}0&-T\\-I&0\end{smallmatrix}\right)$, and $C^ow^o=(-b,0,0)$ lies in $(w^o)^\perp$ with coordinates $(-b,0)$. Thus $((w^o)^\perp,Q^o,\text{compression},C^ow^o)$ is the even model of $(T,-b,-c)$, since $H(-b)=(-c)^{\mathsf t}$. For $u\in L'$, the $w^o$-component of $C^ou$ is $Q^o(C^ou,w^o)/Q^o(w^o)=Q^o(u,C^ow^o)/d$, which is integral if $C^ow^o\in L'$. Hence $C^oL\subset L$ iff $L'$ is stable under the compression and contains $C^ow^o$. This bijection proves $I^o(S)=I^e(T,-b,-c)$.
\end{proof}

\begin{lemma}[Support]\label{lem:support}
If $\nu\notin\cO$, then $J^e=I^e=0$ (even case), resp.\ $J^o=I^o=0$ (odd case).
\end{lemma}
\begin{proof}
Even linear: a contributing $\Lambda$ has $x\in\Lambda$ and $r\Lambda\subset\cO$, so $rx\in\cO$. Odd linear: $S$ integral on $\Lambda\oplus\cO f$ forces $d\in\cO$. Orthogonal: a vector in a selfdual lattice has integral norm, and $\nu=\pm Q(w)$.
\end{proof}

\section{The induction}\label{sec:induction}

We combine the common representation of \Cref{cor:common} with the unit descents and support statement of \S\ref{sec:descent}.

\begin{theorem}\label{thm:LieFL}
For every $n\ge0$ and every even datum $(S,x,r)$ of size $n$, $J^e(S,x,r)=I^e(S,x,r)$. For every $n\ge1$ and every odd datum $S$ of size $n$, $J^o(S)=I^o(S)$.
\end{theorem}

\begin{proof}
We argue by induction in the order $E_0,O_1,E_1,O_2,E_2,\dots$, where $E_n$ (resp.\ $O_n$) denotes the even (resp.\ odd) identity of size $n$. $E_0$ is $1=1$, and $O_1$ is proved directly: $S=(d)$, $J^o(S)=\cW(d)$ equals $1$ if $\val d\in2\Z_{\ge0}$ and $0$ otherwise (\ref{S} and \Cref{lem:hessenberg}); on the orthogonal side the selfdual lattices for the form $-dx^2$ are $\varpi^m\cO$ with $2m=-\val d$, and $1\in\varpi^m\cO$ iff $m\le0$. So both sides equal $\one_{\val d\in2\Z_{\ge0}}$.

Consider a step $E_n$ ($n\ge1$) or $O_n$ ($n\ge2$), and assume the preceding identity on its full enhanced regular locus. Choose an arbitrary spectral parameter $P$, represented by $S$ in the even case and by $S_0$ in the odd case. Keep it fixed in Steps 1--3 and use $\Phi_P\in\mathcal O_P$ from \Cref{cor:common}.

\emph{Step 1: support on the fibre.} If $\nu(p)\in\cO^\times$ for $p\in\mathcal A_P$, \Cref{lem:even-odd} or \Cref{lem:odd-even} expresses $J(p)$ and $I(p)$ through the same datum of the preceding stage. They are equal by induction, even if that smaller datum has a different spectral parameter or is outside its semisimple chart. If $\nu(p)\notin\cO$, both vanish by \Cref{lem:support}. Hence $\Phi_P(p)=0$ unless $\nu(p)\in\varpi\cO$.

\emph{Step 2: invariance under all upper unipotents.} By Step 1 and \Cref{cor:common}, $\rho_P(\tilde n(t))\Phi_P=\psi(t\nu)\Phi_P=\Phi_P$ for $t\in\varpi^{-1}\cO$. Also $\rho_P(\tilde w_0)\Phi_P=\Phi_P$. Let $\tilde l(t)$ be the canonical lift of $l(t)=\left(\begin{smallmatrix}1&0\\t&1\end{smallmatrix}\right)$. By \ref{C:unip}, $\tilde w_0\tilde n(t)\tilde w_0^{-1}=\tilde l(-t)$, so $\Phi_P$ is fixed by $\tilde l(t)$ for $t\in\varpi^{-1}\cO$. Put
\[
 x=\tilde n(\varpi^{-1})\,\tilde l(-\varpi)\,
       \tilde n(\varpi^{-1})\,\tilde w_0^{-1}.
\]
Every factor fixes $\Phi_P$. Since $n(a)l(-a^{-1})n(a)=h(a)w_0$, this $x$ lies over $h(\varpi^{-1})$, and \ref{C:unip} gives $x\tilde n(t)x^{-1}=\tilde n(\varpi^{-2}t)$. Invariance for $t\in\varpi^{-m}\cO$ therefore implies invariance for $t\in\varpi^{-m-2}\cO$. Repeating this proves
\[
 \psi(t\nu(p))\Phi_P(p)=\Phi_P(p)
 \qquad(t\in F,\ p\in\mathcal A_P).
\]
No property of the center of $\Mp_2(F)$ is used, and all these operations act on the same fibre.

\emph{Step 3: vanishing on the entire fibre.} If $\nu(p)\ne0$, choose $t$ with $\psi(t\nu(p))\ne1$; then $\Phi_P(p)=0$. The locus $\nu\ne0$ is dense even with $P$ fixed. Indeed, at an even datum with $rx=0$, choose a row $\ell$ with $\ell x=1$ and replace $r$ by $r+t\ell$, keeping $S,x$ fixed. Then $\nu=t$, and the Hankel determinant stays nonzero for sufficiently small $t$. At an odd datum with $d=0$, replace $d$ by $t$ while fixing $S_0,b,c$. The Hankel condition is unchanged, the full matrix remains invertible for small $t$, and again $\nu=t$. By \Cref{prop:local-constancy}, the unmultiplied difference $J-I$ is locally constant, so it vanishes at the original point as well. Thus $\Phi_P=0$ on all of $\mathcal A_P$.

\emph{Step 4: all spectral parameters and the full regular locus.} The parameter $P$ was arbitrary, so Step 3 proves $J=I$ on the union of all semisimple fibres, namely the chart. Now fix any enhanced regular datum, including one with a nonsemisimple spectral matrix or singular $S_0$. By \Cref{prop:local-constancy}, $J-I$ is constant on a neighbourhood in all coordinates of the datum. By \Cref{rem:chart-density} that neighbourhood meets the semisimple chart, where the difference is zero. Hence $J=I$ on the full regular locus. This full-locus conclusion is obtained before the next induction step; no Fourier action outside the semisimple fibres is asserted. This completes the induction.
\end{proof}

\begin{proof}[Proof of \Cref{thm:FL}]
If $z$ is not integral both sides vanish (\Cref{lem:integral}). Otherwise, by \Cref{thm:factorization} and \Cref{thm:LieFL},
\[J^\circ(z)=J^e(S_e,x_e,r_0)J^o(S_o)=I^e(S_e,x_e,r_0)I^o(S_o)=I^\circ(z).\qedhere\]
\end{proof}

\chapter{Proof of the Fourier compatibility}\label{ch:fourier}

We prove \Cref{thm:fourier} on one arbitrary spectral fibre. At first sight the orbital integrals are complicated objects: they are integrals over a noncompact reductive group, and on the linear side the integrand is a spherical Shalika coefficient carrying a metaplectic sign. The point of the proof is that most of this structure is \emph{passive} for the action of $\SL_2(F)$ (even case) or $\Mp_2(F)$ (odd case), and that once the passive variables are filtered out the orbital integrals become simple objects, as functions of active parameters. The active orbit parameter is $a$ in even parity; in odd parity we use the pullbacks $R_P^\delta$ and $R_{W,P}^\delta$ in $(a,y)$. 

In even model case, the linear orbital integral becomes a \emph{twisted hyperbolic integral} $\int\phi(u,a/u)\eta(u)\,d^\times u$ of a Schwartz function in two variables over $K$, and the orthogonal orbital integral becomes a \emph{norm integral} of a Schwartz function on a quadratic extension $E/K$. In both cases the Weyl element acts on the Schwartz function by a Fourier transform. Consequently the only identity that has to be proved is a rank one statement over $K$: the \emph{binary kernel identity} of \S\ref{sec:binary} (\Cref{prop:binary}), an elementary Gaussian computation showing that the hyperbolic and the norm radial transforms have the same kernel $\mathcal K_\eta$. No computation of the orbital integrals themselves is needed.

 The odd model case is similar, though the analogue of the Schwartz function above has to be extracted from the exceptional representation itself through a Heisenberg quotient (\S\ref{sec:odd}). This is the only place where the Shalika coefficient is analyzed, with the required description given in \Cref{prop:jacobi}.

Throughout, $F$ is as in \Cref{ch:induction}.

\section{The binary kernel identity}\label{sec:binary}

Recall our setting: $K=\prod_iK_i$ is a finite \'etale $F$-algebra, $\psi_K=\psi\circ\Tr_{K/F}$, and $E=K[\sqrt D]=\prod_iE_i$ a quadratic \'etale $K$-algebra (each $E_i$ a quadratic field or $K_i\times K_i$), with character $\eta=\eta_{E/K}=\prod_i\eta_{E_i/K_i}$. All measures are $\psi_K$-selfdual and $d^\times t=dt/|t|$. Every statement below is a product of statements over the field factors $K_i$, and the reader may assume that $K$ is a field.

For $\phi\in\cS(K^2)$ define the \emph{twisted hyperbolic integral}
\[L\phi(a)=\int_{K^\times}\phi(u,a/u)\eta(u)\,d^\times u\qquad(a\in K^\times),\]
and let $\mathcal F_{\rm hyp}\phi(x',y')=\int\phi(u,v)\psi_K(uy'+vx')\,du\,dv$ be the Fourier transform for the hyperbolic form $uv$. For $\beta\in K^\times$ let $q_\beta(z)=\beta\N_{E/K}(z)$ on $E$, with polar form $2q_\beta(z,z')$. For $\varphi\in\cS(E)$ let $N_\beta\varphi$ be the density of the pushforward of $\varphi(z)\,dz$ under $q_\beta$:
\[\int_E\varphi(z)g(q_\beta(z))\,dz=\int_KN_\beta\varphi(a)g(a)\,da\qquad(g\in C_c^\infty(K)).\]
$N_\beta\varphi$ is supported on $q_\beta(E)$, and on $K^\times$ on $\beta\N(E^\times)$. The isometry classes of the binary forms $q_\beta$ are indexed by $\beta\in K^\times/\N(E^\times)$.

The point of this section is that the Fourier transforms act on the images of $L$ and of $N_\beta$ by one and the same integral operator, whose kernel is the oscillatory integral $\mathcal K_\eta$.

The basic analytic facts we use are proved in \Cref{app:technical} (\Cref{lem:shell,lem:L1}): the annular integrals defining $\mathcal K_\eta$ in \eqref{eq:kernel} stabilize and are bounded by $C(1+\val(a)^+)$ locally uniformly in $s$; $L\phi$ and $N_\beta\varphi$ are compactly supported, satisfy the same logarithmic bound near $0$, and are in $L^1(K)$; and functions in the image of $L$ or $N_\beta$ that vanish on $K^\times$ vanish as distributions (there are no hidden terms at $a=0$). Consequently the radial kernel integrals below converge absolutely; Gaussian integrals are understood in the oscillatory sense of \eqref{eq:gauss} and are justified by compact cutoffs.

\begin{lemma}\label{lem:kernel-props}
$\mathcal K_\eta(s,a)=\mathcal K_\eta(a,s)$, and $\mathcal K_\eta(s,a)=0$ unless $\eta_i(s_i)=\eta_i(a_i)$ for every $i$. In particular, for $s\in\beta\N(E^\times)$ only $a\in\beta\N(E^\times)$ contribute to $\int\cdot\,\mathcal K_\eta(s,a)\,da$.
\end{lemma}
\begin{proof}
Assume first that $K$ is a field. The substitution $t\mapsto t^{-1}$ gives the symmetry, and the substitution $t\mapsto s/(at)$ gives\[\mathcal K_\eta(s,a)=\eta(s)\eta(a)\mathcal K_\eta(s,a).\] For a product the kernel is the product of the kernels of the factors. For nonsplit $E_i$ the kernel of $\eta_i$ is $\N(E_i^\times)$; for split $E_i$, $\eta_i=1$ and $\N(E_i^\times)=K_i^\times$.
\end{proof}

\begin{lemma}[Weil indices of binary norm forms]\label{lem:weil-index}
$\gamma(tq_\beta)=\eta(t)\gamma(q_\beta)$ for $t\in K^\times$, and $\gamma(q_\beta)^2=\eta(-1)$. (Weil indices are taken with respect to $\psi_K$.)
\end{lemma}
\begin{proof}
For a field $K$ with $E=K(\sqrt D)$, $q_\beta\simeq\langle\beta,-\beta D\rangle$, so $\gamma(q_\beta)=\gamma(\psi_\beta)\gamma(\psi_{-\beta D})$ with $\gamma(\psi_a)$ the Weil index of $ax^2$. We use $\gamma(\psi_{ab})\gamma(\psi)=\gamma(\psi_a)\gamma(\psi_b)\hil ab$ and $\gamma(\psi_a)\gamma(\psi_{-a})=1$. The first gives $\gamma(\psi_a)^2=\hil aa\gamma(\psi)^2=\hil a{-1}\gamma(\psi)^2$, and then
$\gamma(\psi_{t\beta})\gamma(\psi_{-t\beta D})=\gamma(\psi_\beta)\gamma(\psi_{-\beta D})\frac{\gamma(\psi_t)^2}{\gamma(\psi)^2}\hil t\beta\hil t{-\beta D}=\gamma(q_\beta)\hil t{-1}\hil t{-D}=\gamma(q_\beta)\hil tD$. Also $\gamma(q_\beta)=\gamma(\psi_\beta)/\gamma(\psi_{\beta D})=\gamma(\psi)/(\gamma(\psi_D)\hil\beta D)$, whose square is $\gamma(\psi)^2/\gamma(\psi_D)^2=\hil D{-1}=\eta(-1)$. In the split case $q_\beta$ is hyperbolic, $\gamma=1$ and $\eta=1$.
\end{proof}

\begin{lemma}[Hyperbolic kernel]\label{lem:hyp-kernel}
For $\phi\in\cS(K^2)$ and $s\in K^\times$, $L(\mathcal F_{\rm hyp}\phi)(s)=\int_KL\phi(a)\,\mathcal K_\eta(s,a)\,da$.
\end{lemma}
\begin{proof}
Restrict first the outer variable $x'$ of $L(\mathcal F_{\rm hyp}\phi)(s)=\int\eta(x')\mathcal F_{\rm hyp}\phi(x',s/x')\,d^\times x'$ to an annulus $-N\le\val(x')\le N$; then everything converges absolutely. In
\[\int_{-N\le\val x'\le N}\eta(x')\int\phi(u,v)\psi_K(us/x'+vx')\,du\,dv\,d^\times x'\]
put $v=a/u$ (so $dv=da/|u|$; the set $u=0$ has measure zero) and then $x'=tu$. Since $\eta(x')=\eta(t)\eta(u)$, the expression becomes
\[\int\!\!\int\phi(u,a/u)\eta(u)\,\mathcal K^{[N,u]}_\eta(s,a)\,d^\times u\,da,\]
where $\mathcal K^{[N,u]}_\eta$ is the integral \eqref{eq:kernel} over the shifted annulus $-N-\val u\le\val t\le N-\val u$. As $N\to\infty$ it tends to $\mathcal K_\eta(s,a)$, and by \Cref{lem:shell} it is bounded by $C_s(1+\val(a)^+)$ independently of $N$ and $u$. The function $|\phi(u,a/u)|(1+\val(a)^+)$ is integrable for $d^\times u\,da$: after the substitution $a=uv$ its integral is $\int|\phi(u,v)|(1+\val(uv)^+)\,du\,dv<\infty$. By dominated convergence the limit is $\int L\phi(a)\mathcal K_\eta(s,a)\,da$. On the left, $\mathcal F_{\rm hyp}\phi$ is Schwartz, so for fixed $s\ne0$ the $x'$-integral has compact support and the truncation is eventually irrelevant.
\end{proof}

\begin{lemma}[Norm kernel]\label{lem:norm-kernel}
For $\varphi\in\cS(E)$ and $s\in K^\times$, $N_\beta(\gamma(q_\beta)\mathcal F_{q_\beta}\varphi)(s)=\int_KN_\beta\varphi(a)\,\mathcal K_\eta(s,a)\,da$.
\end{lemma}
\begin{proof}
Pair both sides with $g\in C_c^\infty(K^\times)$ and write $g(a)=\int\hat g(t)\psi_K(ta)\,dt$ with $\hat g(t)=\int g(s)\psi_K(-ts)\,ds$. Put $q=q_\beta$ and
\[
H(t)=\int_E\mathcal F_q\varphi(w)\psi_K(tq(w))\,dw.
\]
Since $\mathcal F_q\varphi$ is Schwartz, $|H(t)|\le\|\mathcal F_q\varphi\|_{L^1(E)}$. Thus ordinary Fubini gives the absolutely convergent identity
\[
\int_E\mathcal F_q\varphi(w)g(q(w))\,dw
=\int_K\hat g(t)H(t)\,dt.
\]
For $t\ne0$, we claim that
\[
H(t)=\gamma(q)\eta(t)|t|^{-1}
\int_E\varphi(z)\psi_K(-q(z)/t)\,dz.
\]
To prove this, first restrict $w$ to a large compact ball before inserting the Fourier integral for $\mathcal F_q\varphi$. Fubini is then legitimate. After diagonalizing $q$, completing the square translates $w$ by $z/t$. For $z$ in the compact support of $\varphi$ and $t$ in a fixed compact annulus, this translation preserves every sufficiently large ball, and the resulting one-dimensional Gaussian integrals stabilize uniformly. Removing the ball cutoff and applying \eqref{eq:gauss} and \Cref{lem:weil-index} gives the claim. On the original side the cutoff is already irrelevant once the ball contains the support of $\mathcal F_q\varphi$.

Now restrict the integral $\int_K\hat g(t)H(t)\,dt$ to $A_N=\{t\in K^\times:-N\le\val t\le N\}$. Substituting the established formula for $H(t)$, using the pushforward defining $N_\beta$, and inserting $\hat g$, we obtain
\[
\gamma(q)\int_K\!\int_K N_\beta\varphi(a)g(s)
\int_{A_N}\eta(t)\psi_K(-ts-a/t)\,\frac{dt}{|t|}\,ds\,da.
\]
These exchanges are absolutely convergent because $A_N$ is compact, $N_\beta\varphi\in L^1(K)$ and $g$ is compactly supported. The substitution $t\mapsto-1/t$ preserves $A_N$ and turns the inner integral into $\eta(-1)\mathcal K_\eta^{[-N,N]}(s,a)$. By \Cref{lem:weil-index}, $\gamma(q)\eta(-1)=\gamma(q)^{-1}$. On the original side the bound on $H(t)$ justifies the limit $N\to\infty$. On the kernel side, \Cref{lem:shell,lem:L1} give dominated convergence, uniformly for $s$ in the compact support of $g$. Hence
\[
N_\beta(\mathcal F_{q_\beta}\varphi)
=\gamma(q_\beta)^{-1}\int_KN_\beta\varphi(a)\mathcal K_\eta(\cdot,a)\,da
\]
as distributions on $K^\times$. Both sides are continuous on $K^\times$ (the left side is locally constant, the right side continuous by dominated convergence), so they agree pointwise.
\end{proof}

\begin{proposition}[Binary norm transfer]\label{prop:binary}
Say that $\phi\in\cS(K^2)$ and a family $(\varphi_\beta)_{\beta\in K^\times/\N(E^\times)}$, $\varphi_\beta\in\cS(E)$, match if $L\phi(a)=N_\beta\varphi_\beta(a)$ for all $\beta$ and all $a\in\beta\N(E^\times)$. If they match, so do $\omega_{\rm hyp}(g)\phi$ and $(\omega_{q_\beta}(g)\varphi_\beta)_\beta$ for every $g\in\SL_2(F)$; here $\SL_2(F)$ acts through the Weil representations of the $F$-quadratic forms $\Tr_{K/F}(uv)$ and $\Tr_{K/F}\circ q_\beta$.
\end{proposition}
\begin{proof}
$\SL_2(F)$ is generated by the $n(t)$ and $w_0$. The element $n(t)$ multiplies both $\phi$ and $\varphi_\beta$ by $\psi_K(ta)$, where $a$ is the value of the respective quadratic form, so it preserves matching. For $w_0$: the Weil index of $\Tr_{K/F}\circ q_\beta$ with respect to $\psi$ equals that of $q_\beta$ with respect to $\psi_K$, the selfdual measures agree, and the hyperbolic form has Weil index $1$. By \Cref{lem:hyp-kernel,lem:norm-kernel}, both radial transforms are given by the kernel $\mathcal K_\eta$. For $s\in\beta\N(E^\times)$, \Cref{lem:kernel-props} shows that only $a\in\beta\N(E^\times)$ contribute, where the inputs agree; the remaining set $K\smallsetminus K^\times$ has measure zero and carries no mass (\Cref{lem:L1}). Hence the outputs agree. (When $E$ is split there is one class and $q_\beta$ is hyperbolic, and the statement is \Cref{lem:hyp-kernel} itself.)
\end{proof}

\section{The even case}\label{sec:even}

Fix an even fibre and retain the companion representative $S$, the fixed $K$-coordinates, and the orbit parameter $a=\beta x$ of \S\ref{sec:fibres}. We calculate the two orbital maps by integration over their centralizers.

\paragraph{The linear side.}
\begin{lemma}\label{lem:even-gauge}
For $t\in K^\times$, $\epsilon(S,t)=\eta(t)$, and $\Delta(S,x)=\eta(x)$ for $x\in K^\times$.
\end{lemma}
\begin{proof}
The first statement is \eqref{eq:eps-centralizer}. In the basis $1,S,\dots,S^{n-1}$, the cyclic frame $C_x=[x,Sx,\dots,S^{n-1}x]$ is the matrix of multiplication by $x$, so $\Delta(S,x)=\epsilon(S,x)=\eta(x)$.
\end{proof}

\begin{proposition}[Linear centralizer reduction]\label{prop:even-lin}
For $\phi\in\cS(F^n\oplus F^{1\times n})$ and $g\in\GL_n(F)$ put $\phi_g(x',r')=\phi(g^{-1}x',r'g)$, viewed as a function on $K\times K$. Then (with $a=\beta x$ as above)
\[R_P(\xi\otimes\phi)(a)=\int_{K^\times\backslash\GL_n(F)}\epsilon(S,g)\cW_\xi(S^g)\,L(\phi_g)(a)\,d\dot g,\]
and the integral over $K^\times\backslash\GL_n(F)$ is over a compact set (\Cref{lem:compact-outer}). Moreover $(\mathcal F_{\rm hyp}\phi)_g=\mathcal F_{\rm hyp}(\phi_g)$.
\end{proposition}
\begin{proof}
Write $h=tg$ with $t\in K^\times$. Since $t$ centralizes $S$, the cocycle relation and \Cref{lem:even-gauge} give $\epsilon(S,tg)=\epsilon(S,t)\epsilon(S,g)=\eta(t)\epsilon(S,g)$, and $\cW_\xi(S^{tg})=\cW_\xi(S^g)$. Also $\phi(h^{-1}x,rh)=\phi_g(t^{-1}x,rt)$, and $rt$ corresponds to $\beta t$. Hence the integral over $K^\times$ is
\[\int_{K^\times}\eta(t)\phi_g(t^{-1}x,\beta t)\,d^\times t=\eta(x)\int_{K^\times}\eta(u)\phi_g(u,a/u)\,d^\times u=\eta(x)L(\phi_g)(a),\]
after the substitution $u=t^{-1}x$. The factor $\eta(x)$ cancels $\Delta(S,x)^{-1}$ by \Cref{lem:even-gauge}. Finally $(x,r)\mapsto(g^{-1}x,rg)$ preserves the hyperbolic pairing and has Jacobian $1$, so it commutes with $\mathcal F_{\rm hyp}$.
\end{proof}

\paragraph{The orthogonal side: the vector invariant.}
Fix a compatible triple $(W,Q,C)$ over $P$ and vary only the cyclic vector $w$. Lemma~\ref{lem:Lie-matching} already identifies the full orthogonal fibre with $K^\times$ through its linear counterpart $(S,x,r)$; we now describe the identification of $(W,Q,C,w)$ with its invariant $a=a(w)$ directly.
The algebra $K$ acts on $W$ through $C^2$, with $\vartheta\in K$ acting by $C^2$. Define $a(w)\in K$ by
\[
Q(w,kw)=\Tr_{K/F}\bigl(k\,a(w)\bigr)\qquad(k\in K).
\]
The trace pairing makes this definition unique. On the explicit even model, for $f\in F[T]$,
\[
Q^e\bigl(w^e,f((C^e)^2)w^e\bigr)
=r f(S)x=\Tr_{K/F}\bigl(f(S)\beta x\bigr).
\]
Thus $a(w^e)=\beta x=a$, and $Q(w)=\Tr_{K/F}(a(w))=\nu$, where $a,\nu$ are the invariants associated to $(S,x,r)$.

\paragraph{The centralizer fibres.} We describe the fibre of the map $w\mapsto a(w)$.

Choose a cyclic vector $w_0$ in the fixed triple and put $\beta_0=a(w_0)\in K^\times$. The algebra
\[
E=K[C]\simeq K[\sqrt S]
\]
makes $W$ a free rank-one $E$-module, so every cyclic vector has the form $w=zw_0$ with $z\in E^\times$. The adjoint involution on $E$ fixes $K$ and sends $C$ to $-C$; under $E\simeq K[\sqrt S]$ it is the nontrivial $K$-involution $z\mapsto\bar z$. Hence, for $k\in K$,
\[
\begin{aligned}
	Q(zw_0,kzw_0)
	&=Q(w_0,\bar z k z\,w_0)\\
	&=\Tr_{K/F}\!\left(k\,\beta_0\N_{E/K}(z)\right).
\end{aligned}
\]
Thus, with
\[
q_{\beta_0}(z)=\beta_0\N_{E/K}(z),
\]
we have
\[
a(zw_0)=q_{\beta_0}(z).
\]
Moreover an element of $E^\times$ preserves $Q$ exactly when $z\bar z=1$, so
\[
Z_{\Ort(W)}(C)=E^1:=\ker\!\left(\N_{E/K}:E^\times\to K^\times\right).
\]
Consequently the fibres of $w\mapsto a(w)$ on this fixed triple are precisely the $E^1$-orbits, and its image is the norm coset $\beta_0\N(E^\times)$. 

\begin{lemma}\label{lem:even-orth}
	For $\varphi\in\cS(\so(W)\oplus W)$,
	\[
	R_{W,P}(\varphi)(C,w)
	=c_{\beta_0}\int_{E^1\backslash\Ort(W)}
	N_{\beta_0}(\varphi_g)(a)\,d\dot g,
	\qquad
	\varphi_g(z)=\varphi(g^{-1}Cg,g^{-1}zw_0),
	\]
	where $w=z_w w_0$, $a=a(w)=q_{\beta_0}(z_w)$, and $c_{\beta_0}>0$ depends only on the fixed measure normalizations. Although $\varphi_g$ depends on the representative of $g\in E^1\backslash\Ort(W)$, the function $N_{\beta_0}(\varphi_g)$ does not. The outer integral has compact effective support. Under the identification $E\xrightarrow{\sim}W$, $z\mapsto zw_0$, the Weil representation $\omega_Q$ on the vector variable is the Weil representation of $q_{\beta_0}$ with respect to $\psi_K$.
\end{lemma}
\begin{proof}
	Disintegrate Haar measure on $\Ort(W)$ along the closed subgroup
	$E^1=Z_{\Ort(W)}(C)$ and the quotient $E^1\backslash\Ort(W)$. Writing
	$w=z_w w_0$ and using inversion invariance of Haar measure on $E^1$, the inner integral is
	\[
	\int_{E^1}\varphi_g(ez_w)\,de.
	\]
	For $a\in\beta_0\N(E^\times)$ choose $z_a\in E^\times$ with
	$q_{\beta_0}(z_a)=a$. Since
	$q_{\beta_0}^{-1}(a)=E^1z_a$, disintegration of the additive measure $dz$ along $q_{\beta_0}$ gives
	\[
	N_{\beta_0}\phi(a)
	=c(a)^{-1}\int_{E^1}\phi(ez_a)\,de
	\]
	for some positive factor $c(a)$. This factor is constant on the norm coset. Indeed, for $\lambda\in E^\times$ and $\phi_\lambda(z)=\phi(\lambda z)$,
	\[
	N_{\beta_0}(\phi_\lambda)(a)
	=N_{\beta_0}\phi\bigl(\N_{E/K}(\lambda)a\bigr),
	\]
	because multiplication by $\lambda$ scales $dz$ by
	$|\N_{E/K}(\lambda)|_K$, the same factor by which multiplication by
	$\N_{E/K}(\lambda)$ scales $da$. On the other hand,
	\[
	\int_{E^1}\phi_\lambda(ez_a)\,de
	=\int_{E^1}\phi(e\lambda z_a)\,de.
	\]
	Thus $c(a)=c(\N_{E/K}(\lambda)a)$; since $\beta_0\N(E^\times)$ is one orbit under multiplication by $\N(E^\times)$, $c(a)$ is constant there. Denote the resulting normalization by $c_{\beta_0}>0$, so that
	\[
	\int_{E^1}\phi(ez_a)\,de
	=c_{\beta_0}N_{\beta_0}\phi(a).
	\]
	Substituting this identity into the quotient disintegration gives the stated formula. If $g$ is replaced by $e_0g$ with $e_0\in E^1$, then $\varphi_{e_0g}(z)=\varphi_g(e_0^{-1}z)$. Since $q_{\beta_0}(e_0^{-1}z)=q_{\beta_0}(z)$ and multiplication by $e_0$ preserves the chosen additive measure on $E$, the pushforward $N_{\beta_0}(\varphi_g)$ is unchanged. Thus the integrand is well defined on $E^1\backslash\Ort(W)$.
	
	The effective support of the outer integral is compact: the orbit map
	$E^1\backslash\Ort(W)\to\so(W)$, $g\mapsto g^{-1}Cg$, identifies the quotient with the closed regular semisimple orbit of $C$, and the projection of $\operatorname{supp}(\varphi)$ to $\so(W)$ is compact. Finally, the identity above with $k=1$ gives
	$Q(zw_0)=\Tr_{K/F}(q_{\beta_0}(z))$. Hence the identification
	$z\mapsto zw_0$ carries the Weil representation for $(W,Q)$ to that of the $K$-quadratic form $q_{\beta_0}$ with additive character $\psi_K$ (equivalently, of the $F$-quadratic form $\Tr_{K/F}\circ q_{\beta_0}$ with character $\psi$).
\end{proof}

\begin{proposition}[Even Fourier compatibility]\label{prop:even}
	On the fixed even fibre $\mathcal A_P^e$, for test vectors $v$ on either side:
	\[R_P(\rho(\tilde n(t))v)(a)=\psi(t\Tr_{K/F}a)R_P(v)(a),\qquad R_P(\rho(\tilde w_0)v)(s)=\int_KR_P(v)(a)\mathcal K_\eta(s,a)\,da .\]
\end{proposition}
\begin{proof}
	The unipotent formula is immediate from \eqref{eq:weil}, since $\nu=rx=Q^e(w^e)=\Tr_{K/F}a$ and the orbital maps commute with multiplication by functions of the invariant. For the Weyl element on the linear side: by \Cref{prop:even-lin} it suffices to compute $L(\mathcal F_{\rm hyp}\phi_g)$, which is $\int L(\phi_g)(a)\mathcal K_\eta(\cdot,a)\,da$ by \Cref{lem:hyp-kernel}; the compact outer integral commutes with the radial kernel by the uniform bounds of \Cref{lem:shell} (the finitely many slices $\phi_g$ that occur range over a finite dimensional space by smoothness). On the orthogonal side use \Cref{lem:even-orth} and \Cref{lem:norm-kernel}: $\omega_Q(\tilde w_0)$ acts on the slices as $\gamma(q_{\beta_0})\mathcal F_{q_{\beta_0}}$.
\end{proof}

\section{The odd case}\label{sec:odd}

Let $n\ge2$ and fix an odd fibre with its representative $S_0$ and the coordinates of \S\ref{sec:fibres}. The linear source is $\theta_{2n}$, so we first extract the varying coordinates $(b,c,y)$ from this representation through a twisted Jacquet module. The evaluation formula \Cref{prop:jacobi} then expresses the exceptional coefficient in these coordinates, and the centralizer reduction yields the orbital map $R_P^\delta$. The Jacobi calculation requires only invertibility of $S_0$; regular semisimplicity is used for the subsequent centralizer reduction. The case $n=1$ is treated at the end.

\subsection*{The middle embedding}
Order the coordinates of $F^{2n}$ as $I_1=(1,\dots,n-1)$, $I_2=(n,2n)$, $I_3=(n+1,\dots,2n-1)$. Let $P_\iota$ be the permutation matrix with ordered columns $(e_n,e_{2n},e_1,\dots,e_{n-1},e_{n+1},\dots,e_{2n-1})$, so that conjugation by $P_\iota$ carries the block $\diag(\GL_2,1_{2n-2})$ to the coordinates $I_2$. Put $\iota(\tilde g)=\kappa(P_\iota)(\tilde g\oplus1)\kappa(P_\iota)^{-1}$ for $\tilde g\in\tGL_2(F)$. This is a homomorphism of covers; on $\Mp_2(F)$ it is a homomorphism from the double cover, not a splitting. Since $P_\iota$ and $w_0$ are integral, $\iota(\tilde w_0)=\iota(\kappa(w_0))\in\kappa(\GL_{2n}(\cO))$ (we have $s(w_0)=\kappa(w_0)$ by \ref{C:compact}, because $\sigma(w_0,w_0^{-1})=1$ and $w_0^{-1}$ is a simple Weyl representative). By \ref{C:unip}, $\iota(\tilde n(t))$ is the canonical lift of $I+tE_{n,2n}$. For $t_\tau =s(\diag(\tau ,1))$ we have $\iota(t_\tau )=\chi_\iota(\tau )\,s(\diag(I_{n-1},\tau ,I_n))$ for a sign $\chi_\iota(\tau )$; since both $\tau \mapsto t_\tau $ and $\tau \mapsto s(\diag(I_{n-1},\tau ,I_n))$ are homomorphisms (\ref{C:torus}), $\chi_\iota$ is a quadratic character of $F^\times$. (In fact $\chi_\iota=1$, but we do not need this.)

\subsection*{The Heisenberg quotient}
In the order $(I_1,I_2,I_3)$ consider the unipotent group
\[\mathcal H=\left\{h(X,Y,Z)=\begin{pmatrix}I&X&Z\\0&I_2&Y\\0&0&I\end{pmatrix}:X\in\Mat_{n-1,2},\ Y\in\Mat_{2,n-1},\ Z\in\Mat_{n-1}\right\},\]
with group law $h(X,Y,Z)h(X',Y',Z')=h(X+X',Y+Y',Z+Z'+XY')$. Let $N_0=\{h(0,0,Z)\}$, $\psi_{S_0}(h(0,0,Z))=\psi(\tr(S_0Z))$ and
\[J_{S_0}=(\theta_{2n})_{N_0,\psi_{S_0}},\]
the twisted Jacquet module. The subgroup $\{h(0,0,Z):\tr(S_0Z)=0\}$ acts trivially on $J_{S_0}$, and $\mathcal H$ modulo this subgroup is a Heisenberg group with commutator pairing $((X,Y),(X',Y'))\mapsto\psi(\tr(S_0(XY'-X'Y)))$, which is nondegenerate because $S_0$ is invertible. The middle group $\iota(\tGL_2)$ normalizes $\mathcal H$, acting by $X\mapsto Xg^{-1}$, $Y\mapsto gY$, and centralizes $N_0$; so it acts on $J_{S_0}$.

The Schr\"odinger model of this Heisenberg group is $\cS(\Mat_{n-1,2})$ with $[\omega_{S_0}(h(X,Y,Z))f](T)=\psi(\tr(S_0(Z+TY)))f(T+X)$, and $[\omega_{S_0}(g)f](T)=|\det g|^{(n-1)/2}f(Tg)$ is a linear action of $\GL_2(F)$ compatible with the action on $\mathcal H$. By the categorical Stone--von Neumann theorem (\Cref{lem:SvN}),
\[J_{S_0}\simeq\cS(\Mat_{n-1,2})\otimes M_{S_0}\]
as $\mathcal H$-modules, where $M_{S_0}$ is the space of coinvariants of $J_{S_0}$ under $\{h(0,Y,0)\}$, and then the middle cover acts by $\omega_{S_0}\otimes\mu_{S_0}$ for a unique representation $\mu_{S_0}$ of $\tGL_2(F)$ on $M_{S_0}$ (\Cref{lem:SvN}).

\begin{lemma}\label{lem:Mmodule}
	As representations of $\tGL_2(F)$, $\mu_{S_0}\simeq\theta_2\otimes\chi_{S_0}$ with $\chi_{S_0}(\tilde g)=\hil{\det S_0}{\det g}$. In particular $\mu_{S_0}|_{\Mp_2(F)}\simeq\theta_2|_{\Mp_2(F)}$.
\end{lemma}
\begin{proof}
	First let $S_0=I$. In the order $(I_2,I_1,I_3)$, the group $N_0\{h(0,Y,0)\}$ is the unipotent radical of the maximal parabolic with blocks $(n+1,n-1)$, and $\psi_{I}$ is trivial on its first two rows and equals $\psi(\tr Z)$ on the last $n-1$ rows. This is the character of rank $n-1$ in Kaplan's description of the Jacquet modules of exceptional representations; by \cite[Theorem 3.3]{Kaplan} (equivalently \cite[Theorem 4.5]{KM} in the present cocycle convention), the corresponding twisted Jacquet module of $\theta_{2n}$ is, as a representation of the full stabilizer, $|\det|^{(n-1)/2}\theta_2$ on the $\tGL_2$ block tensored with a one dimensional space (the full rank Jacquet module of $\theta_{2n-2}$). The coinvariants under $\{h(0,Y,0)\}$ of $\cS(\Mat_{n-1,2})\otimes M_I$ are $M_I$ tensored with the line of evaluation at $T=0$, on which $\omega_I(g)$ acts by $|\det g|^{(n-1)/2}$ (\Cref{lem:SvN}). Hence $\mu_I\simeq\theta_2$. For general $S_0$ let $A_0=\diag(S_0,1,I_n)$ in the original coordinates. Conjugation by $A_0$ maps $h(X,Y,Z)$ to $h(S_0X,Y,S_0Z)$, so $v\mapsto\theta_{2n}(s(A_0))v$ induces an isomorphism $J_{S_0}\to J_I$ intertwining the $\mathcal H$-actions up to this automorphism. The lifts $s(A_0)$ and $\iota(\tilde g)$ commute up to the scalar $\hil{\det S_0}{\det g}$ (they lie over commuting block diagonal matrices in disjoint coordinates; compute the commutator after a permutation using \ref{C:block}). This gives the twist $\chi_{S_0}$, which is trivial on $\Mp_2(F)$.
\end{proof}

\subsection*{The rank two exceptional representation}
By \Cref{lem:theta2}, $\theta_2$ has a model $R$ on the space of functions $f$ on $\tGL_2(F)$ with values in $\cS(F)^+$ (even Schwartz functions) induced from the subgroup over $G^\square=\{g:\det g\in F^{\times2}\}$, with $\psi$-Whittaker functional $\lambda(f)=f(1)(1)$. Fix representatives $\mathcal D$ of $F^\times/F^{\times2}$ and put $\phi_\delta=|\delta|^{-1/4}f(t_\delta)$ for $\delta\in\mathcal D$. Then (\Cref{lem:theta2}):
\begin{enumerate}[label=(\roman*)]
	\item $\Mp_2(F)$ preserves each component $\delta$ and acts on it by the Weil representation of $\delta y^2$: $\tilde n(t)$ by $\psi(t\delta y^2)$ and $\tilde w_0$ by $\gamma(\psi_\delta)\mathcal F_\delta$, where $\mathcal F_\delta\phi(y')=\int\phi(y)\psi(2\delta yy')|\delta|^{1/2}dy$;
	\item $\lambda(R(t_\tau )f)=|\tau |^{1/4}\phi_\delta(y)$ if $\tau =\delta y^2$.
\end{enumerate}

\subsection*{Jacobi coordinates and the evaluation formula}

Write $T=(T_1,T_2)$ for the two columns of $T\in\Mat_{n-1,2}$, and pass to the mixed model
\[\Phi(b,c)=\int_{F^{n-1}}f(S_0^{-1}b,T_2)\psi(-cT_2)\,dT_2\qquad(b\in F^{n-1},\ c\in F^{1\times n-1}).\]
A direct computation from the formulas for $\omega_{S_0}$ gives, for $g\in\SL_2(F)$,
\begin{equation}\label{eq:mixed}
	\begin{aligned}
		\omega_{S_0}(n(t))\Phi(b,c)
		&=\psi(t\,cS_0^{-1}b)\Phi(b,c),\\
		\omega_{S_0}(w_0)\Phi(b,c)
		&=\int\Phi(b',c')\psi(c'S_0^{-1}b+cS_0^{-1}b')
		\frac{db'\,dc'}{|\det S_0|}.
	\end{aligned}
\end{equation}
For instance, $(Tw_0)=(-T_2,T_1)$, and inverting the partial Fourier transform and substituting $b'=-S_0T_2'$ gives the second formula. The measure $|\det S_0|^{-1}db'\,dc'$ is selfdual for the pairing $c'S_0^{-1}b+cS_0^{-1}b'$. Thus $\omega_{S_0}|_{\SL_2}$ is the Weil representation of the hyperbolic form $cS_0^{-1}b$ (Weil index $1$). The two root subgroups $\{h(X,0,0)\}$ with $X$ supported in the first column, and $\{h(0,Y,0)\}$ with $Y$ supported in the second row, act by translations: $\Phi(b,c)\mapsto\Phi(b+S_0u,c)$, resp.\ $\Phi(b,c)\mapsto\Phi(b,c-vS_0)$; the other two root subgroups act by the characters $\psi(c\alpha)$ and $\psi(\beta b)$.

Fix an isomorphism $M_{S_0}\simeq\theta_2\otimes\chi_{S_0}$ as in \Cref{lem:Mmodule}; it is unique up to a scalar because $\theta_2$ is irreducible. Every $\xi\in\theta_{2n}$ has an image in $J_{S_0}\simeq\cS(F^{n-1}\oplus F^{1\times n-1})\otimes M_{S_0}$, and taking the $\delta$-components in the model $R$ of $\theta_2$ we obtain, for each $\delta\in\mathcal D$, a function
\[\Phi_{\delta,\xi}\in\cS(F^{n-1}\oplus F^{1\times n-1})\otimes\cS(F)^+,\qquad(b,c,y)\mapsto\Phi_{\delta,\xi}(b,c,y),\]
depending linearly on $\xi$. We call these the \emph{Jacobi coordinates} of $\xi$ at $S_0$. By (i) and \eqref{eq:mixed}:

\begin{lemma}\label{lem:jacobi-action}
	$\Phi_{\delta,\theta_{2n}(\iota(g))\xi}=\omega_\delta(g)\Phi_{\delta,\xi}$ for $g\in\Mp_2(F)$, where $\omega_\delta$ is the Weil representation of the quadratic form $q_\delta(b,c,y)=cS_0^{-1}b+\delta y^2$.
\end{lemma}

Next we relate the Jacobi coordinates to Shalika coefficients. For $S=\left(\begin{smallmatrix}S_0&b\\c&d\end{smallmatrix}\right)$ with $\tau =d-cS_0^{-1}b\ne0$ use the Gaussian factorization
\[S=LDU,\qquad L=\begin{pmatrix}I&0\\cS_0^{-1}&1\end{pmatrix},\quad D=\diag(S_0,\tau ),\quad U=\begin{pmatrix}I&S_0^{-1}b\\0&1\end{pmatrix},\]
and define the sign $k(S)$ by $\tilde L\,s(D)\,\tilde U=k(S)\,s(S)$.

\begin{proposition}[Evaluation formula]\label{prop:jacobi}
	There are constants $c(S_0,\delta)$, depending only on $S_0$ and $\delta$, such that for every $\xi\in\theta_{2n}$, every $(b,c)$ and every $y\ne0$, with $\tau =\delta y^2$ and $d=cS_0^{-1}b+\tau $,
	\[\cW_\xi\begin{pmatrix}S_0&b\\c&d\end{pmatrix}=c(S_0,\delta)\,k(S)\,\Phi_{\delta,\xi}(b,c,y).\]
\end{proposition}
\begin{proof}
	Since $C\mapsto m_C$ lifts to a homomorphism of covers carrying canonical lifts to canonical lifts, $s(m_S)=k(S)\,\widetilde{m_L}\,s(m_D)\,\widetilde{m_U}$.
	
	\emph{The factor $m_L$.} $m_L=d_L\cdot\diag(I_n,L^{-1})$. The diagonal Shalika character $\beta$ is trivial on the canonical lift of the unipotent element $d_L$: by \ref{C:unip} and the action of the diagonal torus, canonical lifts of unipotent elements are commutators $[\tilde t,\tilde u]$ of elements of $\tGL_n^{\rm diag}$, and $\beta$ kills commutators. Hence $\ell(\theta(\widetilde{m_L})\,\cdot)=\ell(\theta(\widetilde{\diag(I,L^{-1})})\,\cdot)$. The element $\diag(I_n,L^{-1})$ is $h(0,Y,0)$ with $Y$ supported in the second row (coordinates $(2n,I_3)$), and it commutes with $m_D$; its canonical lift commutes with $s(m_D)$ (the commutator is a $\pm1$-valued homomorphism on a unipotent group).
	
	\emph{The factor $m_D$.} In the original coordinates $m_D=\diag(S_0,\tau ,I_n)=A_0\cdot\diag(I_{n-1},\tau ,I_n)$ with $A_0=\diag(S_0,1,I_n)$, so by \ref{C:block} and the definition of $\chi_\iota$, $s(m_D)=\hil{\det S_0}{\tau }\chi_\iota(\tau )\,s(A_0)\,\iota(t_\tau )$.
	
	\emph{The factor $m_U$.} $m_U=h(X,0,0)$ with $X$ supported in the first column (coordinates $(I_1,n)$).
	
	Hence
	\[\ell(\theta(s(m_S))\xi)=k(S)\hil{\det S_0}{\tau }\chi_\iota(\tau )\;\ell_{S_0}\bigl(\theta(\iota(t_\tau ))\theta(\widetilde{h(0,Y,0)}\,\widetilde{h(X,0,0)})\xi\bigr),\qquad\ell_{S_0}=\ell\circ\theta(s(A_0)).\]
	The functional $\ell_{S_0}$ transforms under $N_0$ by $\psi_{S_0}$, so it factors through $J_{S_0}$. It is invariant under the root subgroups acting by $\psi(c\alpha)$ and $\psi(\beta b)$: these lie in the Shalika unipotent radical, off the diagonal of $X$, and are normalized by $A_0$. A functional on $\cS(F^{n-1}\oplus F^{1\times n-1})\otimes M_{S_0}$ invariant under multiplication by all $\psi(c\alpha)$ and $\psi(\beta b)$ is of the form $\Phi\otimes m\mapsto\Phi(0,0)\lambda_{S_0}(m)$. Moreover $\iota(\tilde n(t))$ is the canonical lift of $u(tE_{nn})$ and commutes with $s(A_0)$, so $\ell_{S_0}$ transforms by $\psi(t)$ under it; since $\omega_{S_0}(n(t))$ does not change the value at $(0,0)$, $\lambda_{S_0}$ is a $\psi$-Whittaker functional on $M_{S_0}\simeq\theta_2\otimes\chi_{S_0}$. By uniqueness of Whittaker functionals, $\lambda_{S_0}=c_1(S_0)\lambda$ in the model above.
	
	The two translations move the point $(0,0)$ to $(b,c)$ (with the parametrization of the root groups by $u=S_0^{-1}b$ and $v=-cS_0^{-1}$; they commute). Then $\iota(t_\tau )$ acts by $\omega_{S_0}(t_\tau )\otimes\mu_{S_0}(t_\tau )$; on the mixed model $\omega_{S_0}(t_\tau )\Phi(b,c)=|\tau |^{(n-1)/2}\Phi(\tau b,c)$, so at $(0,0)$ it contributes $|\tau |^{(n-1)/2}$, and on $M_{S_0}$, by (ii) and $\chi_{S_0}(t_\tau )=\hil{\det S_0}{\tau }$, $\lambda(\mu_{S_0}(t_\tau )m)=\hil{\det S_0}\tau |\tau |^{1/4}\phi_\delta(y)$. Altogether
	\[\ell(\theta(s(m_S))\xi)=k(S)\,\chi_\iota(\delta)\,c_1(S_0)\,|\tau |^{\frac{2n-1}4}\,\Phi_{\delta,\xi}(b,c,y),\]
	using that $\chi_\iota(\tau )=\chi_\iota(\delta)$. Dividing by $|\det S|^{(2n-1)/4}=|\det S_0|^{(2n-1)/4}|\tau |^{(2n-1)/4}$ proves the proposition with $c(S_0,\delta)=\chi_\iota(\delta)c_1(S_0)|\det S_0|^{-(2n-1)/4}$.
\end{proof}

\begin{remark}
	Read backwards, \Cref{prop:jacobi} says that $\Phi_{\delta,\xi}(b,c,y)$ is, up to the sign $k(S)$ and a constant, the Shalika coefficient of $\xi$ at $S$, viewed as a function of the off-diagonal entries and of a square root of the Schur complement. The content of the proposition is that this function is Schwartz in $(b,c,y)$, including at $y=0$ where $\det S=0$, and that $\Mp_2(F)$ acts on it by a Weil representation.
\end{remark}

\subsection*{The cyclic and Gaussian normalizations}
The orbital integral uses the cyclic sign $\Delta(S)$, while \Cref{prop:jacobi} uses the Gaussian sign $k(S)$. Put $G(S)=k(S)\Delta(S)$.

\begin{lemma}[Gauge lemma]\label{lem:gauge}
	Let $S_0$ be regular semisimple, $K=F[S_0]$, $\eta=\eta_{K[\sqrt{S_0}]/K}$, and $b_0$ a cyclic vector for $S_0$. There is a sign $C(S_0,b_0,\delta)$ such that for all $t\in K^\times$, all $c$ and all $\tau \in\delta F^{\times2}$,
	\[G\begin{pmatrix}S_0&tb_0\\c&cS_0^{-1}tb_0+\tau \end{pmatrix}=C(S_0,b_0,\delta)\,\eta(t).\]
\end{lemma}
\begin{proof}
	\emph{Change of basis.} For $h\in\GL_{n-1}(F)$ we compare $G(S^{\hat h})$ and $G(S)$. Conjugating $\tilde Ls(D)\tilde U=k(S)s(S)$ by a lift of $\hat h$, canonical lifts go to canonical lifts, $D^{\hat h}=\diag(S_0^h,\tau )$ and $S^{\hat h}$ has Gaussian factors $L^{\hat h},D^{\hat h},U^{\hat h}$. By \Cref{lem:eps}(a),(d) the two conjugation scalars are $\hil{\det D}{\det h}\epsilon(S_0,h)$ and $\hil{\det S}{\det h}\epsilon(S,\hat h)$; as $\det D=\det S$, $k(S^{\hat h})=k(S)\epsilon(S,\hat h)\epsilon(S_0,h)$. Since $C_f(S^{\hat h})=\hat h^{-1}C_f(S)$, the cocycle relation gives $\Delta(S^{\hat h})=\epsilon(S,\hat h)\Delta(S)$. Hence
	\begin{equation}\label{eq:G-change}
		G(S^{\hat h})=\epsilon(S_0,h)\,G(S).
	\end{equation}
	Applying this with a fixed $h_0$ we may assume that $S_0$ is a companion matrix and $b_0=e_1$.
	
	\emph{The cyclic sign.} Let $P_1=[e_n,e_1,\dots,e_{n-1}]$. For $b=e_1$ the projection of $S^je_n$ to $F^{n-1}$ is $S_0^{j-1}e_1=e_j$ plus a combination of $e_1,\dots,e_{j-1}$, so $P_1^{-1}C_f(S)$ is upper unipotent. For an upper unipotent $R$, $\epsilon(X,R)=1$ by \ref{C:unip} ($s(X)s(R)=s(XR)=s(R)s(X^R)$). Hence $\Delta(S)=\epsilon(S,P_1)$.
	
	\emph{Independence of $c$.} Write $S=S(c,\tau )=L(c)S(0,\tau )$ with $L(c)=\left(\begin{smallmatrix}I&0\\cS_0^{-1}&1\end{smallmatrix}\right)$, $S(0,\tau )=\left(\begin{smallmatrix}S_0&e_1\\0&\tau \end{smallmatrix}\right)=DU$. Then $k(S(0,\tau ))=1$ by \ref{C:unip}. The matrix $P_1^{-1}L(c)P_1=I+e_1(0,cS_0^{-1})$ is upper unipotent. Conjugating $\tilde L(c)s(S(0,\tau ))=k(S)s(S)$ by a lift of $P_1$ and using \ref{C:unip} on the left,
	\[s(P_1^{-1}L(c)P_1)\cdot\hil{\det S}{\det P_1}\epsilon(S(0,\tau ),P_1)\,s(S(0,\tau )^{P_1})=k(S)\hil{\det S}{\det P_1}\epsilon(S,P_1)s(S^{P_1}),\]
	and the left side equals $\hil{\det S}{\det P_1}\epsilon(S(0,\tau ),P_1)s(S^{P_1})$. Therefore $G(S(c,\tau ))=k(S)\epsilon(S,P_1)=\epsilon(S(0,\tau ),P_1)=G(S(0,\tau ))$.
	
	\emph{Independence of the representative of the square class.} $S(0,a^2\tau )=\diag(I,a^2)S(0,\tau )$. By \ref{C:squares}, $s(\diag(I,a^2))s(S(0,\tau ))=s(S(0,a^2\tau ))$, and also $s(\diag(a^2,I))s(X)=s(\diag(a^2,I)X)$ for all $X$, and $\epsilon(\diag(I,a^2),P_1)=1$ since $\diag(I,a^2)^{P_1}=\diag(a^2,I)$. Applying $\T_{P_1}$ we get $\epsilon(S(0,a^2\tau ),P_1)=\epsilon(S(0,\tau ),P_1)$.
	
	\emph{Dependence on $t$.} For $t\in K^\times$, $S^{\hat t}$ has entries $(S_0,t^{-1}b,ct,d)$ and the same Schur complement. By \eqref{eq:G-change} and \eqref{eq:eps-centralizer}, $G(S^{\hat t})=\eta(t)G(S)$. Taking $S$ with $b=tb_0$ and using the previous steps, $G(S)=\eta(t)G(S^{\hat t})=\eta(t)G(S(0,\tau ))$ for the matrix $S(0,\tau )$ with $b=b_0$. Put $C(S_0,b_0,\delta)=\epsilon(S(0,\tau ),P_1)$ (for any $\tau \in\delta F^{\times2}$).
\end{proof}

\subsection*{The odd centralizer reduction}
Use the fixed coordinates of \S\ref{sec:fibres}: $b=ub_0$, $c=c_\gamma$ with $c_\gamma S_0^{-1}(kb_0)=\Tr_{K/F}(\gamma k)$, and $a=\gamma u$. On the chosen component $\delta$, the full matrix varies according to
\[
S=\begin{pmatrix}S_0&ub_0\\c_\gamma&\Tr_{K/F}(a)+\delta y^2\end{pmatrix}.
\]
The fixed $K$-hyperbolic pairing established there shows that $\omega_\delta$ is the hyperbolic Weil representation over $K$ tensored with the scalar oscillator for $\delta y^2$. Define
\[
\tilde\Phi_{\delta,\xi}(u,\gamma,y)
=\Phi_{\delta,\xi}(ub_0,c_\gamma,y).
\]
The reference identifications are transported to $S_0^g$ by $b_0\mapsto g^{-1}b_0$ and $F[S_0^g]\simeq K$; the chosen square-root coordinates $(a,y)$ do not change in this outer integration.

\begin{proposition}\label{prop:odd-lin}
	For $S$ with block $S_0$ and orbit parameters $(a,\tau)$, write $\tau=\delta y^2\ne0$. Then
	\[R_P^\delta(\xi)(a,y)=\int_{Z\backslash\GL_{n-1}(F)}C_g\;L\bigl(\tilde\Phi^{(g)}_{\delta,\xi}\bigr)(a,y)\,d\dot g,\]
	where $Z$ is the centralizer of $S_0$, $\tilde\Phi^{(g)}$ denotes the Jacobi coordinates at $S_0^g$ (with cyclic vector $g^{-1}b_0$), $L$ acts in the variables $(u,\gamma)$ over $F[S_0^g]\simeq K$, $C_g=c(S_0^g,\delta)C(S_0^g,g^{-1}b_0,\delta)$ is a constant, and the outer integral is over a compact set (\Cref{lem:compact-outer}).
\end{proposition}
\begin{proof}
	First, by the cocycle relation and $C_f(S^{\hat h})=\hat h^{-1}C_f(S)$, $\Delta(S)^{-1}\epsilon(S,\hat h)=\Delta(S^{\hat h})^{-1}$, so
	\[R_P^\delta(\xi)(a,y)=\int_{\GL_{n-1}(F)}\Delta(S^{\hat h})\cW_\xi(S^{\hat h})\,dh.\]
	Write $h=gt$ with $t$ in the centralizer of $S_0^g$ and put $S'=S^{\hat g}$. By \Cref{prop:jacobi} at $S_0'=S_0^g$ and \Cref{lem:gauge},
	\[\Delta(S'^{\hat t})\cW_\xi(S'^{\hat t})=G(S'^{\hat t})c(S_0',\delta)\Phi'_{\delta,\xi}(t^{-1}b',c't,y)=C_g\,\eta(t^{-1}u')\,\tilde\Phi'_{\delta,\xi}(t^{-1}u',t\gamma_{c'},y),\]
	where $b'=u'g^{-1}b_0$ and the primes refer to $S_0'$ (identified with $K$ by conjugation by $g$). Substituting $v=t^{-1}u'$, the integral over $t$ is $C_g\int\eta(v)\tilde\Phi'(v,a/v,y)\,d^\times v=C_g\,L\tilde\Phi'(a,y)$, because $\gamma_{c'}u'=a$.
\end{proof}

\paragraph{Uniformity in the outer integration.}
The weighted Jacobi coordinates in \Cref{prop:odd-lin} can all be transported to the single quotient at $S_0$. For $g\in\GL_{n-1}(F)$ put $\hat g=\diag(g,1)$, $D(g)=d_{\hat g}$ and
\[
\xi_g=\beta\bigl(\widetilde{D(g)}\bigr)^{-1}
\theta_{2n}\bigl(\widetilde{D(g)}\bigr)\xi.
\]
This is independent of the lift because both factors are genuine, and defines a smooth ordinary $\GL_{n-1}(F)$-action. By \Cref{lem:eps}(b) and the transformation rule for $\Delta$,
\[
\Delta(S^{\hat g})^{-1}\cW_\xi(S^{\hat g})
=\Delta(S)^{-1}\cW_{\xi_g}(S).
\]
Apply \Cref{prop:jacobi,lem:gauge} to both sides, using $b=ub_0$ and the row $c=c_\gamma$. With $C_0=c(S_0,\delta)C(S_0,b_0,\delta)$, the common factor $\eta(u)$ cancels and gives
\[
C_g\,\tilde\Phi^{(g)}_{\delta,\xi}(u,\gamma,y)
=C_0\,\tilde\Phi_{\delta,\xi_g}(u,\gamma,y).
\]
Initially this holds for $u,\gamma\in K^\times$ and $y\ne0$; since both sides are Schwartz functions, it extends to all $(u,\gamma,y)$ by density. Choose compact representatives of the compact outer support supplied by \Cref{lem:compact-outer}. Smoothness makes $\xi_g$ range over a finite set, so these weighted Schwartz functions have common compact support and uniform bounds. The same support for the outer integral works along the $\iota(\Mp_2(F))$-orbit of $\xi$, again by \Cref{lem:compact-outer}. Moreover, $D(g)$ commutes with the middle $\Mp_2(F)$ in the cover: after permuting the coordinates the blocks are disjoint, and the middle block has determinant one, so the block commutator is trivial. Thus the transport $\xi\mapsto\xi_g$ commutes with this action. These facts justify uniformly the kernel interchanges and the extension to $y=0$ below.

\paragraph{The orthogonal side: the core and the kernel line.}
\phantomsection\label{par:odd-norm-coordinates}
For $n\ge2$, fix a compatible triple $(W,Q,C)$ on the odd fibre. We first describe the invariants $(a,\tau)$ directly from the quadruple $(W,Q,C,w)$. Since $C$ is skewadjoint, the kernel of $C$ is one dimensional. Thus $W$  has an orthogonal decomposition into the core $W_{\rm c}$ and the kernel line $W_0$
\[
W=W_{\rm c}\perp W_0,\qquad
W_{\rm c}=\operatorname{im}C,\qquad W_0=\ker C,\qquad \dim W_0=1,
\]
with $C$ invertible on $W_{\rm c}$. Write $w=w_{\rm c}+w_{\ker}$. Use the form $-Q$, so that its value at $w$ is the scalar $\nu$, and let $K$ act on $W_{\rm c}$ through $C^2$. Define
\[
\Tr_{K/F}\bigl(k\,a(w)\bigr)=-Q(w_{\rm c},kw_{\rm c})\quad(k\in K),
\qquad \tau(w)=-Q(w_{\ker}).
\]
Then $\nu=-Q(w)=\Tr_{K/F}(a(w))+\tau(w)$. We skip the verification that the definitions here match the counterpart from the linear parameters $S$.


\paragraph{The fibres of $w\mapsto(a(w),\tau(w))$.}
Apply the even centralizer-fibre description to the core $(W_{\rm c},-Q|_{W_{\rm c}},C|_{W_{\rm c}})$. Choose a cyclic vector $w_{{\rm c},0}\in W_{\rm c}$ and put $\beta_0=a(w_{{\rm c},0})\in K^\times$. Then $W_{\rm c}$ is a free rank-one $E$-module, every cyclic core vector is $z w_{{\rm c},0}$ with $z\in E^\times$, and
\[
a(z w_{{\rm c},0})=\beta_0\N_{E/K}(z).
\]
Choose $e_0\in W_0$ and $\delta\in\mathcal D$ so that $-Q(e_0)=\delta$. Every cyclic vector of $C$ is then uniquely of the form
\[
w=z w_{{\rm c},0}+y e_0,
\qquad z\in E^\times,\quad y\in F^\times,
\]
and
\[
a(w)=\beta_0\N_{E/K}(z),\qquad
\tau(w)=\delta y^2.
\]
The centralizer preserves the orthogonal decomposition and is
\[
Z_{\Ort(W)}(C)=E^1\times\{\pm1\},
\]
acting by $(z,y)\mapsto(ez,\pm y)$. Consequently the fibres of
\[
w\longmapsto(a(w),\tau(w))
\]
on this fixed triple are precisely the $E^1\times\{\pm1\}$-orbits, and the image is
\[
\beta_0\N(E^\times)\times\delta F^{\times2}.
\]
For $n=1$ only the kernel line remains: $C=0$ and, after choosing $e_0$ with $-Q(e_0)=\delta$, one has $w=ye_0$ and $\tau(w)=\nu=\delta y^2$.

For the centralizer reduction put $Z_C=E^1\times\{\pm1\}$ and
\[
\varphi_g(z,t)=\varphi\bigl(g^{-1}Cg,g^{-1}(z w_{{\rm c},0}+t e_0)\bigr).
\]
Disintegrating Haar measure along $Z_C$ as in \Cref{lem:even-orth} gives
\[
R_{W,P}^\delta(\varphi)(a,y)
=c_{\beta_0}\int_{Z_C\backslash\Ort(W)}
\frac12\sum_{\varepsilon=\pm1}
(N_{\beta_0}\otimes1)(\varphi_g)(a,\varepsilon y)\,d\dot g.
\]
Here $N_{\beta_0}$ acts in the $E$-variable. The Weil action on each slice is the tensor product of the Weil representation of the norm-form core $\beta_0\N_{E/K}$ and that of the line $\delta t^2$. The sign average produces even functions of $y$, matching $\cS(F)^+$ on the linear side. Compactness and the measure convention are as in \Cref{lem:even-orth}.

\begin{proposition}[Odd Fourier compatibility]\label{prop:odd}
	On the fixed odd fibre $\mathcal A_P^o$, for each $\delta\in\mathcal D$ and test vector $v$ on either side, $R_P^\delta(\rho(\tilde n(t))v)=\psi(t\nu)R_P^\delta(v)$ with $\nu=\Tr_{K/F}(a)+\delta y^2$. The function $R_P^\delta(\rho(\tilde w_0)v)$ is obtained from $R_P^\delta(v)$ by applying $\mathcal K_\eta$ in $a$ and $\gamma(\psi_\delta)\mathcal F_\delta$ in $y$.
\end{proposition}
\begin{proof}
	On the linear side, $\rho(g)$ acts on each Jacobi coordinate $\tilde\Phi^{(g')}_{\delta,\xi}$ by $\omega_\delta(g)$ (\Cref{lem:jacobi-action}), and the constants $C_{g'}$ and the outer integral do not involve the variables $(u,\gamma,y)$. The unipotent formula follows since $q_\delta=\Tr_{K/F}(\gamma u)+\delta y^2$ has value $\Tr_{K/F}(a)+\tau =d=\nu$. For $\tilde w_0$, $\omega_\delta(\tilde w_0)$ is $\mathcal F_{\rm hyp}$ over $K$ in $(u,\gamma)$ tensored with $\gamma(\psi_\delta)\mathcal F_\delta$ in $y$; \Cref{lem:hyp-kernel}, applied for each fixed $y$ and integrated against the $y$-transform (all integrals converge absolutely by \Cref{lem:shell,lem:L1}), gives the formula for $L\tilde\Phi$; the common compact support and uniform bounds just proved allow the outer integral to commute with this operator. On the orthogonal side use the fibre description above and \Cref{lem:norm-kernel}; on compact representatives of the outer support, the slices of the fixed Schwartz function range over a finite set, so the same uniform argument applies.
	
	It remains to justify the use of the Fourier transform in $y$, which involves the values at $y=0$, i.e.\ at $\det S=0$, outside the regular fibre. On the linear side these values are those of the Schwartz functions $\Phi_{\delta,\xi}$, so $L\tilde\Phi(a,y)$ is defined and continuous in $y$ for fixed $a\in K^\times$: the integrand $\tilde\Phi(v,a/v,y)$ has compact effective support in $v\in K^\times$, uniformly for $y$ in a compact set and for the weighted slices in the outer integration, by the preceding transport argument. The same holds on the orthogonal side by the finite-slice argument above. The pulled-back functions $R_P^\delta$ and $R_{W,P}^\delta$ therefore extend continuously to $y=0$, uniquely from $y\ne0$; this does not define orbital values at $\tau=0$.
\end{proof}

\paragraph{The case $n=1$.} There is no Heisenberg factor and $\iota$ is the identity of $\tGL_2(F)$. By \Cref{lem:theta2}(ii),
\[
R_P^\delta(\xi)(y)=\cW_\xi(\delta y^2)=\phi_\delta(y).
\]
On a fixed orthogonal line choose $e_0$ with $-Q(e_0)=\delta$ and write $\varphi(t)=\varphi(te_0)$. Then
\[
R_{W,P}^\delta(\varphi)(y)=\tfrac12\bigl(\varphi(y)+\varphi(-y)\bigr).
\]
The Weyl element acts by $\gamma(\psi_\delta)\mathcal F_\delta$ on both.

\section{Proof of \texorpdfstring{\Cref{thm:fourier}}{the Fourier compatibility theorem}}\label{sec:proof-fourier}

Fix an arbitrary spectral parameter $P$ and use throughout the fibre $\mathcal A_P$ of \S\ref{sec:fibres}.

(1) Convergence. At a point of this fibre, \Cref{lem:compact-outer} bounds the outer integration modulo the centralizer $K^\times$. In the even case this follows from the bounded support of $\cW_\xi(S^g)$. In the odd case the bound on $S_0^g$ is independent of $b,c,d$ and uniform along the $\iota(\Mp_2(F))$-orbit of $\xi$. The transport argument in \S\ref{sec:odd} supplies uniform Schwartz bounds for the weighted Jacobi slices on this compact outer support. The centralizer integrals are $L$ of Schwartz functions and converge absolutely for $a\in K^\times$: in each field factor the coordinates $v$ and $a/v$ cannot both remain bounded as $v$ tends to zero or infinity. On the orthogonal side, compact outer support and the norm pushforward give the same conclusion; nonsplit norm-one factors are compact and split factors give the same two-coordinate argument. This proves convergence at the fixed $P$ without any integration over spectral parameters.

(2) and (3) are \Cref{prop:even,prop:odd}, together with the scalar calculation for $n=1$. The norm-sector support is \Cref{lem:kernel-props}. The $L^1$ extensions in \Cref{lem:L1} and the fixed-$S_0$ extension to $y=0$ in \S\ref{sec:odd} justify the radial integral formulas; no extra boundary values or distributions are introduced.

(4) Suppose $R_P(v)=R_{{\rm orth},P}(\boldsymbol\varphi)$ on this fibre alone. By (2), the transforms by $\tilde n(t)$ still match there. By (3), their Weyl transforms are obtained from the same orbital function by the same integral operator on $\mathcal A_P$, so these also match. In odd parity use the equal pullbacks $R_P^\delta(v)=R_{{\rm orth},P}^\delta(\boldsymbol\varphi)$ for each $\delta$. In both parities the kernel preserves the norm sectors. The boundary assertions just cited ensure that equality on the regular fibre determines the radial transforms. Thus the matching relation on this fibre is stable under the generators. As in the source representations, the subgroup generated by $\tilde n(t)$ and $\tilde w_0$ maps onto $\SL_2(F)$, and in odd parity also contains the kernel of $\Mp_2(F)\to\SL_2(F)$ (otherwise it would split the metaplectic cover). Hence the matching relation is stable under the full group in either parity. No equality on another spectral fibre was used.

(5) Even linear: $\mathcal F_{\rm hyp}\one_{\cO^n\oplus\cO^{1\times n}}=\one_{\cO^n\oplus\cO^{1\times n}}$ because $\psi$ has conductor $\cO$ and the measure is selfdual. Odd linear: $\iota(\tilde w_0)\in\kappa(\GL_{2n}(\cO))$ fixes $\xi^\circ$. Orthogonal: for a selfdual lattice $L$, diagonalize the form with unit coefficients. Its Weil index is one, and $\mathcal F_{\pm Q}\one_L=\vol(L)\one_{L^\#}=\one_L$. The factor $\one_{\{C:CL\subset L\}}$ is unchanged. On spaces without a selfdual lattice the basic component is zero. These are identities of the original test vectors, so they remain true after restricting their orbital maps to any $\mathcal A_P$. \qed

\appendix
\chapter{Convergence of the unfolding}\label{app:convergence}

This appendix proves the analytic statements used in \S\ref{sec:unfolding}: absolute convergence of the regular parts on both sides, the entire dependence of the linear regular part on $s$ (so that $s=\frac12$ may be substituted under all sums and integrals), and convergence of the local orbital integrals. We also record the continuation of the two Fourier axes, and explicit classes of test functions for which the full expressions converge: regular support at one place on the orthogonal side (\Cref{prop:A-orth}), and the auxiliary-place class of \Cref{prop:A-full} on the linear side.

Throughout, $k$ is a number field, absolute values are normalized so that the product formula holds, $C_k=k^\times\backslash\A^\times$ and $C_k^1=\ker|\cdot|\subset C_k$, which is compact. We normalize the Haar measure on $C_k$ so that $\int_{C_k}F(|t|)\,d^\times t=\kappa_k\int_0^\infty F(\rho)\,d\rho/\rho$.

\section{Fourier expansions and the axes}

\begin{lemma}\label{lem:A-fourier}
Let $\mathcal C$ be a compact subset of $\GL_{2n}(\A)$.
\begin{enumerate}[label=(\alph*)]
\item The Fourier coefficients $\Theta^\pm_A(g)$, $g\in\mathcal C$, vanish unless $A$ lies in a fixed lattice $\mathcal L\subset\Mat_n(k)$, and for every $N$, $|\Theta^\pm_A(g)|\le C_N(1+\|A\|_\infty)^{-N}$. The Fourier series of $\Theta^\pm$ along $[U]$ converge absolutely and uniformly on $[U]\mathcal C$.
\item $|C_Y(g)|\le\int_{[U]}|\Theta^+\Theta^-(\tilde u(X)g)|\,dX$, uniformly in $Y$.
\item For every $N$ there is $C_N$ with
\[\sum_{r\ne0,\,v\ne0}|\widehat\Phi_g(tr,t^{-1}v)|\le C_N\min(|t|,|t|^{-1})^N\qquad(g\in\mathcal C,\ t\in\A^\times).\]
\end{enumerate}
Consequently, for fixed $g$ and $\Re s>1$ all interchanges in the proof of \Cref{lem:poisson} are justified, and the doubly nonzero part of the Mellin integral,\[\int_{C_k}\eta(t)|t|^{n(2s-1)}\sum_{r\ne0,v\ne0}C_{I-vr}(g)\widehat\Phi_g(tr,t^{-1}v)\,d^\times t,\]is entire in $s$, locally uniformly in $g\in\mathcal C$.
\end{lemma}
\begin{proof}
(a) The finite components of $\Theta^\pm$ are fixed by an open compact subgroup $\mathcal U$ of $U(\A_f)$, uniformly for $g\in\mathcal C$; the Fourier coefficient at $A$ is invariant under $\mathcal U$ and is multiplied by $\psi(\tr AX)$, so it vanishes unless $\psi(\tr A\,\cdot)$ is trivial on $\mathcal U$, which confines $A$ to a lattice. At the archimedean places, each derivative in a coordinate of $X$ multiplies the coefficient by the corresponding frequency coordinate, while all derivatives of $\Theta^\pm(\tilde u(X)g)$ are bounded on the compact set $[U]\mathcal C$. A lattice has polynomial counting growth, so taking $N$ larger than its real dimension gives absolute uniform convergence. (b) is immediate. (c) Choose a continuous section of the norm, so that $t$ may be taken in a fixed compact set times that section modulo $k^\times$. Dominate $\Phi$ by a product of an archimedean Schwartz function and a lattice indicator. For $|t|\ge1$ use arbitrary decay in the nonzero coordinate $tr$ and sum the possibly contracting coordinate $t^{-1}v$ over a lattice at polynomial cost; for $|t|\le1$ interchange the roles. Compactness of $C_k^1$ makes the bound uniform. The final assertion follows from (b), (c) and differentiation under the integral sign.
\end{proof}

\begin{proposition}[The axes]\label{prop:A-tate}
For a Schwartz function $\phi$ on $\A^n$ put $Z(\phi,\eta;u)=\int_{C_k}\eta(t)|t|^{nu}\sum_{x\in k^n\smallsetminus0}\phi(tx)\,d^\times t$ $(\Re u>1)$. The terms of \Cref{lem:poisson} with $r=0$, together with the term $-\Phi(0)C_{I_n}$, and the terms with $v=0$, $r\neq0$, contribute \footnote{In the mixed variables $(r,v)\in k^{1\times n}\times k^n$ of \Cref{lem:poisson}, we call the two coordinate subspaces $\{r=0\}$ and	$\{v=0\}$ the axes. Their contribution is partitioned	into the terms with $r=0$ and arbitrary $v$, and the terms with $v=0$ and $r\ne0$, so that the origin is	counted only once.}
\[|\det g|^sC_{I_n}(g)\bigl(Z(\phi_{1,g},\eta;2s)+Z(\phi_{2,g},\eta;2s-1)\bigr),\qquad\phi_{1,g}(x)=\Phi((0,x)g),\ \phi_{2,g}(r)=\widehat\Phi_g(r,0).\]
Moreover
\[Z(\phi_{1,g},\eta;2s)+Z(\phi_{2,g},\eta;2s-1)=\mathcal E_g(s)+\one_{\eta=1}\frac{\kappa_k}n\Bigl(-\frac{\Phi(0)}{2s}+\frac{|\det g|^{-1}\widehat\Phi(0)}{2s-2}\Bigr)\]
with $\mathcal E_g$ entire. In particular $s=\frac12$ is a regular point of the sum of the two axes, although its value there need not vanish.
\end{proposition}
\begin{proof}
For $r=0$ the frequency is $I_n-v\cdot0=I_n$. Poisson summation in $v$ gives $|t|^{n(2s-1)}\sum_v\widehat\Phi_g(0,t^{-1}v)=|t|^{2ns}\sum_x\Phi((0,tx)g)$, and subtracting $|t|^{2ns}\Phi(0)$ removes $x=0$. For $v=0$ the term is $|t|^{n(2s-1)}\widehat\Phi_g(tr,0)$, and $\widehat\Phi_g(r,0)=\int\Phi((r,x)g)\,dx$ is the integral over a whole fibre. This proves the first formula. For the second, split $Z$ at $|t|=1$ and apply Poisson summation on $|t|<1$:
\[Z(\phi,\eta;u)=Z^{\ge1}(\phi,\eta;u)+Z^{\ge1}(\hat\phi,\eta;1-u)+\one_{\eta=1}\frac{\kappa_k}n\Bigl(\frac{\hat\phi(0)}{u-1}-\frac{\phi(0)}u\Bigr),\]
where $Z^{\ge1}$ is the integral over $|t|\ge1$ (entire), and the constant terms vanish for nontrivial $\eta$ because $\eta$ integrates to zero over $C_k^1$. Apply this to the two axes: the possible residues at $s=\frac12$ cancel, since
\[
\widehat\varphi_{1,g}(0)
=
\int_{\mathbb A^{1\times n}}\Phi((0,x)g)\,dx
=
\widehat\Phi_g(0,0)
=
\varphi_{2,g}(0).
\]
The remaining constants are $\phi_{1,g}(0)=\Phi(0)$ and $\hat\phi_{2,g}(0)=|\det g|^{-1}\widehat\Phi(0)$.
\end{proof}

\section{The local orbital integrals}

\begin{lemma}[Schwartz majorants for Shalika coefficients]\label{lem:A-majorant}
Over any local field there is $\alpha>0$ (depending on the representation) such that every vector $\xi$ has a nonnegative Schwartz function $\varphi_\xi$ on $\Mat_n$ with $|\ell(\theta(s(m_C))\xi)|\le|\det C|^{-\alpha}\varphi_\xi(C)$; one $\varphi$ works for a compact family of translates of $\xi$. At an unramified place one may take $\varphi=\one_{\Mat_n(\cO)}$ for the spherical vector, for every $\alpha\ge0$.
\end{lemma}
\begin{proof}
The first statement is \cite[Claim 3.2]{KM}. For a compact family of translates: at a nonarchimedean place the family is finite by smoothness; at an archimedean place, for each $N$ the linear map $\xi\mapsto|\det C|^\alpha(1+\|C\|)^N\ell(\theta(s(m_C))\xi)$ from the Fr\'echet space of smooth vectors to bounded functions is everywhere defined and has closed graph, hence is continuous, and continuous seminorms are bounded on compact sets. The unramified statement follows from \ref{S}: $|\cW(C)|=q^{\sum(i-1)d_i}$ and $|\det C|^{(2n-1)/4}=q^{-\sum_i(n-\frac12)d_i}$, so $|\ell(\theta(s(m_C))\xi^\circ)|=q^{-\sum_i(n-i+\frac12)d_i}\le\one_{\Mat_n(\cO)}(C)$.
\end{proof}

\begin{proposition}\label{prop:A-local}
Let $z\in X^{\rm reg}(F)$ over a local field. The integral $O_z(s)$ of \eqref{eq:rawOI} converges absolutely and is entire in $s$; over a nonarchimedean field the $h$-integration is over a compact set. The bounds are uniform for $g$ in a compact set.
\end{proposition}
\begin{proof}
Nonarchimedean case: the finitely many vectors $\theta^\pm(\tilde g)\xi^\pm$, $g\in{\rm supp}f$, are fixed by an open compact subgroup of $U(F)$; by the character relation, $\ell(\theta^+(s(m_{A^h}))\theta^+(\tilde g)\xi^+)\ne0$ forces $A^h$ into a fixed compact set. The Schwartz factor $\widehat\Phi_g(rh,h^{-1}v)$ confines $h^{-1}v$ and $rh$ to compact sets. By \Cref{lem:free}, $C_{h^{-1}v}(A^h)=h^{-1}C_v$ and $R_{rh}(A^h)=R_rh$; the left sides are bounded, and $C_v,R_r$ are invertible, so $h^{-1}$ and $h$ are bounded. Archimedean case: on the orbit the determinants of $A^h,B^h$ are constant, so \Cref{lem:A-majorant} bounds the integrand by a Schwartz function of $(A^h,B^h,h^{-1}v,rh)$; by the same two frame identities this decays faster than any power of $\|h\|+\|h^{-1}\|$, while Haar measure and $|\det h|^{2\Re s-1}$ grow polynomially.
\end{proof}

\section{Global summability}

\begin{theorem}[Summability on the cyclic locus]\label{thm:A-sum}
For $m\ge0$ let $\mathcal V_m=\Mat_n^{m+1}\oplus k^n\oplus k^{1\times n}$, with points $(A,M_1,\dots,M_m,v,r)$ and $\GL_n$ acting by conjugation on the matrices and by $(v,r)\mapsto(h^{-1}v,rh)$. Let $\mathcal V_m^{\rm cyc}$ be the locus where the Hankel matrix $R_rC_v$ of $(A,v,r)$ is invertible. For $\Psi\in\cS(\mathcal V_m(\A))$ and $R>0$,
\[\sum_{x\in\mathcal V_m^{\rm cyc}(k)/\GL_n(k)}\int_{\GL_n(\A)}|\Psi(x\cdot h)|\,|\det h|^{\tau}\,(1+|\log|\det h||)^R\,dh<\infty,\]
uniformly for real $\tau$ with $|\tau|\le R$. The same holds for any $\GL_n(k)$-stable subset of $\mathcal V^{\rm cyc}_m(k)$.
\end{theorem}
\begin{proof}
\emph{Canonical coordinates.} Let $P_A=T^n+a_1T^{n-1}+\dots+a_n$, $b_i=rA^iv$, $\mathbf H=(b_{i+j})_{i,j<n}$ (with $b_i$ for $i\ge n$ defined by the recursion from $P_A$) and $\mathbf D_\ell=(rA^iM_\ell A^jv)_{i,j<n}$. These are polynomial invariants. By \Cref{lem:free}, conjugation by $C_v$ gives the unique representative
\[x_c=(A_c,\ \mathbf H^{-1}\mathbf D_1,\dots,\mathbf H^{-1}\mathbf D_m,\ e_1,\ (b_0,\dots,b_{n-1}))\]
of the orbit, with $A_c$ the companion matrix; so the invariants $\mathbf c=(a,b,\mathbf D)$ parametrize the rational orbits, with $\Delta(\mathbf c)=\det\mathbf H\ne0$. Changing the representative by a rational matrix does not change the integral, as $|\det|_\A=1$ on $\GL_n(k)$.

\emph{Reduction to integral data.} Dominate $\Psi$ by a product of an archimedean Schwartz function and the indicator of a lattice; multiplying all coordinates by a suitable $d\in k^\times$ (an equivariant bijection of rational cyclic points) we may assume the finite part is supported on integral points. Then the contributing $\mathbf c$ are integral.

\emph{Finite places.} At a finite place $w$, if $x_c\cdot h$ is integral, $\Lambda=h\cO_w^n$ satisfies $A_c\Lambda\subset\Lambda$, $e_1\in\Lambda$, $r\Lambda\subset\cO_w$, so $\cO_w^n=\cO_w[A_c]e_1\subset\Lambda\subset\mathbf H^{-1}\cO_w^n$. The quotient has length $\ell_w=\val_w\Delta(\mathbf c)$, so there are at most $q_w^{n\ell_w}$ such $\Lambda$, each of volume one, and $1\le|\det h|_w\le q_w^{\ell_w}$. The product over the finite places is at most $|N_{k/\Q}\Delta(\mathbf c)|^{n+R}$ (with logarithms absorbed by enlarging $R$).

\emph{Archimedean places.} Put $y=x_c\cdot h$. Then $C_{e_1}(y)=h^{-1}$ and $R(y)=\mathbf Hh$, so $h^{-1}$ is polynomial in $y$ and $h=\mathbf H^{-1}R(y)$. By Cramer's rule, a two-sided height $\mathrm{Ht}(h)$ bounding $\|h\|$ and $\|h^{-1}\|$ satisfies $\mathrm{Ht}(h)\le C(1+\|\mathbf c\|)^D\prod_{w|\infty}\max(1,|\Delta(\mathbf c)|_w^{-1})^D(1+\|y\|)^D$. Haar measure on $\GL_n(k_\infty)$ has polynomial growth in $\mathrm{Ht}$, so $\int\mathrm{Ht}(h)^{-P}dh<\infty$ for large $P$, and $|\det h|^\tau$ is bounded by a power of $\mathrm{Ht}(h)$. Since $1+\|\mathbf c\|\ll(1+\|y\|)^{D_0}$, splitting the Schwartz decay of $\Psi$ gives
\[\int_{\GL_n(k_\infty)}|\Psi_\infty(x_c\cdot h)||\det h|^\tau\,dh\le C_N(1+\|\mathbf c\|)^{-N}\prod_{w|\infty}\max(1,|\Delta(\mathbf c)|_w^{-1})^{D_R}.\]

\emph{No escape to the boundary.} $\Delta(\mathbf c)$ is a nonzero algebraic integer, so $\prod_{w|\infty}|\Delta(\mathbf c)|_w\ge1$ and therefore $\prod_{w|\infty}\max(1,|\Delta|_w^{-1})\le\prod_{w|\infty}\max(1,|\Delta|_w)$, which is polynomial in $1+\|\mathbf c\|$; so is $|N\Delta(\mathbf c)|$. The orbit integrals therefore decay faster than any power of $1+\|\mathbf c\|$, and $\mathbf c$ runs over a lattice.
\end{proof}

The cyclicity is doing real work here: the frame $C_v$ bounds $h^{-1}$ and the frame $R_r$ bounds $h$; with one framing vector only, a parabolic direction could escape.

\section{The regular parts}

\begin{theorem}\label{thm:A-regular}
The integral \eqref{eq:Jreg} converges absolutely for every $s$, locally uniformly, and all rearrangements in the proof of \Cref{thm:linear-unfold} are justified. In particular $J^{\rm reg}$ is entire.
\end{theorem}
\begin{proof}
By \Cref{lem:factor}, \Cref{lem:A-majorant} at every place (with a common exponent $\alpha$, enlarging the finitely many ramified majorants if necessary), and the product formula $\prod_v|\det A^h|_v^{-\alpha}=1$ for rational $A$, the product $|\Theta^+_A\Theta^-_B(d_hg)|$ is bounded by $\varphi(A^h)\varphi'(B^h)$ for adelic Schwartz functions $\varphi,\varphi'$, uniformly for $g$ in the compact support of $f_1$ (and for the compact integration hidden in $f_2$). Also $|\widehat\Phi_{d_hg}(r,v)|=|\det h|^{-1}|\widehat\Phi_g(rh,h^{-1}v)|$ is bounded by a Schwartz function of $(h^{-1}v,rh)$ times $|\det h|^{-1}$. After combining $t$ and $h$ as in Step 2 of the proof of \Cref{thm:linear-unfold} (which only uses $|\mu|=1$ and $|\eta|=1$), the absolute value of the integrand is dominated by $\Psi(z\cdot h)|\det h|^{2\Re s-1}$ with $\Psi$ Schwartz on $\mathcal V_1(\A)$ ($M_1=B$). \Cref{thm:A-sum} (restricted to $B=A+vr-I$, $A,B$ invertible) shows that the sum over rational orbits of the integrals over $\GL_n(\A)$ converges, locally uniformly in $s$, also after inserting powers of $\log|\det h|$. Tonelli's theorem justifies the disintegration of $[\GL_n]$ along its center, the unfolding of the free orbits and the factorization; differentiation under the integral sign (or Morera) gives holomorphy.
\end{proof}

\begin{remark}
The product formula $\prod_v|\det A^h|_v^{-\alpha}=1$ is used only for \emph{rational} invertible $A$ and $B$. Locally, the factor $|\det C|^{-\alpha}$ in \Cref{lem:A-majorant} is not bounded as $\det C\to0$, so the estimate is not uniform near the singular locus. This is why the singular frequencies $\det A\det B=0$ are not covered by \Cref{thm:A-regular} and remain boundary terms.
\end{remark}

\begin{proposition}\label{prop:A-orth}
In \Cref{prop:orth-unfold} the regular part converges absolutely, only finitely many rational regular orbits contribute, each local orbital integral converges absolutely, and the global orbital integrals factor. If moreover, at one place $v_1$, the support of $f_{v_1}$ consists of pairs $(x,y)$ with $xy^{-1}$ regular, then the full expression $I(f)$ converges absolutely and equals $I^{\rm reg}(f)$.
\end{proposition}
\begin{proof}
$F_W(x)=\int|f(xu,u)|\,du$ is compactly supported, since $u$ and $xu$ range over compact sets. For a regular $g$, let $E_g=[e,ge,\dots,g^{2n}e]$. If $h\in\Ort(W)$ (extended by $e\mapsto e$) and $g'=h^{-1}gh$, then $E_{g'}=h^{-1}E_g$, i.e.\ $h^{-1}=E_{g'}E_g^{-1}$; so $g'$ in a compact set forces $h$ into a compact set (also $h$ is bounded since $h$ preserves the form). Adelically, $g$ and $E_g^{-1}$ are rational and integral at almost all places, so $h$ is in the standard maximal compact subgroup outside a finite set of places. This proves convergence of each orbit, adelically and locally. The characteristic polynomial of $g$ and the moments $(e,g^ie)$, $i\le 2n$, are polynomial invariants which on the cyclic locus determine the $\Ort(W)(k)$-orbit (they give the matrix of $g$ in the basis $g^ie$ and, since $(g^ie,g^je)=(e,g^{j-i}e)$ and $(e,g^{-m}e)=(g^me,e)$ for $g$ orthogonal, also the Gram matrix in this basis); on a compact adelic set they take finitely many rational values, so only finitely many rational $\Ort(W)$-orbits contribute. Tonelli then justifies the unfolding of the nonnegative kernel, and factorization follows for factorizable $f$. For the last assertion: the collapsed element of a term $f(\gamma\cdot(h_1,h_2))$ is $\widehat h_1^{-1}\gamma_1\widehat\gamma_2^{-1}\widehat h_1$, and regularity is a polynomial condition defined over $k$, so if $\gamma_1\gamma_2^{-1}$ is not regular then $f_{v_1}$ vanishes at the $v_1$-component of every such term. Hence all nonregular terms of the kernel vanish identically, and $I(f)=I^{\rm reg}(f)$ with absolute convergence.
\end{proof}

\section{Full kernels}

\begin{proposition}\label{prop:A-full}
Fix a finite place $v_0$ of odd residual characteristic and compact open sets $\mathcal U_A,\mathcal U_B\subset\GL_n(k_{v_0})$, $\mathcal U_r\subset k_{v_0}^{1\times n}$, $\mathcal U_v\subset k_{v_0}^n$ such that $A,B$ are invertible and the Hankel matrix of $(A,v,r)$ is invertible on $\mathcal U_A\times\mathcal U_B\times\mathcal U_r\times\mathcal U_v$ (here $B$ is independent of $A$). Replace the $v_0$-components of $\Theta^\pm$ by $\int\rho_\pm(X)\theta^\pm(\tilde u(X))\xi^\pm_{v_0}\,dX$ with $\rho_\pm\in C_c^\infty(\Mat_n(k_{v_0}))$ whose signed Fourier transforms
\[
\widehat\rho_\pm(C)=\int_{\Mat_n(k_{v_0})}\rho_\pm(X)\psi_{v_0}\bigl(\pm\tr(CX)\bigr)\,dX
\]
are supported in $\mathcal U_A$, resp.\ $\mathcal U_B$; let $\Phi_{v_0}$ be the inverse partial Fourier transform of a function in $C_c^\infty(\mathcal U_r\times\mathcal U_v)$; and let $f_{1,v_0},f_{2,v_0}$ be supported in a small compact open subgroup fixing these vectors and $\Phi_{v_0}$. Then the full expression $J(f_1,f_2,\Phi;s)$ converges absolutely for all $s$, is entire, and equals $J^{\rm reg}(s)$: the axes, the singular frequencies and the noncyclic terms all vanish.
\end{proposition}
\begin{proof}
At $v_0$, right convolution by $\rho_\pm$ multiplies the Fourier coefficient of frequency $C$ evaluated at $d_h$ by $\widehat\rho_\pm(C^h)$, since $d_hu(X)=u(hXh^{-1})d_h$. The allowed translations at $v_0$ fix the filtered vectors, so the same support restrictions hold at $d_hg$. Thus the surviving coefficients satisfy $A^h\in\mathcal U_A$ and $B^h\in\mathcal U_B$; in particular both are invertible. The conditions on $\Phi_{v_0}$ give $\Phi_{v_0}(0)=0$ and $\widehat\Phi_{v_0}(r,0)=\widehat\Phi_{v_0}(0,v)=0$, so the axes vanish; and the support conditions persist under the allowed translations. So every surviving term has $(A^h,B^h,trh,t^{-1}h^{-1}v)$ in a regular cyclic set at $v_0$. Now bound the absolute value of the whole expanded expression as in \Cref{thm:A-regular}, using \Cref{thm:A-sum} with $B$ independent of $A$. This justifies all interchanges, before the $U$-integration imposes $A-B+vr=I_n$. It also verifies the two hypotheses of \Cref{lem:kernels}(c). Indeed, by the triangle inequality the absolute value $|\Xi_{f_2}(ag)E(ag,s)|$ is at most the sum of the absolute values of the terms of its expansion, so the majorant bounds $\int\int|f_1(g^{-1})\Xi_{f_2}(ag)E(ag,s)|\,da\,dg$. And unfolding the absolute kernels as in the proof of \Cref{lem:kernels} bounds the absolute triple integral defining $J$ by
\[\int_{\GL_{2n}(\A)^2}|f_1(y^{-1})f_2(b)|\int_{[S]_Z}|E(ay,s)||\Xi(ayb)|\,da\,dy\,db,\]
which is finite because the majorant is uniform for $y,b$ in the compact supports of the test functions. The Mellin integral of the original theta series is entire because both constant terms in Poisson summation vanish. This does not mean that the Mellin integral may be interchanged, term by term, with the original lattice sum over $w\in k^{1\times2n}$ at $s=\frac12$; the interchanges are made in the order of \Cref{lem:poisson}: partial Poisson summation for fixed $t$ first, and then integration of the localized expansion, which has the absolute majorant above.
\end{proof}

\begin{remark}\label{rem:A-full-nonempty}
	The class of Proposition \ref{prop:A-full} is not empty, and its localizations can be chosen around a regular rational datum of the prescribed discriminant. Fix \(\delta\in k^\times\) with \(\eta=\eta_\delta\). Let \(A_*\) be the companion matrix of \(T^n+\delta\), with subdiagonal entries \(1\), and put
	\[
	v_*=e_1,\qquad
	r_*=\delta e_n^t,\qquad
	B_*=A_*+v_*r_*-I_n.
	\]
	If \(N\) denotes the nilpotent matrix with subdiagonal entries \(1\) and all other entries zero, then
	\[
	A_*=N-\delta e_1e_n^t,\qquad B_*=N-I_n.
	\]
	Consequently,
	\[
	\det A_*=(-1)^n\delta,\qquad
	\det B_*=(-1)^n,\qquad
	\det A_*\det B_*=\delta.
	\]
	The Hankel matrix
	\[
	\bigl(r_*A_*^{i+j}v_*\bigr)_{0\leq i,j<n}
	\]
	is antitriangular with antidiagonal entries \(\delta\), and is therefore nondegenerate. Thus \(z_*=(A_*,v_*,r_*)\) belongs to \(X_\delta^{\mathrm{reg}}(k)\). By continuity, there are compact open neighborhoods \(U_A\ni A_*\), \(U_B\ni B_*\), \(U_r\ni r_*\), and \(U_v\ni v_*\) at \(v_0\) satisfying the conditions of Proposition \ref{prop:A-full}; they may also be chosen so that \(\det A\det B\in\delta k_{v_0}^{\times2}\) throughout \(U_A\times U_B\).

	The filtered vectors can be taken nonzero: $\ell\circ\theta^+(s(m_A))$ composed with the filter $\xi\mapsto\int\rho_+(X)\theta^+(\tilde u(X))\xi\,dX$ equals $\widehat\rho_+(A)\,\ell\circ\theta^+(s(m_A))$, so if $\widehat\rho_+(A_*)\ne0$ and $\xi^+_{v_0}$ is chosen with $\ell(\theta^+(s(m_{A_*}))\xi^+_{v_0})\ne0$, then the filtered local vector has a nonzero full-rank coefficient at $A_*$; the same argument applies to the negative-character coefficient at $B_*$. At the archimedean places and the remaining exceptional finite places, choose suitable local data, with nonzero full-rank theta coefficients at \(A_*\) and \(B_*\). At all other sufficiently good finite places, take the unramified basic data. With these choices, \Cref{lem:factor} shows that the resulting
	filtered global forms $\Theta^+$ and $\Theta^-$ have nonzero
	Fourier coefficients at $A_*$ and $B_*$, respectively. 
\end{remark}

\begin{remark}\label{rem:divergence}
Without hypotheses such as those of \Cref{prop:A-orth} (regular support at one place) or \Cref{prop:A-full}, the full expressions can diverge. For $n=1$ and $W$ the split plane, $\SO(W)=\mathbb G_m$ and $[\Ort(W)]$ has infinite volume. The orbit of $(1,1)$ has stabilizer the diagonal $\Ort(W)$ in $\mathcal H_W$; for $f\ge0$ with $\int_{\SO(W)(\A)}f(\widehat u,u)\,du>0$, its contribution is therefore $\vol([\Ort(W)])\int_{\SO(W)(\A)}f(\widehat u,u)\,du=\infty$. Its collapsed element is $g=1$, which is not regular.
\end{remark}

\chapter{Technical lemmas}\label{app:technical}

\section{The metaplectic cover and the spherical coefficient}

\begin{lemma}\label{lem:unip-kappa}
If $u\in\GL_j(\cO)$ is unipotent, its canonical lift $\tilde u$ equals $\kappa(u)$.
\end{lemma}
\begin{proof}
Both lie over $u$, so $\kappa(u)=\epsilon\tilde u$ with $\epsilon=\pm1$, and they commute. The canonical lift is a homomorphism on a unipotent group containing $u$, so $\tilde u^m$ is the canonical lift of $u^m$; also $\kappa(u)^m=\kappa(u^m)$. Writing $u=1+N$, $u^{p^k}=\sum_i\binom{p^k}iN^i\to1$ as $k\to\infty$ ($p$ the residual characteristic, $N^j=0$). By continuity of both splittings, $\tilde u^{p^k}\to1$ and $\kappa(u)^{p^k}\to1$, hence $\epsilon^{p^k}=\kappa(u)^{p^k}\tilde u^{-p^k}\to1$. As $p$ is odd, $\epsilon^{p^k}=\epsilon$, so $\epsilon=1$.
\end{proof}

\begin{lemma}\label{lem:compact}
Assume the data are unramified.
\begin{enumerate}[label=(\alph*)]
\item $\beta(\kappa(d_k))=1$ for $k\in\GL_j(\cO)$.
\item If $\tilde C$ lies over $C$ with $\val\det C$ even, then $\cB(\kappa(k_1)\tilde C\kappa(k_2))=\cB(\tilde C)$ for $k_1,k_2\in\GL_j(\cO)$.
\item $\cB^-=\cB$; in particular $\cW^-=\cW$.
\item For every $C$, the function $h\mapsto\cB(\T_h\tilde C)$, equivalently $h\mapsto\epsilon(C,h)\cW(C^h)$, is right $\GL_j(\cO)$-invariant.
\end{enumerate}
\end{lemma}
\begin{proof}
(a) $\ell(\theta(\kappa(d_k))\xi^\circ)=\ell(\xi^\circ)=1$, and the left side is $\beta(\kappa(d_k))\ell(\xi^\circ)$.

(b) By \ref{C:compact}, $\kappa(k_2)\oplus1=\kappa(\diag(k_2,I))$ fixes $\xi^\circ$, which gives right invariance. For left invariance write $\kappa(k_1)\oplus1=\kappa(\diag(k_1,I))=\kappa(d_{k_1})\kappa(\diag(I,k_1^{-1}))$. By \ref{C:block}, $s(\diag(I,K))s(\diag(C,I))=\hil{\det C}{\det K}s(\diag(C,I))s(\diag(I,K))$, so lifts of $\diag(I,k_1^{-1})$ and $\diag(C,I)$ commute up to $\hil{\det C}{\det k_1}$, which is $1$ because $\det k_1\in\cO^\times$ and $\val\det C$ is even. Hence $\ell(\theta(\kappa(k_1)\tilde C\oplus1)\xi^\circ)=\beta(\kappa(d_{k_1}))\ell(\theta(\tilde C\oplus1)\theta(\kappa(\diag(I,k_1^{-1})))\xi^\circ)$, and (a) and the invariance of $\xi^\circ$ finish the proof.

(c) If $C$ is not integral, or has an odd Smith exponent, both sides vanish by \ref{S}. Otherwise $C=k_1tk_2$ with $k_i\in\GL_j(\cO)$, $t=\diag(\varpi^{2d_1},\dots,\varpi^{2d_j})$, $d_1\ge\dots\ge d_j\ge0$, and $\tilde C=\zeta\kappa(k_1)s(t)\kappa(k_2)$ for a sign $\zeta$. By (b), $\cB^\pm(\tilde C)=\zeta\cB^\pm(s(t))=\zeta q^{\sum(i-1)d_i}$.

(d) $\T_{hk}\tilde C=\T_k(\T_h\tilde C)=\hil{\det C}{\det k}\kappa(k)^{-1}\T_h(\tilde C)\kappa(k)$. If $\val\det C$ is even the symbol is $1$ and (b) applies; if it is odd, $\cB$ vanishes on all elements over conjugates of $C$. The equivalence with $\epsilon(C,h)\cW(C^h)$ is \eqref{eq:summand-T}.
\end{proof}

\begin{lemma}[Unit blocks]\label{lem:unitblock}
Let $\tilde X$ lie over $C\in\GL_a(F)$ with $\val\det C$ even. Then $\cB_{a+b}(\tilde X\oplus1)=\cB_{a+b}(1\oplus\tilde X)=\cB_a(\tilde X)$, where $\cB_m$ is the function \eqref{eq:Wcoef} for $\theta_{2m}$.
\end{lemma}
\begin{proof}
The Smith exponents of $\diag(C,I_b)$ and $\diag(I_b,C)$ are those of $C$ together with $b$ zeros, so if one side vanishes by \ref{S} all do. Otherwise write $\tilde X=\zeta\kappa(k_1)s(t)\kappa(k_2)$ as in \Cref{lem:compact}(c). Since $\tilde Y\mapsto\tilde Y\oplus1$ is a homomorphism compatible with $\kappa$ and $s$, $\tilde X\oplus1=\zeta\kappa(k_1\oplus I)s(t\oplus I)\kappa(k_2\oplus I)$, and by \Cref{lem:compact}(b) and \ref{S}, $\cB_{a+b}(\tilde X\oplus1)=\zeta q^{\sum(i-1)d_i}=\cB_a(\tilde X)$ (the appended exponents are $0$ and come last). Similarly $\cB_{a+b}(1\oplus\tilde X)=\zeta\cB_{a+b}(s(I\oplus t))$. Let $w$ be the permutation (as in \ref{C:weyl}) with $w(I\oplus t)w^{-1}=t\oplus I$. By \ref{C:weyl}, $s(w)s(I\oplus t)s(w)^{-1}=\sigma(w,I\oplus t)s(t\oplus I)$, and every symbol $\hil{-t_k}{t_i}$ is $1$ because all entries are even powers of $\varpi$ ($\hil{x}{y^2}=1$). Since $s(w)=\kappa(w)$, \Cref{lem:compact}(b) gives $\cB(s(I\oplus t))=\cB(s(t\oplus I))$.
\end{proof}

\begin{lemma}[Positive Hessenberg value]\label{lem:hessenberg}
Let $C\in\GL_j(F)$ be integral and upper Hessenberg with all subdiagonal entries equal to $1$. Then $\cW(C)=\one_{\val\det C\in2\Z}$. The same holds for $\cW^-$.
\end{lemma}
\begin{proof}
If $\val\det C$ is odd, $\cW(C)=0$ by \ref{S}. Assume it is even and put $\delta=\det C$. Let $U_C=[e_1,Ce_1,\dots,C^{j-1}e_1]$. The Hessenberg condition gives $C^ie_1=e_{i+1}+(\text{combination of }e_1,\dots,e_i)$ with integral coefficients, so $U_C\in N_j(\cO)$. The matrix $U_C^{-1}CU_C$ is the companion matrix with last column $z=U_C^{-1}C^je_1$, which is integral, and $z_1=(-1)^{j-1}\delta$. Let $w_j=w_{\alpha_1}\cdots w_{\alpha_{j-1}}$, so that $w_je_a=e_{a+1}$ for $a<j$ and $w_je_j=(-1)^{j-1}e_1$, let $t_\delta=\diag(1,\dots,1,\delta)$ and $V=I+\sum_{i<j}z_{i+1}E_{ij}$. Checking on the basis vectors one finds $U_C^{-1}CU_C=w_jt_\delta V$, i.e.\
\[C=U_C\,w_jt_\delta\,V',\qquad V'=VU_C^{-1}\in N_j(\cO).\]
The only inversions of $w_j$ are the pairs $(i,j)$, $i<j$, so by \ref{C:weyl} $\sigma(w_j,t_\delta)=\prod_{i<j}\hil{-\delta}{1}=1$. With \ref{C:unip} we get $s(C)=s(U_C)s(w_j)s(t_\delta)s(V')$, and the outer three factors are $\kappa$-lifts by \ref{C:compact}. By \Cref{lem:compact}(b), $\cW(C)=\cB(s(t_\delta))$. Write $\delta=u\varpi^{2d}$, $u\in\cO^\times$; by \ref{C:torus}, $s(t_\delta)=s(t_u)s(t_{\varpi^{2d}})$ and $s(t_u)=\kappa(t_u)$, so $\cW(C)=\cB(s(t_{\varpi^{2d}}))$. Finally $s(w_j)s(t_{\varpi^{2d}})s(w_j)^{-1}=s(\diag(\varpi^{2d},1,\dots,1))$ by \ref{C:weyl} (the symbol is $1$ as before), and $\cB$ of this is $q^0=1$ by \ref{S}.
\end{proof}

\section{Selfdual lattices}

\begin{lemma}\label{lem:selfdual}
Let $(W,Q)$ be a nondegenerate quadratic space over $F$ (odd residual characteristic). Then $\Ort(W)$, and also $\SO(W)$, acts transitively on the set of selfdual lattices in $W$.
\end{lemma}
\begin{proof}
A selfdual lattice $L$ has an orthogonal basis consisting of vectors of unit norm (diagonalize the unimodular form over $\cO$, using $2\in\cO^\times$). If $L,L'$ are selfdual, their Gram determinants are units with the same class in $F^\times/F^{\times2}$ (the discriminant of $W$), hence the same class in $\cO^\times/\cO^{\times2}$. Unimodular quadratic $\cO$-lattices of the same rank and discriminant are isometric (reduce modulo $\varpi$, where nondegenerate forms over the residue field are classified by rank and discriminant, and lift by Hensel's lemma). An isometry $L\to L'$ extends to an element of $\Ort(W)$. Composing with the reflection in a unit vector of $L$ (which preserves $L$) we may make its determinant $1$.
\end{proof}

\section{Local constancy}

\begin{lemma}\label{lem:hess-continuous}
Let $\mathcal H_j$ be the set of invertible upper Hessenberg matrices with subdiagonal entries $1$. Then $s|_{\mathcal H_j}$ is continuous.
\end{lemma}
\begin{proof}
The factorization $C=U_Cw_jt_{\det C}V'$ in the proof of \Cref{lem:hessenberg} holds on all of $\mathcal H_j$ (integrality was not used for it), with $U_C,V'$ polynomial in $C$, and $s(C)=s(U_C)s(w_j)s(t_{\det C})s(V')$. The canonical lift is continuous on $N_j(F)$, and $\delta\mapsto s(t_\delta)$ is continuous: on the open set $\delta\in u\varpi^a$, $u\in\cO^\times$, it equals $s(t_{\varpi^a})\kappa(t_u)$.
\end{proof}

\begin{proposition}\label{prop:local-constancy}
The functions $J^e$, $I^e$ (on the set of even data of size $n$) and $J^o$, $I^o$ (on the set of odd data of size $n$) are locally constant as functions of all coordinates of the datum, not merely along a fixed spectral fibre \footnote{Here the domains are the full enhanced regular loci of Definitions \ref{def:even} and \ref{def:odd}: the relevant Hankel determinant and the determinant of the full matrix $S$ remain nonzero.	No semisimplicity of the spectral matrix is assumed, and $S_0$ may be singular in the odd case. Thus the assertion	includes points outside the semisimple charts of Section \ref{sec:fibres};	it makes no assertion across the boundary where enhanced	regularity fails.}.
\end{proposition}
\begin{proof}
\emph{Even linear.} Fix a datum $p_0=(S,x,r)$. By \Cref{lem:models-basic}(c) the contributing lattices lie between $C_x(S)\cO^n$ and $R_r(S)^{-1}\cO^n$. These matrices depend continuously on the datum and are invertible, and $g\mapsto g\cO^n$ is locally constant; hence on a neighbourhood $\mathcal U$ of $p_0$ both lattices are constant, and the contributing lattices belong to a fixed finite list $h_1\cO^n,\dots,h_N\cO^n$. For each $i$ the summand is
\[\Delta(S,x)^{-1}\epsilon(S,h_i)\cW(S^{h_i})\one(h_i^{-1}x)\one(rh_i)=\cB\bigl(\T_{h_i}(\tilde S)\bigr)\one(h_i^{-1}x)\one(rh_i),\qquad\tilde S=\Delta(S,x)^{-1}s(S).\]
Now $\tilde S=\T_{C_x^{-1}}\bigl(s(S^{C_x})\bigr)$, because $\T_{C_x^{-1}}(s(S^{C_x}))=\epsilon(S^{C_x},C_x^{-1})s(S)$ and $\epsilon(S^{C_x},C_x^{-1})=\epsilon(S,C_x)^{-1}$ by the cocycle relation. The companion matrix $S^{C_x}$ lies in $\mathcal H_n$ and depends continuously on the datum, so by \Cref{lem:hess-continuous} and continuity of $(g,\tilde X)\mapsto\T_g\tilde X$, the lift $\tilde S$ depends continuously on the datum. Finally $\cB$ is locally constant on the cover, being a matrix coefficient of a smooth representation times a locally constant function. Hence each summand, and the finite sum, is locally constant on $\mathcal U$. (The individual factors $\Delta(S,x)$ and $\cW(S^{h_i})$ need not be locally constant; see \Cref{rem:wall}.)

\emph{Odd linear.} Contributing lattices satisfy $S_0\Lambda\subset\Lambda$, $b\in\Lambda$, $c\Lambda\subset\cO$, hence lie between $C_b(S_0)\cO^{n-1}$ and $R_c(S_0)^{-1}\cO^{n-1}$, which are invertible by regularity of $(S_0,b,c)$; no inverse of $S_0$ is used. The same argument applies with $\tilde S=\Delta(S)^{-1}s(S)=\T_{C_f^{-1}}(s(S^{C_f}))$, $S^{C_f}$ being a companion matrix.

\emph{Orthogonal.} Work in the fixed coordinate space $F^d$, where $d=\dim W$. The Gram matrix $Q$ ($Q^e$ or $Q^o$), the operator $C$ and the vector $w$ depend continuously on the datum (by \Cref{lem:Q0}, $H=C_x^{-\mathsf t}(R_rC_x)C_x^{-1}$, with the analogous formula for $(S_0,b,c)$ in the odd case). Every counted lattice lies between
\[
 M=E\cO^d,\qquad M^{\#_Q}=Q^{-1}E^{-\mathsf t}\cO^d,
 \qquad E=[w,Cw,\dots,C^{d-1}w].
\]
Both bounding lattices are constant on a sufficiently small neighbourhood, since $E$ and $Q$ remain invertible. Hence the candidate lattices form a fixed finite list (possibly empty). For a fixed candidate $L=g_L\cO^d$, its three defining conditions are
\[
\begin{aligned}
 L=L^{\#_Q}&\iff g_L^{\mathsf t}Qg_L\in\GL_d(\cO),\\
 CL\subseteq L&\iff g_L^{-1}Cg_L\in\Mat_d(\cO),\\
 w\in L&\iff g_L^{-1}w\in\cO^d.
\end{aligned}
\]
The sets on the right are open and closed, so their indicators are locally constant. Shrinking the neighbourhood for the finitely many candidates proves local constancy of the count, even though $Q$ varies with the datum. The neighbourhood and finite list may depend on the point; no uniform spectral neighbourhood or spectral cutoff is required.
\end{proof}

\begin{remark}\label{rem:wall}
The cyclic normalization is essential here. Let $\nu\in\cO^\times$ be a nonsquare and $C_N=\left(\begin{smallmatrix}\nu&0\\\varpi^N&1\end{smallmatrix}\right)$ for $N\ge1$, $C_\infty=\diag(\nu,1)$, $x=(1,1)^{\mathsf t}$ (cyclic for all these, as $\nu\not\equiv1$). A Kubota computation gives $\cW(C_N)=\hil{\varpi^N}{\nu}=(-1)^N$ while $\cW(C_\infty)=1$, so the raw coefficient is not continuous at $C_\infty=\lim C_N$; but $\epsilon(C_N,C_x(C_N))=(-1)^N$ as well, and the normalized summand is $1$ for all $N$.
\end{remark}

\section{Compactness of the outer integration}

\begin{lemma}\label{lem:compact-outer}
Let $S\in\GL_n(F)$ be regular semisimple, $K=F[S]$, and $\xi\in\theta_{2n}$. The set of $g\in K^\times\backslash\GL_n(F)$ with $\cW_\xi(S^g)\ne0$ is relatively compact. In the odd case, for $S_0$ regular semisimple with centralizer $K^\times$, the set of $g\in K^\times\backslash\GL_{n-1}(F)$ for which $\cW_{\xi}(S^{\hat g})\ne0$ for some $S$ with upper left block $S_0$ is relatively compact, uniformly for $\xi$ in an orbit of $\iota(\Mp_2(F))$.
\end{lemma}
\begin{proof}
Let $\xi$ be fixed by $\tilde u(X)$ for $X\in\varpi^N\Mat_n(\cO)$. Then $\ell(\theta(s(m_C))\xi)=\ell(\theta(\tilde u(CX))\theta(s(m_C))\xi)=\psi(\tr CX)\ell(\theta(s(m_C))\xi)$, so $\cW_\xi(C)\ne0$ forces $C\in\varpi^{-N}\Mat_n(\cO)$. If no \(g\) contributes, the assertion is immediate. Otherwise, for one contributing \(g_0\), the matrix
\[
g_0^{-1}(\varpi^NS)g_0
\]
is integral. Hence \(\varpi^NS\) is integral over \(\mathcal O\); since it generates \(K\) over \(F\), the ring
\[
R=\mathcal O[\varpi^NS]
\]
is an order in \(K\). Every contributing lattice \(g\mathcal O^n\) is stable under \(R\).

Identify $F^n$ with $K$ by a cyclic vector of $S$; then $g\cO^n$ is an $R$-stable lattice in $K$. Its $\cO_K$-span is $\lambda\cO_K$ for some $\lambda\in K^\times$; replacing $g$ by $\lambda^{-1}g$, the lattice lies between $\mathfrak f\cO_K$ and $\cO_K$, where $\mathfrak f$ is the conductor of $R$. There are finitely many such lattices. In the odd case apply the same argument to the elements $\tilde u(X)$ with $X$ supported in the upper left $(n-1)\times(n-1)$ block (the group $N_0$): $\tr(S^{\hat g}X)=\tr(S_0^gX_0)$, so $S_0^g$ is bounded independently of $b,c,d$. Since $\iota(\Mp_2(F))$ centralizes $N_0$, the subgroup of $N_0$ fixing $\xi$ fixes all $\theta(\iota(g))\xi$.
\end{proof}

\section{Estimates for the binary kernel}

\begin{lemma}\label{lem:shell}
Let $K$ be a nonarchimedean field of odd residual characteristic, $\eta$ a quadratic character of $K^\times$, $\mathcal C\subset K^\times$ compact and $R\ge0$. For $s\in\mathcal C$ and $\val(a)\ge-R$, every partial annular integral
\[\mathcal K_\eta^{[p,q]}(s,a)=\int_{p\le\val t\le q}\eta(t)\psi_K(at+s/t)\,d^\times t\]
satisfies $|\mathcal K_\eta^{[p,q]}(s,a)|\le C(1+\val(a)^+)$ with $C$ independent of $p,q$; for $s,a\ne0$ it stabilizes as $p\to-\infty$, $q\to\infty$.
\end{lemma}
\begin{proof}
Assume first that $\psi_K$ has conductor $\cO_K$. Since $\eta$ is quadratic and the residual characteristic is odd, $\eta$ is trivial on $1+\varpi_K\cO_K$. Hence $\int_{\cO_K^\times}\eta(u)\psi_K(\lambda u)\,du=0$ if $\val\lambda<-1$: on each coset $u_0+\varpi_K\cO_K$, $\eta$ is constant and the additive integral vanishes. The same holds with $u$ replaced by $u^{-1}$. Put $t=\varpi_K^ku$, $S=\val s$, $A=\val a$. If $A+k\ge0$ and $S-k<-1$, the term $at$ is trivial under $\psi_K$ and the shell integral vanishes; if $S-k\ge0$ and $A+k<-1$ the same holds with the roles exchanged. As $S$ is bounded and $A\ge-R$, only shells with $k$ in an interval of length $\le C(1+A^+)$ can contribute, each of bounded volume. Stabilization follows. A general conductor is reduced to this by scaling $s,a$, and different selfdual measures change only the constants.
\end{proof}

\begin{lemma}\label{lem:L1}
For $\phi\in\cS(K^2)$, $\int_K|L\phi(a)|\,da\le\int_{K^2}|\phi|$, $L\phi$ has compact support, and $|L\phi(a)|\le C_\phi(1+\val(a)^+)$. For $\varphi\in\cS(E)$, $N_\beta\varphi$ is an $L^1$ function with compact support and $|N_\beta\varphi(a)|\le C_\varphi(1+\val(a)^+)$. Consequently the pairings with $\mathcal K_\eta(s,\cdot)$ converge absolutely for $s\ne0$, locally uniformly, and a function in the image of $L$ or $N_\beta$ that vanishes on $K^\times$ is zero as a distribution. For an \'etale algebra the same holds factor by factor.
\end{lemma}
\begin{proof}
The change of variables $(u,a)\mapsto(u,a/u)$ has Jacobian $|u|^{-1}$, giving the $L^1$ bound. If $\phi$ is supported in $\varpi^{-R}\cO\times\varpi^{-R}\cO$, the integrand of $L\phi(a)$ is nonzero only for $-R\le\val u\le\val a+R$, which gives the logarithmic bound; the support of $L\phi$ lies in the compact image of ${\rm supp}\,\phi$ under $(u,v)\mapsto uv$. On the norm side, if $E=K\times K$ this is the same computation. If $E$ is a field, the differential of $q_\beta$ at $z\ne0$ is the nonzero linear form $u\mapsto2q_\beta(z,u)$, so locally one coordinate can be replaced by $q_\beta$, which shows absolute continuity away from $0$; near $0$, $\varphi$ is constant and the pushforward of the measure of a small ball is bounded by homogeneity (a shell $\varpi_E^k\cO_E^\times$ and its image shell have volumes scaling by the same power of $q_K$). The pairing with the kernel is dominated by $\sum_{m\ge0}q_K^{-m}(1+m)^2<\infty$ near $0$. Finally $K\smallsetminus K^\times=\{0\}$ has measure zero.
\end{proof}

\section{The rank one oscillator and the rank two exceptional representation}\label{sec:oscillator}

Let $\mathcal F_\psi f(x)=\int_Ff(y)\psi(2xy)\,dy$ with $dy$ selfdual for $\psi(2xy)$, let $\gamma(\psi_a)$ be the Weil index of $ax^2$, so that
\begin{equation}\label{eq:gauss1}
\int_F\psi(at^2+2bt)\,dt=\gamma(\psi_a)|a|^{-1/2}\psi(-b^2/a),
\end{equation}
and put $\mu_\psi(a)=\gamma(\psi_a)/\gamma(\psi)$. Then $\gamma(\psi_{a^{-1}})=\gamma(\psi_a)$, $\gamma(\psi_a)\gamma(\psi_{-a})=1$ and $\mu_\psi(a)^2=\hil a{-1}$.

\begin{lemma}\label{lem:oscillator}
Let $\rho_\psi$ be the oscillator (Weil) representation of $\Mp_2(F)$ on $\cS(F)$ for the form $x^2$, identified with the cover of \S\ref{sec:actions} (the identification is unique because $\SL_2(F)$ is perfect). Then
\[\rho_\psi(\tilde n(t))f(x)=\psi(tx^2)f(x),\qquad\rho_\psi(s(w_0))=\gamma(\psi)\mathcal F_\psi,\]\[\rho_\psi(s(h(a)))f(x)=\mu_\psi(a)^{-1}|a|^{1/2}f(ax).\]
\end{lemma}
\begin{proof}
$\rho_\psi(\tilde n(t))$ is a homomorphic lift of $n(t)$, hence equals $N(t)f=\psi(tx^2)f$ up to a character $F\to\{\pm1\}$, which is trivial since $F$ is $2$-divisible. By Weil's construction, $W=\rho_\psi(s(w_0))=\pm\gamma(\psi)\mathcal F_\psi$. Let $l(t)=\left(\begin{smallmatrix}1&0\\t&1\end{smallmatrix}\right)$. On the lower unipotents the BLS cocycle is trivial: $\sigma(l(t),l(t'))=\hil{(t+t')/t}{(t+t')/t'}=1$ by the Steinberg relation, so $s(l(t))$ is the canonical lift, and $s(w_0)\tilde n(t)s(w_0)^{-1}=s(l(-t))$. A direct computation with the cocycle (here $x(g)$ denotes the lower left entry of $g$ if it is nonzero and the lower right entry otherwise, and $\sigma(g,g')=\hil{x(gg')/x(g)}{x(gg')\det g/x(g')}$ \cite[(1.2)]{KM}) gives
\[s(n(a))s(l(-a^{-1}))s(n(a))=s(h(a)w_0),\qquad s(h(a))s(w_0)=\hil{-1}{a}s(h(a)w_0).\]
Hence $\rho_\psi(s(h(a)w_0))=N(a)\,W N(a^{-1})W^{-1}\,N(a)$. By \eqref{eq:gauss1}, $\mathcal F_\psi N(b)\mathcal F_\psi^{-1}$ has kernel $\gamma(\psi_b)|b|^{-1/2}\psi(-(x-y)^2/b)$, so with $b=a^{-1}$ the operator $N(a)\mathcal F_\psi N(a^{-1})\mathcal F_\psi^{-1}N(a)$ has kernel $\gamma(\psi_a)|a|^{1/2}\psi(2axy)$. At $a=1$ this says $W=N(1)WN(1)W^{-1}N(1)=\gamma(\psi)\mathcal F_\psi$, which fixes the sign of $W$ (the right side does not change if $W$ is replaced by $-W$). For general $a$, $\rho_\psi(s(h(a)))=\hil{-1}a\rho_\psi(s(h(a)w_0))W^{-1}=\hil{-1}a\mu_\psi(a)|a|^{1/2}D_a$, with $D_af(x)=f(ax)$; and $\hil{-1}a\mu_\psi(a)=\mu_\psi(a)^{-1}$.
\end{proof}

\paragraph{General quadratic forms.}
For a nondegenerate symmetric bilinear form $Q$ on a finite dimensional $F$-space $U$, write $Q(u)=Q(u,u)$ and take $du$ selfdual for $\psi(2Q(u,u'))$, as in \S\ref{sec:actions}. The Weil index $\gamma(Q)$ is normalized by the Gaussian identity, understood as an oscillatory integral:
\begin{equation}\label{eq:gauss}
\int_U\psi\bigl(tQ(u)+2Q(u,u')\bigr)\,du=\gamma(tQ)\,|t|^{-\dim U/2}\,\psi\bigl(-Q(u')/t\bigr)\qquad(t\in F^\times).
\end{equation}
If $Q\simeq\langle a_1,\dots,a_m\rangle$, then $\gamma(Q)=\prod_i\gamma(\psi_{a_i})$. The representation $\omega_Q$ is the tensor product of the one dimensional oscillator representations $\rho_{\psi_{a_i}}$ of \Cref{lem:oscillator}; thus their formulas give \eqref{eq:weil} for the specified lifts. This supplies the normalization used in \S\ref{sec:actions}; the Gaussian identity is used in the binary norm-kernel calculation of \S\ref{sec:binary}.

\begin{lemma}\label{lem:theta2}
Let $G^\square=\{g\in\GL_2(F):\det g\in F^{\times2}\}$ and $\tilde G^\square$ its inverse image. By \cite[(2.13)--(2.14), Proposition 2.16]{Takeda}, $\theta_2\simeq R=\mathrm{Ind}_{\tilde G^\square}^{\tGL_2}\rho$, where $\rho$ is an extension of the even part $\rho_\psi^+$ of $\rho_\psi$ to $\tilde G^\square$; the $\psi$-Whittaker functional is $\lambda(f)=f(1)(1)$, and in the unramified case the spherical vector satisfies $f^\circ(1)=\one_\cO$. Let $t_a=s(\diag(a,1))$, $\mathcal D$ a set of representatives of $F^\times/F^{\times2}$, and $\phi_\delta=|\delta|^{-1/4}f(t_\delta)\in\cS(F)^+$. Then $f\mapsto(\phi_\delta)_{\delta\in\mathcal D}$ is an isomorphism onto $\bigoplus_\delta\cS(F)^+$, and
\begin{enumerate}[label=(\roman*)]
\item $\Mp_2(F)$ preserves each component and acts on it by $\tilde n(t)\mapsto\psi(t\delta y^2)$ and $s(w_0)\mapsto\gamma(\psi_\delta)\mathcal F_\delta$, $\mathcal F_\delta\phi(y')=\int\phi(y)\psi(2\delta yy')|\delta|^{1/2}dy$;
\item $\lambda(R(t_s)f)=|s|^{1/4}\phi_\delta(y)$ if $s=\delta y^2$.
\end{enumerate}
\end{lemma}
\begin{proof}
The elements $t_\delta$ represent $\tilde G^\square\backslash\tGL_2(F)$, which gives the isomorphism.

\emph{The torus.} By \ref{C:torus}, $t_at_b=t_{ab}$ and $t_a$, $t_b$ commute. A lift of $aI_2$ commutes with $\Mp_2(F)$ (its commutator with $\tilde g$ is $\hil{\det g}a=1$ for $\det g=1$), and $\diag(a^2,1)=h(a)\cdot aI_2$; so $\rho(t_{a^2})$ is a scalar multiple of $\rho_\psi(s(h(a)))$ by Schur's lemma ($\rho^+_\psi$ is irreducible), i.e.\ $\rho(t_{a^2})=c(a)|a|^{1/2}D_a$ by \Cref{lem:oscillator}. Since $t_{a^2}t_{b^2}=t_{(ab)^2}$ and $D_aD_b=D_{ab}$, $c$ is a character. In the unramified case, by \ref{S} for $j=1$, $\ell(\theta_2(t_{a^2})\xi^\circ)=|a|^{1/2}$ for $a\in\cO\smallsetminus0$, while $\lambda(R(t_{a^2})f^\circ)=\rho(t_{a^2})\one_\cO(1)=c(a)|a|^{1/2}$. Hence $c=1$ on $\cO\smallsetminus0$, which generates $F^\times$; so $\rho(t_{a^2})=|a|^{1/2}D_a$. (This normalization is a property of the representation $\theta_2$ and does not depend on the data being unramified.)

(i) For $g\in\Mp_2(F)$, $[R(g)f](t_\delta)=\rho(t_\delta gt_\delta^{-1})f(t_\delta)$. By \ref{C:unip}, $t_\delta\tilde n(t)t_\delta^{-1}=\tilde n(\delta t)$. From the cocycle, $t_\delta s(w_0)t_\delta^{-1}=s(h(\delta)w_0)=\hil{-1}\delta s(h(\delta))s(w_0)$. By \Cref{lem:oscillator} the latter acts by $\hil{-1}\delta\mu_\psi(\delta)^{-1}|\delta|^{1/2}D_\delta\gamma(\psi)\mathcal F_\psi=\gamma(\psi_\delta)\mathcal F_\delta$, using $\hil{-1}\delta\mu_\psi(\delta)^{-1}\gamma(\psi)=\mu_\psi(\delta)\gamma(\psi)=\gamma(\psi_\delta)$.

(ii) $\lambda(R(t_s)f)=f(t_s)(1)=f(t_{y^2}t_\delta)(1)=(\rho(t_{y^2})f(t_\delta))(1)=|y|^{1/2}f(t_\delta)(y)=|y|^{1/2}|\delta|^{1/4}\phi_\delta(y)$, and $|y|^{1/2}|\delta|^{1/4}=|s|^{1/4}$.
\end{proof}

\section{The Stone--von Neumann theorem}

\begin{lemma}\label{lem:SvN}
Let $\mathcal H$ be a Heisenberg group over $F$ with center $Z\simeq F$ and nondegenerate commutator pairing on $\mathcal H/Z=\mathcal X\oplus\mathcal Y$ (a polarization), let $\chi$ be a nontrivial character of $Z$, and let $\omega$ be the Schr\"odinger representation on $\cS(\mathcal X)$ with central character $\chi$.
\begin{enumerate}[label=(\alph*)]
\item Every smooth representation $V$ of $\mathcal H$ with central character $\chi$ is isomorphic to $\cS(\mathcal X)\otimes M$ with $\mathcal H$ acting on the first factor, where $M=\mathrm{Hom}_{\mathcal H}(\cS(\mathcal X),V)$; moreover $M\simeq V_{\mathcal Y}$, the $\mathcal Y$-coinvariants, via $\cS(\mathcal X)_{\mathcal Y}=\C$ (evaluation at $0$).
\item If a group $G$ acts on $V$ and on $\mathcal H$ (fixing $\chi$) compatibly, and $\omega$ extends to a representation of $G\ltimes\mathcal H$ on $\cS(\mathcal X)$, then $G$ acts on $M$ by $\mu(g)\varphi=g\circ\varphi\circ\omega(g)^{-1}$, and $V\simeq\cS(\mathcal X)\otimes M$ as $G\ltimes\mathcal H$-modules, $G$ acting by $\omega\otimes\mu$.
\end{enumerate}
\end{lemma}
\begin{proof}
(a) The algebra $\mathcal A$ of locally constant functions on $\mathcal H$, compactly supported modulo $Z$, transforming by $\chi$ under $Z$, acts on $V$, and smooth $\mathcal H$-modules with central character $\chi$ are the nondegenerate $\mathcal A$-modules. Fourier transform in $\mathcal Y$ identifies $\mathcal A$ with the algebra of locally constant compactly supported kernels on $\mathcal X\times\mathcal X$, acting on $\cS(\mathcal X)$ by integral operators. Choose $\phi_1,\phi_2\in\cS(\mathcal X)$ with $\int\phi_1\phi_2=1$ and let $e$ be the rank one kernel $\phi_1\otimes\phi_2$. Then $e^2=e$, $e\mathcal Ae=\C e$, and $\mathcal Ae\mathcal A=\mathcal A$ (every compactly supported locally constant kernel is a finite sum of rank one kernels, each of which factors through $e$). By Morita theory for algebras with local units, $V\mapsto eV$ and $N\mapsto\mathcal Ae\otimes N$ are inverse equivalences between nondegenerate $\mathcal A$-modules and vector spaces; since $\mathcal Ae\simeq\cS(\mathcal X)$, $V\simeq\cS(\mathcal X)\otimes eV$, and $eV\simeq M$. A Schwartz function vanishing at $0$ is a finite sum of functions $(\chi(\langle T,Y\rangle)-1)f(T)$, so $\cS(\mathcal X)_{\mathcal Y}$ is the line of evaluation at $0$, and $V_{\mathcal Y}\simeq M$.
(b) $\mu(g)\varphi$ is again $\mathcal H$-equivariant because $\omega(g)$ and the action on $V$ implement the same automorphism of $\mathcal H$, and the evaluation map $\cS(\mathcal X)\otimes M\to V$ intertwines $\omega(g)\otimes\mu(g)$ with the action of $g$.
\end{proof}

In \S\ref{sec:odd}, (b) is applied with $G=\tGL_2(F)$ acting through $\iota$ and $\omega$ the linear action $f(T)\mapsto|\det g|^{(n-1)/2}f(Tg)$; under the identification $M\simeq V_{\mathcal Y}$ the action of $g$ on $V_{\mathcal Y}$ is $\mu(g)$ times the action $|\det g|^{(n-1)/2}$ of $\omega(g)$ on the evaluation line.

\backmatter

\end{document}